\documentclass[preprint,12pt]{elsarticle}

\usepackage{amssymb}
\usepackage{amsmath}
\usepackage{amsthm}

\allowdisplaybreaks[2]
\numberwithin{equation}{section}

\usepackage[a4paper,margin=2.0cm,headheight=15pt,footskip=1cm]{geometry}
\usepackage{microtype}

\biboptions{sort&compress}
\usepackage[colorlinks=true, linkcolor=blue, urlcolor=blue, citecolor=blue]{hyperref}

\usepackage{braket}

\newtheorem{Thm}{Theorem}[section]
\newtheorem{Lem}[Thm]{Lemma}
\newtheorem{Cor}[Thm]{Corollary}
\newtheorem{Prop}[Thm]{Proposition}
\newtheorem{Def}[Thm]{Definition}
\newtheorem{Rem}[Thm]{Remark}

\makeatletter
\newif\ifels@app
\def\@seccntformat#1{\ifels@app \appendixname~\fi
  \csname the#1\endcsname\@seccntDot\hskip 0.5em}
\renewcommand{\appendix}{%
  \els@apptrue
  \par
  \setcounter{section}{0}%
  \setcounter{subsection}{0}%
  \setcounter{equation}{0}%
  \gdef\thefigure{\@Alph\c@section.\arabic{figure}}%
  \gdef\thetable{\@Alph\c@section.\arabic{table}}%
  \gdef\thesection{\@Alph\c@section}%
  \@addtoreset{equation}{section}%
  \gdef\theequation{\@Alph\c@section.\arabic{equation}}%
}
\makeatother

\journal{Journal de Math\'ematiques Pures et Appliqu\'ees}

\begin{document}

\begin{frontmatter}


\title{Separated solutions for a critical competitive system with sublinear coupling in domains with shrinking holes}

\author[bnu]{Zijuan Gao}
\ead{gaozijuan3012@163.com}

\author[muc]{Qing Guo}
\ead{guoqing0117@163.com}

\author[bnu]{Chengxiang Zhang\corref{cor1}}
\ead{zcx@bnu.edu.cn}
\cortext[cor1]{Corresponding author}

\affiliation[bnu]{organization={Laboratory of Mathematics and Complex Systems (Ministry of Education), School of Mathematical Sciences, Beijing Normal University},
            city={Beijing},
            postcode={100875},
            country={P.~R.~China}}

\affiliation[muc]{organization={College of Science, Minzu University of China},
            city={Beijing},
            postcode={100081},
            country={China}}

\begin{abstract}
We construct nonnegative, fully nontrivial separated solutions for the
critical competitive system of two coupled elliptic equations with a
sublinear competitive coupling, in a bounded domain with two shrinking
holes, for $N\ge5$ and every fixed $\beta<0$; each component concentrates
at one of the holes and vanishes near the other.  The coupling is not $C^1$
at the origin, so a direct Lyapunov--Schmidt reduction is unavailable, and
the reduction has to be closed in a Lebesgue space although the matching
requires $L^\infty$ control.  We therefore solve a nonlinear outer problem
by minimization, derive the dead-core property, solve the two scalar inner
problems by a nonlinear contraction, match the inner and outer solutions,
and close the scale and position equations by a Brouwer degree argument.
This gives a Coron-type construction in the high-dimensional
sublinear-coupling regime.
\end{abstract}


\end{frontmatter}

\section{Introduction}\label{sec1}

	Competitive interactions in elliptic systems drive the components apart,
	and it is natural to ask how the geometry of the domain selects their
	concentration points.  In the spirit of Coron's problem, we consider the
	critical competitive elliptic system
	\begin{equation}\label{eq1.1}
		\begin{cases}
			-\Delta u_1 =  \mu_1 |u_1|^{2^{*}-2}u_1 + \beta |u_1|^{\frac{2^{*}}{2}-2}u_1|u_2|^{\frac{2^{*}}{2}} & \text{in } \Omega_\varepsilon, \\[3mm]
			-\Delta u_2 =  \mu_2 |u_2|^{2^{*}-2}u_2 + \beta |u_2|^{\frac{2^{*}}{2}-2}u_2|u_1|^{\frac{2^{*}}{2}} & \text{in } \Omega_\varepsilon, \\[3mm]
			u_1=u_2=0 & \text{on } \partial\Omega_\varepsilon,
		\end{cases}
	\end{equation}
	where $N\ge5$, $2^*=2N/(N-2)$, $\mu_1,\mu_2>0$, $\beta<0$, and
	\[
	\Omega_\varepsilon
	=\Omega\setminus\bigcup_{i=1}^2B_{r_i\varepsilon}(a_i),
	\]
	where $\Omega\subset\mathbb R^N$ is a smooth bounded domain,
	$a_1,a_2\in\Omega$ are distinct, and $r_1,r_2>0$.  We seek
	fully nontrivial solutions, that is, solutions with both components
	nontrivial.
	The parameter $\varepsilon>0$ tends to zero and determines the radii of the
	two holes.
	
	The scalar critical problem on a star-shaped domain has no nontrivial solution
	by the Pohozaev identity.  Coron~\cite{Coron1984} showed that a small hole can
	overcome this obstruction, and Bahri and Coron~\cite{BahriCoron1988} related
	existence to the topology of the domain.  Concentrating solutions in domains
	with shrinking holes and their refinements were subsequently studied in
	\cite{LiYanYang2004,ClappMussoPistoia2009,MussoPistoia2007,GeMussoPistoia2010}.
	The subcritical theory for competitive and mixed-coupling Schr\"odinger
	systems is well developed; see, for example,
	\cite{SatoWang2015a,SatoWang2015b,Soave2015,SoaveTavares2016,ByeonSatoWang2016}.
	For critical systems, variational and phase-separation questions in higher
	dimensions have also been investigated; see, for example,
	\cite{ChenZou2012,ChenZou2015}.
	
	Coron-type constructions for weakly coupled critical systems in domains
	with shrinking holes were obtained by Pistoia and Soave
	\cite{PistoiaSoave2018} for $N=3,4$, where the coupling is $C^1$; among
	their solutions are ones in which the components concentrate at distinct
	holes.  For the related Brezis--Nirenberg system in dimension four,
	spiked solutions concentrating at different points of a fixed domain
	were constructed in \cite{PistoiaTavares2017}.  In $N=4$, Pistoia, Soave
	and Tavares \cite{PistoiaSoaveTavares2020} constructed non-synchronized
	positive solutions in domains that are symmetric about the centre of a
	single shrinking hole: the bubbles of the different components
	concentrate at the same point, at different scales.  In contrast, the
	solutions constructed below live in the sublinear-coupling regime
	$N\ge5$, in arbitrary smooth bounded domains without any symmetry
	assumption: each component concentrates at its own hole and vanishes
	identically in a neighbourhood of the other one.
	
	For $N\ge5$ the situation changes qualitatively.  The coupling is not $C^1$
	at the origin: the factor $|u_1|^{\frac{2}{N-2}}$ is not differentiable in
	$u_1$, while in the small component $u_2$ the derivative vanishes at
	$u_2=0$ and is only H\"older continuous of exponent $\frac{2}{N-2}$.  Hence
	a direct linearization is degenerate, and the standard contraction
	scheme, which relies on a $C^2$ energy functional, does not apply; see also
	\cite[Remark~1.9]{PistoiaSoaveTavares2020}, where the $C^1$ threshold $N\le4$
	is discussed.  A separated ansatz makes this obstruction unavoidable: near
	the concentration point of one component the other component is small and
	is expected to vanish under the competitive interaction.  Consequently the
	coupling cannot be treated as an ordinary smooth perturbation in a global
	Lyapunov--Schmidt reduction, and in particular the Lyapunov--Schmidt scheme
	of \cite{PistoiaSoave2018} is not applicable in the regime $N\ge5$.
	
	A second obstruction concerns the norms in which the reduction can be closed.
	The error terms give no smallness in $L^\infty$, so a fixed point has to be
	sought in $L^q$ with $q=2N/(N-4)$; but the pointwise estimates for the outer
	minimizer require $L^\infty$ boundary data on the inner boundary.  A nonlinear
	solution operator recovers this local $L^\infty$ regularity from the inner
	equation, first by a bootstrap in $L^p$, $p>N/2$, and then by comparison
	against the explicit barriers, whereas a linear inverse supplies only the
	$L^2\to L^q$ gain of the energy theory.  Thus both the outer minimizer and
	the inner solver are nonlinear.
	
	Before stating our main result, we introduce some standard notation.
	Throughout the paper we fix a number $0<\theta<\frac{N-4}{N-2}$.  It is the loss
	exponent in all the error estimates below: the constants may depend on
	$\theta$, but never on $\varepsilon$, and since $\theta$ is arbitrary it
	may be taken as small as we wish.
	Let 
	\begin{equation}\label{eq1.2}
		U_{\delta,\xi}(x):=\alpha_N\frac{\delta^{\frac{N-2}{2}}}{(\delta^2+|x-\xi|^2)^{\frac{N-2}{2}}},\quad \alpha_N:=\left(N(N-2)\right)^{\frac{N-2}{4}}.
	\end{equation}
	It is well known that the family \(\{U_{\delta,\xi} : \delta > 0, \, \xi \in \mathbb{R}^N\}\) contains all the solutions to the critical problem
	\begin{equation}\label{eq1.3}
		\begin{cases}
			-\Delta u = u^{2^*-1}, & u > 0 \quad \text{in } \mathbb{R}^N, \\
			u \in D^{1,2}(\mathbb{R}^N),
		\end{cases}
	\end{equation}
	where \(D^{1,2}(\mathbb{R}^N)\) is the completion of \(C^\infty_c(\mathbb{R}^N)\) with respect to the norm \(\|u\|_{D^{1,2}(\mathbb{R}^N)} := \|\nabla u\|_{L^2(\mathbb{R}^N)}\).

	We shall use repeatedly the projection $P_{\varepsilon}$ of $D^{1,2}(\mathbb{R}^N)$ onto $H_0^1(\Omega_{\varepsilon})$ with respect to the Dirichlet inner product: for $v\in D^{1,2}(\mathbb{R}^N)$, the function $P_{\varepsilon} v\in H_0^1(\Omega_{\varepsilon})$ is characterized by
	\begin{equation*}
		\int_{\Omega_{\varepsilon}} \nabla (P_{\varepsilon} v) \cdot \nabla \varphi
		= \int_{\Omega_{\varepsilon}} \nabla v \cdot \nabla \varphi
		\qquad \text{for every } \varphi\in H_0^1(\Omega_{\varepsilon}).
	\end{equation*}
	Equivalently, whenever $-\Delta v$ is a function, $P_{\varepsilon} v$ is the unique solution of
	\begin{equation*}
		-\Delta(P_{\varepsilon} v)=-\Delta v \ \text{ in } \Omega_{\varepsilon},
		\qquad P_{\varepsilon} v=0 \ \text{ on } \partial \Omega_{\varepsilon}.
	\end{equation*}
	In particular, since $-\Delta U_{\delta,\xi}=U_{\delta,\xi}^{2^*-1}$ in $\mathbb{R}^N$, the projected bubble $P_{\varepsilon} U_{\delta, \xi}$ satisfies
	\begin{equation*}
		-\Delta(P_{\varepsilon} U_{\delta, \xi}) = U_{\delta, \xi}^{2^*-1} \ \text{ in } \Omega_{\varepsilon},
		\qquad P_{\varepsilon} U_{\delta, \xi} = 0 \ \text{ on } \partial \Omega_{\varepsilon}.
	\end{equation*}

	\begin{Thm}\label{thm1.1}
		Let $N\ge5$, let $\Omega\subset\mathbb R^N$ be smooth and bounded, let
		$\mu_1,\mu_2>0$, $\beta<0$, let $a_1\ne a_2$ be fixed points of
		$\Omega$, and let $\theta\in\bigl(0,\frac{N-4}{N-2}\bigr)$ be fixed.  Then there exists
		$\varepsilon_0=\varepsilon_0(\theta)>0$ such
		that, for every $0<\varepsilon<\varepsilon_0$, problem~\eqref{eq1.1}
		admits a nonnegative, fully nontrivial solution of the form
		\[
		u_{i,\varepsilon}
		=\mu_i^{-\frac1{2^*-2}}P_\varepsilon U_{\delta_{i,\varepsilon},
			\xi_{i,\varepsilon}}+\Phi_{i,\varepsilon},
		\qquad i=1,2,
		\]
		where
		\[
		\delta_{i,\varepsilon}=d_{i,\varepsilon}\sqrt\varepsilon,
		\qquad
		\xi_{i,\varepsilon}=a_i+\delta_{i,\varepsilon}\tau_{i,\varepsilon},
		\]
		and
		\[
		\tau_{i,\varepsilon}\longrightarrow0,
		\qquad
		d_{i,\varepsilon}\longrightarrow
		d_i^0:=\left(\frac{r_i^{N-2}}{H(a_i,a_i)}\right)^{\frac1{2(N-2)}}>0.
		\]
		Here $H$ is the regular part with the normalization fixed in Section~\ref{sec3}.  Moreover,
		\begin{equation}\label{estimate}
			\|\Phi_{1,\varepsilon}\|_{H_0^1(\Omega_\varepsilon)}
			+\|\Phi_{2,\varepsilon}\|_{H_0^1(\Omega_\varepsilon)}
			=O\left(\varepsilon^{\frac{N-2\theta}{4}}\right).
		\end{equation} 
		For some $c_1>0$, one has
		\[ 
		u_{1,\varepsilon}=0
		\quad\hbox{in }B_{c_1\varepsilon^{\frac{1}{N-2}}}(a_2)\cap\Omega_\varepsilon,
		\qquad
		u_{2,\varepsilon}=0
		\quad\hbox{in }B_{c_1\varepsilon^{\frac{1}{N-2}}}(a_1)\cap\Omega_\varepsilon.
		\]
	\end{Thm}

	We call the solutions of Theorem~\ref{thm1.1} {separated} because
	their two concentration points converge to the distinct points $a_1$ and
	$a_2$.  This term refers to the concentration points; it does not mean
	that the supports of the two components are globally disjoint.  The
	dead-core conclusion is a local statement near the opposite holes.  This
	terminology should not be confused with the segregated (phase-separating)
	solutions of \cite{ChenZou2012,ChenZou2015}, whose supports become
	asymptotically disjoint; here the supports overlap, and only the
	concentration points are separated.

	Compared with \cite{PistoiaSoave2018}, Theorem~\ref{thm1.1} adds two
	features: the dead core, a phenomenon typical of sublinear competitive
	interactions (see \cite{TVY}), and the absence of any symmetry assumption
	on $\Omega$.  Beyond the two obstructions discussed above, the projection
	$P_\varepsilon$ has to be controlled near two shrinking holes at the same
	time, and the inner boundary data have to be matched with the
	finite-dimensional scale and position parameters without losing
	uniformity in $\varepsilon$.

	We keep the standard Coron scale $\delta_i\sim\sqrt\varepsilon$ and
	address these difficulties in four steps.  First, a nonlinear minimization
	problem is solved outside the two inner neighbourhoods, where the
	projection $P_\varepsilon$ is controlled near both holes at once; the
	comparison estimates then yield both the required decay and, for each
	component, a dead core around the other hole.  Second, inside the inner
	neighbourhoods the coupling vanishes exactly, so the two equations
	decouple into scalar projected problems, which are solved by a nonlinear
	contraction in $L^q$.  Third, the inner and outer solutions are shown to
	coincide in their overlap regions, which matches the inner boundary data
	uniformly in $\varepsilon$.  Finally, the scale and position equations are
	read off from the reduced residuals: there the coupling contributes only
	lower-order terms, so the leading balance is the same Robin--hole balance
	as in the scalar Coron problem.

	Our construction is closer in spirit to the modified Lyapunov--Schmidt
	scheme of \cite{WZZ23} for logarithmic scalar-field equations.  In
	\cite{GGZ,G-Z}, by contrast, the solution operator is built on the inverse
	of the linearization at the ansatz: the non-smooth coupling is never
	differentiated there, since it is confined to the outer region, where it
	is handled by the minimization.  Here the ansatz error is not small enough
	in $L^p$ for $p>N/2$, so a linear inner inverse cannot produce a global
	$L^\infty$ estimate.  The energy theory gives the $L^2\to L^{2N/(N-4)}$ gain but not
	the $L^\infty$ bound, so the linear inverse is no longer adequate.
	In our construction, both operators are nonlinear.  The outer operator is a minimizer
	map, while the inner operator solves the full nonlinear scalar problem
	by a contraction in $L^{2N/(N-4)}$, the integrability gained from
	an $L^2$ right-hand side under the Dirichlet inverse.  The pointwise
	estimates for the outer minimizer nevertheless require $L^\infty$
	boundary data on the inner boundary; this regularity is recovered from
	the inner equation by the bootstrap and barrier argument described
	above, which is a further reason for the nonlinear inner solver.  The
	relevant estimates are carried out in Section~\ref{sec5}.

	\medskip
	
	The paper is organized as follows.  Section~\ref{sec2} introduces the
	rescaled problem.  Section~\ref{sec3} records the projection and
	comparison estimates.  In Section~\ref{sec4} we solve the outer problem and
	prove the dead-core property.  Section~\ref{sec5} treats the two scalar inner
	problems, and Section~\ref{sec6} identifies the inner and outer solutions in
	the overlap regions.  The finite-dimensional equations are solved in
	Section~\ref{sec7}.  The appendices collect the projection formulas and the
	parameter-continuity argument.

	\section{The rescaled problem}\label{sec2}
	By applying the scaling transformation $x = \sqrt{\varepsilon}y$ and defining the rescaled functions
	\begin{equation}\label{eq2.1}
		v_i(y) := \varepsilon^{\frac{N-2}{4}} u_i(\sqrt{\varepsilon}y), \quad i=1,2,
	\end{equation}
	the original system \eqref{eq1.1} is transformed into the equivalent problem
	\begin{equation}\label{eq2.2}
		\begin{cases}
			-\Delta v_1 = \mu_1 |v_1|^{2^*-2}v_1 + \beta |v_1|^{\frac{2^*}{2}-2}v_1 |v_2|^{\frac{2^*}{2}} & \text{in } \widetilde{\Omega}_\varepsilon, \\[3mm]
			-\Delta v_2 = \mu_2 |v_2|^{2^*-2}v_2 + \beta |v_2|^{\frac{2^*}{2}-2}v_2 |v_1|^{\frac{2^*}{2}} & \text{in } \widetilde{\Omega}_\varepsilon, \\[3mm]
			v_1 = v_2 = 0 & \text{on } \partial\widetilde{\Omega}_\varepsilon,
		\end{cases}
	\end{equation}
	where we write $\widetilde a_i$ for the rescaled hole centers, and where the
	rescaled domain $\widetilde{\Omega}_\varepsilon$ is rigorously defined as
	\begin{equation}\label{eq2.3}
		\widetilde a_i:=\frac{a_i}{\sqrt{\varepsilon}},
		\qquad
		\widetilde{\Omega}_\varepsilon := \left\{ y \in \mathbb{R}^N : \sqrt{\varepsilon}y \in \Omega_\varepsilon \right\} 
		= \frac{1}{\sqrt{\varepsilon}}\Omega \setminus \bigcup_{i=1}^2 B_{r_i\sqrt{\varepsilon}}\left(\widetilde a_i\right).
	\end{equation}
	
	Due to the conformal invariance of the critical Sobolev exponent, the form of the nonlinearities remains unchanged under this transformation.

	Correspondingly, the rescaling acts on the bubble family \eqref{eq1.2} by
	\begin{equation}\label{eq2.4}
		\varepsilon^{\frac{N-2}{4}} U_{\delta,\xi}\bigl(\sqrt{\varepsilon}\,y\bigr)
		= U_{\delta / \sqrt{\varepsilon},\; \xi / \sqrt{\varepsilon}}(y),
		\qquad \delta,\varepsilon>0,\ \xi\in\mathbb R^N .
	\end{equation}
	The free parameters of the
	construction are $(\mathbf d,\boldsymbol\tau)=(d_1,d_2,\tau_1,\tau_2)$,
	equivalently $(\delta_i,\tau_i)$ with $\delta_i=d_i\sqrt{\varepsilon}$; both
	are sets of $N+1$ independent parameters.  The ansatz prescribes the centers,
	which are therefore derived quantities,
	\[
	\xi_i:=a_i+\delta_i\tau_i,
	\qquad
	\widetilde\xi_i:=\frac{\xi_i}{\sqrt{\varepsilon}}=\widetilde a_i+d_i\tau_i ,
	\]
	$\delta_i$ entering them only through this change of variables from
	$(\delta_i,\xi_i)$ to $(\delta_i,\tau_i)$.

	In this paper, we work on the closed parameter box
	\[
	\mathcal P_{\gamma_*} := \left\{ (\mathbf{d}, \boldsymbol{\tau})=(d_1, d_2, \tau_1, \tau_{2}) \in \mathbb{R}^2 \times (\mathbb{R}^N)^2 : \gamma_* \le d_i \le \gamma_*^{-1}, \ |\tau_i| \le \gamma_* \right\},
	\]
	where $\gamma_*>0$ is a fixed small constant, and all estimates below are
	uniform on $\mathcal P_{\gamma_*}$.

	The projection $\widetilde{P}_{\varepsilon}$ onto $H_0^1(\widetilde{\Omega}_{\varepsilon})$ is defined in the same way: for $v\in D^{1,2}(\mathbb{R}^N)$, the function $\widetilde{P}_{\varepsilon} v\in H_0^1(\widetilde{\Omega}_{\varepsilon})$ is characterized by 
	\begin{equation}\label{eq2.5}
		\int_{\widetilde{\Omega}_{\varepsilon}} \nabla (\widetilde{P}_{\varepsilon} v) \cdot \nabla \varphi
		= \int_{\widetilde{\Omega}_{\varepsilon}} \nabla v \cdot \nabla \varphi
		\qquad \text{for every } \varphi\in H_0^1(\widetilde{\Omega}_{\varepsilon}).
	\end{equation}
	
	In particular, since
	$-\Delta U_{\delta,\xi}=U_{\delta,\xi}^{2^*-1}$ in $\mathbb{R}^N$, the projected
	bubbles satisfy
	\[
	-\Delta(\widetilde{P}_{\varepsilon} U_{\delta,\xi})=U_{\delta,\xi}^{2^*-1}
	\ \text{ in }\widetilde{\Omega}_{\varepsilon},
	\qquad
	\widetilde{P}_{\varepsilon} U_{\delta,\xi}=0
	\ \text{ on }\partial\widetilde{\Omega}_{\varepsilon}.
	\]

	For the parameters of the ansatz, $\delta_i=d_i\sqrt{\varepsilon}$ and
	$\xi_i=\sqrt{\varepsilon}\,\widetilde\xi_i$, this projection is related to the original projection by
	\[
	\widetilde{P}_{\varepsilon} U_{d_i,\widetilde\xi_i}
	=\varepsilon^{\frac{N-2}{4}}
	\bigl(P_\varepsilon U_{\delta_i,\xi_i}\bigr)\bigl(\sqrt{\varepsilon}\,\cdot\bigr).
	\]

	The standard inner product and norm in $H^1_0(\widetilde{\Omega}_\varepsilon)$ are given by
	\[
	\langle u, v \rangle_{H^1_0(\widetilde{\Omega}_\varepsilon)} := \int_{\widetilde{\Omega}_\varepsilon} \nabla u \cdot \nabla v \, \mathrm{d}x, \quad \| u \|_{H^1_0(\widetilde{\Omega}_\varepsilon)} := \langle u, u \rangle_{H^1_0(\widetilde{\Omega}_\varepsilon)}^{1/2}.
	\]
	For $p> 1$, we denote the usual Lebesgue norm by $\| u \|_{L^p(\widetilde{\Omega}_\varepsilon)} := \left( \int_{\widetilde{\Omega}_\varepsilon} |u|^p \, \mathrm{d}x \right)^{1/p}$.

	\medskip
	
	Our result is stated as follows.
	
	\begin{Thm}\label{thm2.1}
		Under the assumptions of Theorem~\ref{thm1.1}, for each fixed
		$\theta\in\bigl(0,\frac{N-4}{N-2}\bigr)$ there is $\varepsilon_0(\theta)>0$ such that for every $0<\varepsilon<\varepsilon_0(\theta)$,   the rescaled
		problem~\eqref{eq2.2} admits a solution
		\[
		v_{i,\varepsilon}
		=\mu_i^{-\frac1{2^*-2}}
		\widetilde P_\varepsilon
		U_{d_{i,\varepsilon},\widetilde\xi_{i,\varepsilon}}
		+\phi_{i,\varepsilon},
		\qquad
		\widetilde\xi_{i,\varepsilon}
		=\widetilde a_i+d_{i,\varepsilon}\tau_{i,\varepsilon},
		\]
		where $d_{i,\varepsilon}\to d_i^0$ and $\tau_{i,\varepsilon}\to0$, and
		\[
		\|\phi_{1,\varepsilon}\|_{H_0^1(\widetilde\Omega_\varepsilon)}
		+\|\phi_{2,\varepsilon}\|_{H_0^1(\widetilde\Omega_\varepsilon)}
		=O\left(\varepsilon^{\frac{N-2\theta}{4}}\right),
		\]
		the constant in $O(\cdot)$ depending on $\theta$ but not on $\varepsilon$.
		The two components are nonnegative and fully nontrivial, and
		\[
		v_{1,\varepsilon}=0\quad\hbox{in }
		B_{c_1\varepsilon^{-\frac{N-4}{2(N-2)}}}(\widetilde a_2)\cap\widetilde\Omega_\varepsilon,
		\qquad
		v_{2,\varepsilon}=0\quad\hbox{in }
		B_{c_1\varepsilon^{-\frac{N-4}{2(N-2)}}}(\widetilde a_1)\cap\widetilde\Omega_\varepsilon.
		\]
	\end{Thm}
	In the original variables the correction $\phi_{i,\varepsilon}$ corresponds
	to $\Phi_{i,\varepsilon}(x)=\varepsilon^{-\frac{N-2}{4}}
	\phi_{i,\varepsilon}(x/\sqrt\varepsilon)$ in Theorem~\ref{thm1.1}.

	\section{Preliminaries}\label{sec3}
	
	\subsection{Notation and preliminary results}\label{subsec3.1}

	Recall that $U_{d_i, \tilde{\xi}_i}$ denotes the \textbf{standard bubble} in the rescaled variables, which is the family \eqref{eq1.2} with the rescaled parameters, and that $\widetilde{P}_{\varepsilon} U_{d_i, \tilde{\xi}_i}$ denotes its projection into $H_0^1(\widetilde{\Omega}_{\varepsilon})$, defined as in \eqref{eq2.5}. We shall frequently use the fact that $0 \le \widetilde{P}_{\varepsilon} U_{d_i, \tilde{\xi}_i} \le U_{d_i, \tilde{\xi}_i}$ in $\widetilde{\Omega}_{\varepsilon}$, which is a direct consequence of the maximum principle;
	see, e.g., \cite{PP07}.
	
	Since $\{U_{\delta, \xi}\}$ constitutes the set of all positive solutions to the limit problem \eqref{eq1.3}, it is straightforward to verify that any solution to
	\[
	\begin{cases} 
		-\Delta v = \mu v^{2^*-1}, & v > 0 \quad \text{in } \mathbb{R}^N, \\ 
		v \in D^{1,2}(\mathbb{R}^N),
	\end{cases}
	\]
	with $\mu > 0$ is given by $\mu^{-\frac{1}{2^*-2}} U_{d, \tilde{\xi}}$ for some $d > 0$ and $\tilde{\xi} \in \mathbb{R}^N$.
	
	Returning to the linearized equation around the standard bubble, we recall some classical properties. The linearized operator
	\[
	L_{d,\tilde{\xi}}:=-\Delta-(2^*-1)\,U_{d, \tilde{\xi}}^{2^*-2},
	\]
	acting on $D^{1,2}(\mathbb{R}^N)$, is non-degenerate (see, e.g., \cite{CaffarelliGidasSpruck1989}): its kernel consists exactly of the solutions of the linearized equation
	\begin{equation}\label{eq3.1}
		-\Delta \phi = (2^*-1)\, U_{d, \tilde{\xi}}^{2^*-2} \phi \quad \text{in } \mathbb{R}^N,
		\qquad \phi\in D^{1,2}(\mathbb{R}^N),
	\end{equation}
	and is $(N+1)$-dimensional, spanned by the derivatives of $U_{d, \tilde{\xi}}$ with respect to its parameters,
	\begin{equation}\label{eq3.2}
		\psi_{d, \tilde{\xi}}^0 := \frac{\partial U_{d, \tilde{\xi}}}{\partial d}, \quad 
		\psi_{d, \tilde{\xi}}^h := \frac{\partial U_{d, \tilde{\xi}}}{\partial \tilde{\xi}_h} \quad (h=1,\dots,N).
	\end{equation}
	Explicitly, these basis functions are given by
	\[
	\psi_{d, \tilde{\xi}}^0(y) = \alpha_N \frac{N-2}{2} d^{\frac{N-4}{2}} \frac{|y-\tilde{\xi}|^2 - d^2}{(d^2 + |y-\tilde{\xi}|^2)^{N/2}}, \quad
	\psi_{d, \tilde{\xi}}^h(y) = \alpha_N (N-2) d^{\frac{N-2}{2}} \frac{y_h - \tilde{\xi}_h}{(d^2 + |y-\tilde{\xi}|^2)^{N/2}}.
	\]
	Since $\widetilde{\Omega}_{\varepsilon}$ and the homogeneous Dirichlet condition
	do not depend on the parameters $(d_i,\widetilde\xi_i)$, the Laplacian and the
	boundary condition commute with $\partial_h$, and we have the fundamental identity
	\[
	\widetilde{P}_{\varepsilon} \psi_{d_i, \tilde{\xi_i}}^h = \partial_h (\widetilde{P}_{\varepsilon} U_{d_i, \tilde{\xi_i}}), \quad h=0,\dots,N,
	\]
	where $\partial_0 = \partial_{d_i}$ and $\partial_h = \partial_{\tilde{\xi_i}_h}$ for $h \ge 1$.
	Throughout the paper, the scale derivative $\partial_{d_i}$ is taken at fixed
	center $\widetilde\xi_i$, whereas the translation derivative
	$\partial_{\widetilde\xi_{i,h}}$ is taken at fixed scale.

	\medskip

	Choose the fixed compact parameter range to contain the points
	$(d_i^0,0)$ in Theorem~\ref{thm1.1}, and work on a slightly larger compact
	range when taking limits.
	Choose a radial function
	$\chi_0\in C_c^\infty(B_2(0))$ such that $0\le\chi_0\le1$ and
	$\chi_0=1$ in $B_1(0)$.  For $i=1,2$, set
	\[
	Z_i^h(y):=\chi_0(y-\widetilde\xi_i)\,
	\widetilde P_\varepsilon
	\psi_{d_i,\widetilde\xi_i}^h,
	\qquad h=0,\ldots,N.
	\]
	The functions $Z_i^h$, $h=0,\ldots,N$, are the truncations of the derivatives
	$\partial_h\widetilde P_\varepsilon U_{d_i,\widetilde\xi_i}$ of the $i$-th bubble; we write $Z^h:=Z_1^h$ and $Y^h:=Z_2^h$, and define
	\[
	K_i:=\operatorname{span}\{Z_i^0,\ldots,Z_i^N\},
	\qquad
	K_i^\perp:=\left\{\phi\in H_0^1(\widetilde\Omega_\varepsilon):
	\int_{\widetilde\Omega_\varepsilon}\nabla\phi\cdot\nabla Z_i^h=0,
	\ h=0,\ldots,N\right\}.
	\]
	Finally, put
	\[
	\mathbf K_{\mathbf d,\boldsymbol\tau,\varepsilon}^\perp
	:=K_1^\perp\times K_2^\perp.
	\]
	The finite matrices needed below are checked at the point where they are
	used.  Their limits retain the fixed cutoff $\chi_0$; it does not disappear
	when $\varepsilon\to0$.
	
	Finally, let $G(x,z)$ be the Dirichlet Green function of $-\Delta$ in $\Omega$.
	\[
	G(x,z)=\frac{|x-z|^{2-N}-H(x,z)}{N(N-2)\omega_N}.
	\]
	where $\omega_N$ is the volume of the unit ball in $\mathbb{R}^N$, and the normalization of the regular part $H(\cdot,z)$ is fixed by
	\[
	\begin{cases}
		-\Delta_x H(x,z) = 0 & \text{in } \Omega, \\
		H(x,z) = |x-z|^{2-N} & \text{on } \partial\Omega.
	\end{cases}
	\]

	Let $\mathcal S>0$ be the optimal Sobolev constant in $\mathbb R^N$, normalized by
	\[
	\mathcal S\|w\|_{L^{2^*}(\mathbb R^N)}^2
	\le \|\nabla w\|_{L^2(\mathbb R^N)}^2,
	\qquad w\in D^{1,2}(\mathbb R^N).
	\]
	By zero extension, the same inequality holds with the same constant for every
	$w\in H_0^1(\Omega_0)$ and every open set $\Omega_0\subset\mathbb R^N$; see, e.g., \cite{Adams1975}.
	To obtain pointwise decay estimates, we use the following comparison principle,
	whose proof can be found in \cite[Lemma 2.5]{GGZ}.
	\begin{Lem}\label{lem3.1} 
		Let $\Omega_0 \subset \widetilde{\Omega}_{\varepsilon}$ be a bounded open set with Lipschitz boundary. 
		Suppose $\phi \in L^{N/2}(\Omega_0)$ satisfies $\|\phi\|_{L^{N/2}(\Omega_0)} < \mathcal{S}$.
		Let $v_1, v_2 \in H^1(\Omega_0)$ with $v_1 \geq 0$ a.e.\ in $\Omega_0$. If  
		\[
		-\Delta (v_2 - v_1) - \phi(v_2 - v_1) \leq 0 \quad \text{in } \Omega_0
		\]  
		in the weak sense, 
		and $(v_2-v_1)|_{\partial\Omega_0}\le0$ in the sense of Sobolev traces,
		then $v_2 \leq v_1$ a.e.\ in $\Omega_0$.  
	\end{Lem}

	
	\begin{Lem}\label{lem3.2}
		For $i=1,2$,  
		there exists $C>0$, independent of $\varepsilon$ and uniform for
		$(\mathbf d,\boldsymbol\tau)\in \mathcal P_{\gamma_*}$, such that
		
		\begin{equation}\label{eq3.3}
			\begin{aligned}
				U_{d_i,\tilde{\xi}_i}^{2^*-1} - \bigl(\widetilde{P}_{\varepsilon} U_{d_i,\tilde{\xi}_i}\bigr)^{2^*-1}
				\leq &C\varepsilon^{\frac{N-2}{2}}\left[ \frac{1}{(d_i + |y - \widetilde \xi_{i}|)^4} + \frac{1}{|y - \widetilde a_i|^{N-2}(d_i + |y - \widetilde \xi_{i}|)^4} \right]
				.
			\end{aligned}
		\end{equation}
	\end{Lem}
	
	\begin{proof}
		By the maximum principle, $0 \le \widetilde{P}_{\varepsilon} U_{d_i,\tilde{\xi}_i} \le U_{d_i,\tilde{\xi}_i}$ in $\widetilde{\Omega}_{\varepsilon}$.
		By the Mean Value Theorem, we have 
		\[
		0 \le U_{d_i,\tilde{\xi}_i}^{2^*-1} - \bigl(\widetilde{P}_{\varepsilon} U_{d_i,\tilde{\xi}_i}\bigr)^{2^*-1}
		\le C \, U_{d_i,\tilde{\xi}_i}^{2^*-2}(y)\,
		\bigl|U_{d_i,\tilde{\xi}_i}(y) - \widetilde{P}_{\varepsilon} U_{d_i,\tilde{\xi}_i}(y)\bigr|.
		\]
		By Corollary~\ref{corB.2}, 
		we immediately obtain \eqref{eq3.3}.
	\end{proof}

	\section{The outer problem}\label{sec4}
	For $i=1,2$, let
	\begin{equation*}
		v_i^* := \mu_i^{-\frac{1}{2^*-2}}\widetilde{P}_{\varepsilon} U_{d_i,\widetilde\xi_i} .
	\end{equation*}
	Let $R_0>10$ be sufficiently large, independently of $\varepsilon$.
	Throughout this section, $R_0$ may be increased finitely many times, always
	independently of $\varepsilon$.  The fixed number $\theta\in\bigl(0,\frac{N-4}{N-2}\bigr)$ of
	Section~\ref{sec1} will also serve, through the exponent $N-2-\theta$, to keep the
	auxiliary barrier exponents strictly below the critical one.  We define
	\begin{equation*}
		M_{0}=  B_{R_0}(\widetilde a_1)\setminus B_{r_1\sqrt{\varepsilon}}\left(\widetilde a_1\right),~Q_{0}=B_{R_0}(\widetilde a_2)\setminus B_{r_2\sqrt{\varepsilon}}\left(\widetilde a_2\right).
	\end{equation*}
	Since $|\widetilde\xi_i-\widetilde a_i|=d_i|\tau_i|\le1$ and $3<R_0$, and since $Z^h,Y^h$ of Section~\ref{sec3} vanish on $\partial\widetilde\Omega_\varepsilon$,
	\[
	\operatorname{supp}Z^h\subset B_2(\widetilde\xi_1)\subset B_3(\widetilde a_1)\subset B_{R_0}(\widetilde a_1),
	\qquad
	\operatorname{supp}Y^h\subset B_2(\widetilde\xi_2)\subset B_3(\widetilde a_2)\subset B_{R_0}(\widetilde a_2);
	\]
	as $M_0$ and $Q_0$ are these balls with the small holes removed, the supports are contained in $M_0$ and $Q_0$ respectively.
	
	All the comparison arguments in this section rest on the following standard
	barriers, together with an elementary inequality for two centers that is
	displayed in the proof of Lemma~\ref{lem4.3}.  The barriers are stated in the
	form in which they are used, so that Lemmas~\ref{lem4.2}, \ref{lem4.3} and
	\ref{lem4.8} can quote them without recomputing the Laplacian.
	
	\begin{Lem}\label{lem4.1}
		Let $R_0\ge2\gamma_*^{-1}$, $z\in\mathbb R^N$ and
		$d\in[\gamma_*,\gamma_*^{-1}]$, and put $r(y):=|y-z|$.  For $\alpha>0$ set
		\[
		\varpi_\alpha:=d^{\frac{\alpha}{2}}(d+r)^{-\alpha}.
		\]
		Then there exist constants $c_0,C>0$, depending only on $N$ and $\gamma_*$,
		such that on $\{r\ge R_0/2\}$
		\begin{gather}
			-\Delta \varpi_{N-2}\ge c_0R_0\, \varpi_{N-2}^{2^*-1},
			\label{eq4.1}\\
			-\Delta \varpi_\alpha\ge c_0\alpha R_0\,\varpi_4\varpi_\alpha
			+\alpha(N-2-\alpha)\,d^{-1}\varpi_{\alpha+2},
			\qquad 0<\alpha \leq  N-2,
			\label{eq4.2}\\
			\|\varpi_{N-2}\|_{L^{2^*}(r\ge R_0/2)}^{2^*-2}\le CR_0^{-2},
			\qquad
			\|\varpi_4\|_{L^{N/2}(r\ge R_0/2)}\le CR_0^{-2}.
			\label{eq4.3}
		\end{gather}
	\end{Lem}
	
	\begin{proof}
		We have  
		\begin{equation}\label{eq4.4}
			\begin{aligned}
				-\Delta \varpi_\alpha
				&=\alpha(N-2-\alpha)d^{\frac\alpha2}(d+r)^{-\alpha-2}
				+\alpha(N-1)d^{\frac\alpha2}\frac{d}{r}\,(d+r)^{-\alpha-2}\\
				&=\alpha(N-2-\alpha)d^{-1}\varpi_{\alpha+2}
				+\alpha(N-1)d^{-1}\frac{(d+r)^2}{r}\,\varpi_4\varpi_\alpha.
			\end{aligned}
		\end{equation}
		Since $r\ge R_0/2\ge d$,
		\[
		d^{-1}\frac{(d+r)^2}{r}\ge\frac{r}{d}\ge\frac{R_0}{2d}
		\ge\frac{\gamma_*}{2}R_0 ,
		\]
		so that \eqref{eq4.4} gives \eqref{eq4.2} with
		$c_0:=\frac{N-1}{2}\gamma_*$ and $0<\alpha\le N-2$.  At the endpoint
		$\alpha=N-2$ the coefficient $\alpha(N-2-\alpha)$ vanishes and
		$\varpi_4\varpi_{N-2}=\varpi_{N-2}^{2^*-1}$, so the same inequality is \eqref{eq4.1}.
		
		For the norms, $\varpi_{N-2}^{2^*}=\varpi_4^{N/2}=d^N(d+r)^{-2N}$ and
		$r\ge R_0/2$, so
		\[
		\|\varpi_{N-2}\|_{L^{2^*}}^{2^*}=\|\varpi_4\|_{L^{N/2}}^{N/2}
		=N\omega_N d^N\int_{R_0/2}^{\infty}\frac{s^{N-1}}{(d+s)^{2N}}\,ds
		=N\omega_N\int_{R_0/(2d)}^{\infty}\frac{t^{N-1}}{(1+t)^{2N}}\,dt,
		\]
		and the last integral is at most $C(2d/R_0)^N\le CR_0^{-N}$, which is
		\eqref{eq4.3}.  Note also that
		$\varpi_{N-2}=\varpi_4^{(N-2)/4}$, which identifies the two barriers of the
		same center.
	\end{proof}
	
	Throughout this section $\varpi^{(i)}_\alpha$ denotes the barrier
	$\varpi_\alpha$ of Lemma~\ref{lem4.1} at the center $\widetilde\xi_i$
	and scale $d_i$:
	\[
	\varpi^{(i)}_\alpha
	=d_i^{\frac{\alpha}{2}}(d_i+|y-\widetilde\xi_i|)^{-\alpha},
	\qquad \alpha>0 .
	\]

	After all the requirements on $R_0$ below have been imposed, let
	$R>4R_0$ be fixed, sufficiently large.  For all sufficiently small
	$\varepsilon$, the balls $B_{5R}(\widetilde a_i)$ are disjoint and lie
	in the rescaled outer domain.
	
	Denote
	\begin{equation*}
		M_i:=B_{iR}(\widetilde a_1)\setminus
		B_{r_1\sqrt\varepsilon}\left(\widetilde a_1 \right),
		\qquad
		Q_i:=B_{iR}(\widetilde a_2)\setminus
		B_{r_2\sqrt\varepsilon}\left(\widetilde a_2 \right),
		\qquad i=1,2,3.
	\end{equation*}
	Then $M_0\subset M_i$ and $Q_0\subset Q_i$.  The sets $M_0,Q_0$
	are the prescribed-data regions for the outer problem, whereas $M_3,Q_3$
	are the domains of the inner projected problems.  The collars
	$M_3\setminus M_0$ and $Q_3\setminus Q_0$ are the overlap regions used in
	the gluing argument.
	Set
	\[
	\Lambda_{\varepsilon,1}:=
	\left\{\phi_1\in K_1^\perp:
	\|\nabla\phi_1\|_{L^2(M_0)}\le
	\varepsilon^{\frac{N-2\theta}{4}},
	\|\phi_1\|_{L^\infty(M_0)}\le
	\varepsilon^{\frac{N-2\theta}{4}}
	\right\},
	\]
	where the two norms are computed on the prescribed region $M_0$ only.  Since
	$\phi_1$ vanishes on $\partial B_{r_1\sqrt\varepsilon}(\widetilde a_1)$, the
	Poincar\'e inequality on the fixed region $M_0$ bounds
	$\|\phi_1\|_{L^2(M_0)}$ by $C\|\nabla\phi_1\|_{L^2(M_0)}$, so replacing the
	gradient norm by the full $H^1$ norm changes the class only by a constant.
	Thus
	$\phi_1$ is a global $H_0^1(\widetilde\Omega_\varepsilon)$ function whose
	restriction to $M_0$ is the datum entering \eqref{eq4.6}; since
	$H^1(M_0)\subset L^{2^*}(M_0)$, its trace on $\partial M_0$ is well defined
	and inherits the $L^\infty$ bound.
	Define $\Lambda_{\varepsilon,2}$ analogously, with $Q_0$ in place of
	$M_0$.  Finally, let
	\[
	\Lambda_\varepsilon:=\Lambda_{\varepsilon,1}\times
	\Lambda_{\varepsilon,2}.
	\]
	(The set depends on the parameters through $M_0$ and $Q_0$; the argument is
	suppressed whenever the parameters are fixed, as in Appendix~\ref{secC}.)
	
	With the domain decomposition and the prescribed-data space $\Lambda_\varepsilon$ in place, we now prescribe the data and formulate the outer problem.
	
	For $(\phi_{1,0},\phi_{2,0})\in\Lambda_\varepsilon$ fixed,
	this section aims to solve, on the whole rescaled domain, the system
	\begin{equation}\label{eq4.5}
		\begin{cases}
			- \Delta v_1  = \mu_1 |v_1|^{2^*- 2}v_1 + \beta |v_1|^{\frac{2^*}{2}-2}v_1 |v_2|^{\frac{2^*}{2}}, &~\text{in}~\widetilde\Omega_{\varepsilon} \setminus M_0,\\[3mm]
			-\Delta v_2 = \mu_2 |v_2|^{2^*- 2}v_2 + \beta |v_2|^{\frac{2^*}{2}-2}v_2 |v_1|^{\frac{2^*}{2}}, &~\text{in}~\widetilde\Omega_{\varepsilon} \setminus Q_0,
		\end{cases}
	\end{equation}
	with 
	the prescribed data 
	\begin{equation}\label{eq4.6}
		v_{1}=\phi_{1,0}+\mu_1^{-\frac{1}{2^*-2}}\widetilde{P}_{\varepsilon} U_{d_1,\tilde{\xi}_1} ~\text{in}~ M_0,~~v_{2}=\phi_{2,0}+\mu_2^{-\frac{1}{2^*-2}}\widetilde{P}_{\varepsilon} U_{d_2,\tilde{\xi}_2}~\text{in}~ Q_0,
	\end{equation}
	and we shall look for solutions satisfying the smallness condition
	\begin{equation}\label{eq4.7}
		\begin{aligned}
			(2^{*}-1) \mu_1 \lVert v_{1} \rVert^{2^*-2}_{L^{2^*}(\widetilde\Omega_{\varepsilon}\setminus M_0)}<\frac{\mathcal S}{2}~~\text{and}~~(2^{*}-1) \mu_2 \lVert v_{2} \rVert^{2^*-2}_{L^{2^*}(\widetilde\Omega_{\varepsilon}\setminus Q_0)}<\frac{\mathcal S}{2}.
		\end{aligned}
	\end{equation}
	
	The estimates for this problem do not seem easy to obtain, so
	we do not
	attack \eqref{eq4.5} directly.  We first solve and give estimates for the intermediate problem
	\begin{equation}\label{eq4.8}
		\begin{cases}
			- \Delta v_1  = \mu_1 |v_1|^{2^*- 2}v_1 + \beta |v_1|^{\frac{2^*}{2}-2}v_1 |v_2|^{\frac{2^*}{2}}, &~\text{in}~\widetilde\Omega_{\varepsilon} \setminus (M_0\cup Q_0),\\ 
			-\Delta v_2 = \mu_2 |v_2|^{2^*- 2}v_2 + \beta |v_2|^{\frac{2^*}{2}-2}v_2 |v_1|^{\frac{2^*}{2}}, &~\text{in}~\widetilde\Omega_{\varepsilon} \setminus (M_0\cup Q_0),
		\end{cases}
	\end{equation}
	with the boundary condition \eqref{eq4.6}, the smallness condition
	\eqref{eq4.7}, and with
	\begin{equation}\label{eq4.9}
		v_1=0\quad \text{in}\quad Q_0,\quad v_2=0 \quad \text{in} \quad M_0.
	\end{equation}
	Each equation of \eqref{eq4.8} keeps the component's own hole and removes
	the artificial region around the other center; throughout, $v_1$ is
	extended by zero in $Q_0$ and $v_2$ by zero in $M_0$, so that both are
	defined on all of $\widetilde\Omega_\varepsilon$.  By the dead-core
	Lemma~\ref{lem4.6}, any solution of \eqref{eq4.8} satisfying \eqref{eq4.6},
	\eqref{eq4.9} and \eqref{eq4.7} vanishes in $Q_3$ and in $M_3$,
	respectively, hence on $Q_0\subset Q_3$ and $M_0\subset M_3$, where the
	equations of \eqref{eq4.5} then hold trivially.  Therefore every solution
	of \eqref{eq4.8} solves \eqref{eq4.5} with the same data, and it suffices
	to solve \eqref{eq4.8}.

	\subsection{Pointwise estimates}
	By Kato's inequality, any solution $(v_1,v_2)$ satisfying \eqref{eq4.8} also satisfies 
	\begin{equation}\label{eq4.10}
		\begin{cases}
			-\Delta \lvert v_1\rvert \leq \mu_1\lvert v_1\rvert^{2^{*}-1},  &~\text{in}~\widetilde \Omega_{\varepsilon} \setminus (M_0\cup Q_0),\\[3mm]
			-\Delta \lvert v_2 \rvert\leq \mu_2\lvert v_2\rvert^{2^{*}-1},  &~\text{in}~\widetilde\Omega_{\varepsilon} \setminus (M_0\cup Q_0).
		\end{cases}
	\end{equation}
	
	\begin{Lem}\label{lem4.2}
		For $R_0$ sufficiently large, let
		$(\phi_{1,0},\phi_{2,0})\in\Lambda_\varepsilon$, and let $(v_1,v_2)$ satisfy
		\eqref{eq4.6}, \eqref{eq4.7}, \eqref{eq4.9} and \eqref{eq4.10}.
		Then  
		\[
		|v_1(y)| \leq A_1  \frac{d_{1}^{\frac{N-2}{2}}}{(d_{1} + |y - \widetilde \xi_{1}|)^{N-2}},
		\quad
		|v_2(y)| \leq A_2   \frac{d_{2}^{\frac{N-2}{2}}}{(d_{2} + |y - \widetilde \xi_{2}|)^{N-2}},
		\]
		for $y \in \widetilde\Omega_{\varepsilon} \setminus (M_0\cup Q_0)
		$, where 
		\begin{equation*}
			A_1= 2\left( \frac{4N(N-2)}{\mu_1} \right)^{\frac{N-2}{4}},\quad A_2= 2\left( \frac{4N(N-2)}{\mu_2} \right)^{\frac{N-2}{4}}.
		\end{equation*}
	\end{Lem}
	
	\begin{proof}
		On $\widetilde\Omega_\varepsilon\setminus M_0$ one has
		$|y-\widetilde\xi_1|\ge R_0-1\ge R_0/2$.  The function $\varpi_{N-2}$
		of Lemma~\ref{lem4.1} with $z=\widetilde\xi_1$, $d=d_1$ is
		precisely the factor multiplying $A_1$ on the right-hand side of the
		statement, so \eqref{eq4.1} gives
		\[
		-\Delta \varpi_{N-2}\ge c_0R_0\,\varpi_{N-2}^{2^*-1}
		\qquad\hbox{in }\widetilde\Omega_\varepsilon\setminus M_0.
		\]

		For $(\phi_{1,0},\phi_{2,0})\in\Lambda_\varepsilon$ and $y\in\partial M_0$,
		\begin{equation}\label{eq4.11}
			\begin{aligned}
				\lvert v_1(y)\rvert=&\lvert \phi_{1,0}+\mu_1^{-\frac{1}{2^*-2}}\widetilde{P}_{\varepsilon} U_{d_1,\tilde{\xi}_1}\rvert \leq \lvert \phi_{1,0}\rvert+\mu_1^{-\frac{1}{2^*-2}}U_{d_1,\tilde{\xi}_1}\\
				=&\lvert \phi_{1,0}\rvert+\mu_1^{-\frac{1}{2^*-2}}\big(N(N-2)\big)^{\frac{N-2}{4}}\frac{d_1^{\frac{N-2}{2}}}{\left(d_1^2 + |y - \tilde{\xi}_1|^2\right)^{\frac{N-2}{2}}}
				\\
				\leq& A_1 \varpi_{N-2}(y).
			\end{aligned}
		\end{equation}
		In \eqref{eq4.11}, the second term is controlled using
		$\sqrt2(d_1^2+|y-\widetilde\xi_1|^2)^{1/2}\ge d_1+|y-\widetilde\xi_1|$; for the first
		term, $|\phi_{1,0}|\le\|\phi_{1,0}\|_{L^\infty(M_0)}
		\le\varepsilon^{(N-2\theta)/4}$ while
		$d_1^{(N-2)/2}(d_1+R_0)^{-(N-2)}>0$,
		so that the inequality holds for $\varepsilon$ sufficiently small.

		The preceding computation gives
		\[
		-\Delta(A_1\varpi_{N-2})
		\ge \frac{c_0R_0}{A_1^{2^*-2}}
		(A_1\varpi_{N-2})^{2^*-1}
		\quad\hbox{in }\widetilde\Omega_{\varepsilon}\setminus M_0.
		\]
		By
		\eqref{eq4.3}, choose $R_0$ so large that
		\[
		c_0R_0\ge \mu_1A_1^{2^*-2}
		\quad\hbox{and}\quad
		(2^*-1)\mu_1\|A_1\varpi_{N-2}\|_{L^{2^*}(\widetilde\Omega_{\varepsilon}\setminus M_0)}^{2^*-2}
		<\frac{\mathcal S}{2}.
		\] 
		On the boundary of $\widetilde\Omega_\varepsilon\setminus(M_0\cup Q_0)$ there holds
		$|v_1|\le A_1 \varpi_{N-2}$: on $\partial M_0$ this is \eqref{eq4.11}; on
		$\partial Q_0\cup \partial\widetilde\Omega_\varepsilon$ one has $v_1=0$ by \eqref{eq4.9}.  
		
		Kato's inequality, $\beta<0$, and the supersolution inequality imply
		\[
		-\Delta(|v_1|-A_1\varpi_{N-2})
		\le \widetilde\varphi_1(|v_1|-A_1\varpi_{N-2}),
		\]
		where
		\[
		\widetilde\varphi_1(y)
		:=\mu_1(2^*-1)\int_0^1
		\bigl[t|v_1(y)|+(1-t)A_1\varpi_{N-2}(y)\bigr]^{2^*-2}\,dt.
		\] 
		By  
		\eqref{eq4.7},
		\[
		\|\widetilde\varphi_1\|_{L^{N/2}(\widetilde\Omega_{\varepsilon}\setminus M_0)}
		\le (2^*-1)\mu_1
		\left(\|v_1\|_{L^{2^*}}^{2^*-2}
		+\|A_1\varpi_{N-2}\|_{L^{2^*}}^{2^*-2}\right)
		<\mathcal S.
		\] 
		
		Since $v_1=0$ in $Q_0$, inequality \eqref{eq4.10} extends trivially
		from $\widetilde\Omega_\varepsilon\setminus(M_0\cup Q_0)$ to
		$\widetilde\Omega_\varepsilon\setminus M_0$ (both sides vanish in
		$Q_0$).  Lemma~\ref{lem3.1}, applied on
		$\widetilde\Omega_\varepsilon\setminus(M_0\cup Q_0)$ with
		$v_1=A_1\varpi_{N-2}\ge0$ and $v_2=|v_1|$, now gives the desired
		estimate for $v_1$.  The estimate for $v_2$ is obtained similarly,
		with $\widetilde\Omega_\varepsilon\setminus(M_0\cup Q_0)$ and
		replacing $(d_1,\widetilde\xi_1)$ by $(d_2,\widetilde\xi_2)$ throughout.
	\end{proof}
	
	\begin{Lem}\label{lem4.3}
		Let $(\phi_{1,0},\phi_{2,0})\in\Lambda_\varepsilon$, and let $(v_1,v_2)$ satisfy
		\eqref{eq4.6}, \eqref{eq4.7}, \eqref{eq4.9} and \eqref{eq4.8}.  Upon
		increasing $R_0$ if necessary, still independently of $\varepsilon$, one has
		\begin{align}
			|v_1-v_1^*|
			&\le C R_0^{N/2}\varepsilon^{\frac{N-2\theta}{4}}(d_1+|y-\widetilde\xi_1|)^{-N/2}
			+C\varepsilon(d_2+|y-\widetilde\xi_2|)^{-(N-2-\theta)}, 
			\label{eq4.12}\\
			|v_2-v_2^*|
			&\le C R_0^{N/2}\varepsilon^{\frac{N-2\theta}{4}}(d_2+|y-\widetilde\xi_2|)^{-N/2}
			+C\varepsilon(d_1+|y-\widetilde\xi_1|)^{-(N-2-\theta)},
			\label{eq4.13}
		\end{align}
		for $y \in \widetilde\Omega_{\varepsilon} \setminus (M_0\cup Q_0)$, where
		$C>0$ depends only on $N,\gamma_*,\mu_1,\mu_2,\beta$, $R_0$ and $\theta$, and
		is uniform for $(\mathbf d,\boldsymbol\tau)\in \mathcal P_{\gamma_*}$ and
		$(\phi_{1,0},\phi_{2,0})\in\Lambda_\varepsilon$.
	\end{Lem}
	
	\begin{proof}
		Let $(\phi_1, \phi_2) = (v_1, v_2) - (v_1^*, v_2^*)$.  We use the
		notation of Lemma~\ref{lem4.1}, so that $\varpi^{(i)}_{N-2}$ is
		the barrier of Lemma~\ref{lem4.2} and
		$\varpi^{(i)}_4=d_i^2(d_i+|y-\widetilde\xi_i|)^{-4}$.  Since $d_i\in[\gamma_*,\gamma_*^{-1}]$, the
		factors $d_i^{\pm}$ are bounded and are absorbed into the constants
		without further mention.  Recall that
		$-\Delta\widetilde P_\varepsilon U_{d_i,\widetilde\xi_i}
		=U_{d_i,\widetilde\xi_i}^{2^*-1}$ in
		$\widetilde\Omega_\varepsilon$, by \eqref{eq2.5} and
		$-\Delta U_{d_i,\widetilde\xi_i}=U_{d_i,\widetilde\xi_i}^{2^*-1}$ in
		$\mathbb R^N$.  Hence the reference profiles satisfy
		\[
		-\Delta v_1^*=\mu_1^{-\frac{1}{2^*-2}}
		U_{d_1,\widetilde\xi_1}^{2^*-1},
		\qquad
		-\Delta v_2^*=\mu_2^{-\frac{1}{2^*-2}}
		U_{d_2,\widetilde\xi_2}^{2^*-1}
		\quad\hbox{in }\widetilde\Omega_\varepsilon,
		\]
		with $v_i^*=0$ on $\partial\widetilde\Omega_\varepsilon$.
		Subtracting these identities from \eqref{eq4.8} and using
		\eqref{eq4.6}, \eqref{eq4.9}, we have
		\begin{equation}\label{eq4.14}
			\begin{cases}
				-\Delta \phi_1 = \mu_1 |v_1|^{2^* - 2}v_1 - \mu_1^{-\frac{1}{2^*-2}} U_{d_1,\tilde{\xi}_1}^{2^*-1} + \beta |v_1|^{\frac{2^*}{2}-2}v_1|v_2|^{\frac{2^*}{2}}, & \text{in } \widetilde{\Omega}_{\varepsilon} \setminus (M_0\cup Q_0),\\
				-\Delta \phi_2 = \mu_2 |v_2|^{2^* - 2}v_2 - \mu_2^{-\frac{1}{2^*-2}} U_{d_2,\tilde{\xi}_2}^{2^*-1} + \beta |v_2|^{\frac{2^*}{2}-2}v_2|v_1|^{\frac{2^*}{2}}, & \text{in } \widetilde{\Omega}_{\varepsilon} \setminus (M_0\cup Q_0),\\
				\phi_1 = \phi_{1,0} \text{ in } M_0, \quad
				\phi_1= -v_1^* \text{ in } Q_0,\quad \phi_1=0 \text{ on } \partial \widetilde{\Omega}_\varepsilon\\
				\phi_2 = \phi_{2,0} \text{ in } Q_0,\quad
				\phi_2= -v_2^* \text{ in } M_0, \quad \phi_2=0 \text{ on } \partial \widetilde{\Omega}_\varepsilon.
			\end{cases}
		\end{equation}
		Let $D:=\{y:\,|y-\widetilde\xi_1|\le|y-\widetilde\xi_2|\}$.  On
		$D\cap\widetilde\Omega_\varepsilon\setminus(M_0\cup Q_0)$ one has
		$|y-\widetilde\xi_2|\ge c\varepsilon^{-1/2}$, hence
		$|v_2|^{\frac{2^*}2}\le C\varepsilon^{\frac N2}$ and, by
		Lemma~\ref{lem4.2},
		\begin{equation}\label{eq4.15}
			|v_1|^{\frac{2^*}{2}-1} |v_2|^{\frac{2^*}{2}}
			\le C\varepsilon^{\frac N2}(d_1+r_1)^{-2}
			\le C\varepsilon^{\frac N4}(d_1+r_1)^{-\frac N2-2},
		\end{equation}
		while on $D^c\cap\widetilde\Omega_\varepsilon\setminus(M_0\cup Q_0)$
		one has $|y-\widetilde\xi_1|\ge c\varepsilon^{-1/2}$, hence
		$|v_1|^{\frac{2^*}2-1}=|v_1|^{\frac2{N-2}}\le C\varepsilon$ and, by
		Lemma~\ref{lem4.2},
		$|v_2|^{\frac{2^*}2}=|v_2|^{\frac N{N-2}}
		\le C(d_2+r_2)^{-N}$, so that
		\begin{equation}\label{eq4.16}
			|v_1|^{\frac{2^*}{2}-1} |v_2|^{\frac{2^*}{2}}
			\le C\varepsilon(d_2+r_2)^{-N}
			\le C\varepsilon(d_2+r_2)^{-N+\theta}.
		\end{equation}
		Here we used $d_i^{\frac N2}\le C$,
		$(d_i+r_i)^{\frac N2}\le C\varepsilon^{-\frac N4}$ in
		\eqref{eq4.15}, and in \eqref{eq4.16} the fact that
		$-N\le-N+\theta$, together with $d_2\le\gamma_*^{-1}$ and
		$|v_2|\le A_2d_2^{\frac{N-2}{2}}(d_2+r_2)^{-(N-2)}$.

		In $\widetilde{\Omega}_{\varepsilon} \setminus (M_0\cup Q_0)$, one has $|y-\tilde a_1|\ge R_0$ and $|y-\tilde{\xi}_1|\le C\varepsilon^{-\frac{1}{2}}$.
		Using the mean value theorem and Lemma \ref{lem3.2}, we obtain 
		\begin{equation}\label{eq4.17}
			\begin{aligned}
				&\Bigl|\mu_1 |v_1|^{2^* - 2}v_1 - \mu_1^{-\frac{1}{2^*-2}} U_{d_1,\tilde{\xi}_1}^{2^*-1}\Bigr|\\ 
				= &\Bigl|\mu_1 |v_1|^{2^* - 2}v_1 - \mu_1|v_1^*|^{2^* - 2}v_1^* + \mu_1|v_1^*|^{2^* - 2}v_1^* - \mu_1^{-\frac{1}{2^*-2}} U_{d_1,\tilde{\xi}_1}^{2^*-1}\Bigr| \\
				\leq& \mu_1(2^*-1) \bigl(|v_1|^{2^* - 2} + |v_1^*|^{2^* - 2}\bigr) |\phi_1| + \mu_1^{-\frac{1}{2^*-2}} \Bigl[ U_{d_1,\tilde{\xi}_1}^{2^*-1} - (\widetilde{P}_{\varepsilon} U_{d_1,\tilde{\xi}_1})^{2^*-1} \Bigr] \\
				\leq&2 \mu_1(2^*-1) A_1^{2^*-2} \varpi^{(1)}_4 |\phi_1|
				+ C\varepsilon^{\frac{N-2}{2}} \varpi^{(1)}_4
				+ C \varepsilon^{\frac{N-2}{2}}\frac{\varpi^{(1)}_4}{|y-\tilde{a}_1|^{N-2}}\\
				\leq& \mathcal V |\phi_1| + C\varepsilon^{\frac{N}{4}} \varpi^{(1)}_{\frac N2+2}.
			\end{aligned}
		\end{equation}
		The last step uses $|y-\widetilde a_1|\ge R_0/2$ and
		$d_1+|y-\widetilde\xi_1|\le C\varepsilon^{-1/2}$, which give
		$\varpi^{(1)}_4=d_1^2(d_1+r_1)^{-4}\le C\varepsilon^{-\frac N4+1}\varpi^{(1)}_{\frac N2+2}$;
		hence the two error terms, both of size
		$O(\varepsilon^{\frac{N-2}{2}}\varpi^{(1)}_4)$, are bounded by
		$C\varepsilon^{\frac N4}\varpi^{(1)}_{\frac N2+2}$.  Here
		$\mathcal V:=2\mu_1(2^*-1)A_1^{2^*-2}\varpi^{(1)}_4$ is the potential of the
		linearized equation.
		
		Combining \eqref{eq4.14}--\eqref{eq4.17} and applying Kato's inequality, we obtain 
		\begin{equation}\label{eq4.18}
			-\Delta |\phi_1| \leq \mathcal V|\phi_1|
			+ C\varepsilon^{\frac{N}{4}} \varpi^{(1)}_{\frac N2+2}
			+ C\varepsilon\, \varpi^{(2)}_{N-\theta}
			\quad\text{in }\widetilde \Omega_{\varepsilon} \setminus\Big(M_0\cup Q_0\Big),
		\end{equation}

		Set
		\[
		W_1:=c_1\varpi^{(1)}_{N/2}+c_2\varpi^{(2)}_{N-2-\theta},
		\qquad
		c_1:=R_0^{N/2}\varepsilon^{\frac{N-2\theta}{4}}d_1^{-\frac N4},
		\qquad
		c_2:=\varepsilon\,d_2^{-\frac{N-2-\theta}{2}}.
		\]
		Both exponents $\frac N2$ and $N-2-\theta$ lie in $(0,N-2]$, and
		$|y-\widetilde\xi_i|\ge R_0/2$ on
		$\widetilde\Omega_\varepsilon\setminus(M_0\cup Q_0)$, so
		\eqref{eq4.2} of Lemma~\ref{lem4.1} applies with
		$z=\widetilde\xi_i$ and $d=d_i$.  Multiplying the two bounds by the
		constants $c_i$, which do not affect $\varpi^{(i)}_4$, and adding, we obtain
		\[
		\begin{aligned}
			-\Delta W_1\ge{}&c_0R_0\Bigl[\tfrac N2\varpi^{(1)}_4c_1\varpi^{(1)}_{N/2}
			+(N-2-\theta)\varpi^{(2)}_4c_2\varpi^{(2)}_{N-2-\theta}\Bigr]\\
			&+\tfrac{N(N-4)}4c_1d_1^{-1}\varpi^{(1)}_{\frac N2+2}
			+(N-2-\theta)\theta\,c_2d_2^{-1}\varpi^{(2)}_{N-\theta},
		\end{aligned}
		\]
		the exact terms $\alpha(N-2-\alpha)d_i^{-1}\varpi^{(i)}_{\alpha+2}$ of
		\eqref{eq4.2} producing the last two expressions, the factors $d_i^{-1}$
		being bounded and their coefficients $\frac{N(N-4)}4$, $(N-2-\theta)\theta$
		positive because $N\ge5$ and $0<\theta<N-2$.  The two weights in the bracket
		are comparable: on $D=\{r_1\le r_2\}$ one has $r_2\ge c\varepsilon^{-1/2}$,
		hence $\varpi^{(2)}_4\le C\varpi^{(1)}_4$ and
		$c_2\varpi^{(2)}_{N-2-\theta}\le Cc_1\varpi^{(1)}_{N/2}$, while on $D^c$ the
		same two inequalities hold with the indices exchanged.  We shall use the
		elementary two-center inequality
		\begin{equation}\label{eq4.19}
			a_1q_1+a_2q_2\ge\frac12\,(a_1+a_2)(q_1+q_2)
			\qquad\hbox{whenever}\qquad (a_1-a_2)(q_1-q_2)\ge0,
		\end{equation}
		whose two sides differ by $\frac12(a_1-a_2)(q_1-q_2)$.  The pairs
		$a_i:=\varpi^{(i)}_4$, $q_1:=c_1\varpi^{(1)}_{N/2}$, $q_2:=c_2\varpi^{(2)}_{N-2-\theta}$
		satisfy its hypothesis up to fixed constants by the comparison above, and
		\eqref{eq4.19} gives
		\[
		\bigl(\varpi^{(1)}_4+\varpi^{(2)}_4\bigr)W_1
		\le C\Bigl[\varpi^{(1)}_4c_1\varpi^{(1)}_{N/2}
		+\varpi^{(2)}_4c_2\varpi^{(2)}_{N-2-\theta}\Bigr]
		\]
		throughout $\widetilde\Omega_\varepsilon\setminus(M_0\cup Q_0)$.  Inserting
		this into the previous display and absorbing the constants into $c$ we obtain
		\begin{equation}\label{eq4.20}
			\begin{aligned}
				-\Delta W_1\ge{}&cR_0\bigl(\varpi^{(1)}_4+\varpi^{(2)}_4\bigr)W_1\\
				&+\frac{N(N-4)}{4}c_1d_1^{-1}\varpi^{(1)}_{\frac N2+2}
				+(N-2-\theta)\theta\,c_2d_2^{-1}\varpi^{(2)}_{N-\theta},
			\end{aligned}
		\end{equation}
		In $Q_0$ one has $v_1=0$ by \eqref{eq4.9},
		$|y-\widetilde\xi_1|\ge c\varepsilon^{-1/2}$ and
		$d_2+|y-\widetilde\xi_2|\le2R_0$, so
		\begin{equation}\label{eq4.21}
			|\phi_1|=|v_1^*|\le C\varepsilon^{\frac{N-2}{2}}
			\le c\varepsilon\le c_2\varpi^{(2)}_{N-2-\theta}\le W_1
			\qquad\text{in }Q_0 .
		\end{equation}
		Comparing \eqref{eq4.20} with the right-hand side of \eqref{eq4.18} term
		by term, the first term of \eqref{eq4.20} dominates $\mathcal VW_1$ as soon
		as $cR_0\ge2\mu_1(2^*-1)A_1^{2^*-2}$, while the last two terms dominate the
		two error terms up to a fixed constant $C_1$, since
		$c_1d_1^{-1}\ge cR_0^{N/2}\varepsilon^{\frac{N-2\theta}{4}}$ with
		$\frac N4>\frac{N-2\theta}{4}$, and $c_2d_2^{-1}\ge c\varepsilon$.  We first
		compare the two functions on
		$\partial(\widetilde\Omega_\varepsilon\setminus(M_0\cup Q_0))$, which fixes
		the multiple of $W_1$.  On $\partial M_0$ one has $\phi_1=\phi_{1,0}$ and
		$d_1+|y-\widetilde\xi_1|\le2R_0$, hence
		$W_1\ge c_1\varpi^{(1)}_{N/2}
		\ge cR_0^{N/2}\varepsilon^{\frac{N-2\theta}{4}}(2R_0)^{-N/2}
		=c\varepsilon^{\frac{N-2\theta}{4}}\ge c|\phi_{1,0}|$; on $\partial Q_0$ this
		is \eqref{eq4.21}, and finally $\phi_1=0$ on
		$\partial\widetilde\Omega_\varepsilon$.  Accordingly $|\phi_1|\le\bar CW_1$
		on this boundary, with $\bar C$ independent of $\varepsilon$ and of
		$(\mathbf d,\boldsymbol\tau)$.  Let $C_0:=\max\{\bar C,C_1\}$.  Then
		\[
		-\Delta(|\phi_1|-C_0W_1)\le \mathcal V(|\phi_1|-C_0W_1)
		\quad\text{in }\widetilde\Omega_{\varepsilon}\setminus(M_0\cup Q_0),
		\]
		while $|\phi_1|-C_0W_1\le0$ on the boundary.
		For the potential, \eqref{eq4.3} of Lemma~\ref{lem4.1},
		applied with $z=\widetilde\xi_1$ and $d=d_1$, gives
		$\|\varpi^{(1)}_4\|_{L^{N/2}}\le CR_0^{-2}$; since
		$\mathcal V=2\mu_1(2^*-1)A_1^{2^*-2}\varpi^{(1)}_4$, also
		\[
		\|\mathcal V\|_{L^{{N}/{2}}\bigl(\widetilde{\Omega}_{\varepsilon}
			\setminus (M_0 \cup Q_0)\bigr)}\le C R_0^{-2}.
		\]
		Choosing $R_0$ sufficiently large makes this norm arbitrarily small, which validates the hypothesis of the comparison principle.
		
		By Lemma \ref{lem3.1}, we conclude that
		$$
		|\phi_1| \leq C_0W_1 \quad \text{in } \widetilde \Omega_{\varepsilon} \setminus \Big(M_0\cup Q_0\Big).
		$$
		In $Q_0$ the estimate \eqref{eq4.21} gives $|\phi_1|\le W_1\le C_0W_1$, so
		$$
		|\phi_1| \leq C_0W_1 \quad \text{in } \widetilde \Omega_{\varepsilon}\setminus M_0.
		$$
		A similar argument applies to obtain the corresponding estimate for $\phi_2$.
	\end{proof}

	\begin{Cor}\label{cor4.4}
		Under the assumptions of Lemma~\ref{lem4.3}, there
		exists $C>0$, independent of $\varepsilon$ and uniform for
		$(\mathbf d,\boldsymbol\tau)\in \mathcal P_{\gamma_*}$ and
		$(\phi_{1,0},\phi_{2,0})\in\Lambda_\varepsilon$, and there is a constant
		$c>0$, depending only on $N$ and $\gamma_*$, small enough that, for every
		$y\in\widetilde\Omega_\varepsilon\setminus(M_0\cup Q_0)$ with
		$|y-\widetilde\xi_1|\le c\,\varepsilon^{-1/2}$,
		\begin{equation*}
			|v_1(y)-v_1^*(y)|
			\le C R_0^{N/2}\varepsilon^{\frac{N-2\theta}{4}}
			(d_1+|y-\widetilde\xi_1|)^{-N/2},
		\end{equation*}
		and, for every $y$ with $R_0+1\le|y-\widetilde a_1|$ and $|y-\widetilde\xi_1|\le c\,\varepsilon^{-1/2}$,
		\begin{equation*}
			|\nabla(v_1-v_1^*)(y)|
			\le C R_0^{N/2}\varepsilon^{\frac{N-2\theta}{4}}
			(d_1+|y-\widetilde\xi_1|)^{-N/2}.
		\end{equation*}
		The same two estimates hold for $v_2-v_2^*$, with
		$Q_0,d_2,\widetilde\xi_2,\phi_{2,0}$ in place of
		$M_0,d_1,\widetilde\xi_1,\phi_{1,0}$.
		
		Moreover
		$\|v_1-v_1^*\|_{L^\infty(\partial M_3)}
		\le CR^{-N/2}\varepsilon^{\frac{N-2\theta}{4}}$ and
		$\|v_1-v_1^*\|_{H^{1/2}(\partial M_3)}
		\le CR^{-1/2}\varepsilon^{\frac{N-2\theta}{4}}$,
		and analogously for $v_2-v_2^*$ on $\partial Q_3$.
	\end{Cor}
	
	\smallskip
	\begin{Rem}\label{rem4.5}
		The norms on $\partial M_3$ in
		Corollary~\ref{cor4.4} are computed on the $H^1$ trace of $v_1-v_1^*$.
		By elliptic regularity for \eqref{eq4.5} away from the concentration
		points this trace is of class $C^1$ up to $\partial M_3$, hence it
		coincides with the ordinary restriction and its $L^\infty$ norm is the
		supremum of a continuous function.
	\end{Rem}
	
	\begin{proof}
		We prove the two estimates for $v_1-v_1^*$; the case of $v_2-v_2^*$ is identical.  Let
		$y\in\widetilde\Omega_\varepsilon\setminus(M_0\cup Q_0)$ with $|y-\widetilde\xi_1|\le c\,\varepsilon^{-1/2}$, and let $c$ be as in the statement.
		Since the two centers are separated by a distance comparable to
		$\varepsilon^{-1/2}$, one has $|y-\widetilde\xi_2|\ge c^{\prime}\varepsilon^{-1/2}$ for
		$c$ small enough.  Hence the opposite-center term of \eqref{eq4.12}, whose exponent
		is $N-2-\theta$, is bounded by
		\[
		C\varepsilon\bigl(c^{\prime}\varepsilon^{-1/2}\bigr)^{-(N-2-\theta)}
		=C\varepsilon^{\frac N2-\frac{\theta}{2}}.
		\]
		On the other hand, on the collar
		\[
		(d_1+|y-\widetilde\xi_1|)^{-N/2}\ge(c\varepsilon^{-1/2})^{-N/2}
		=c^{-N/2}\varepsilon^{\frac N4},
		\]
		so the main term of \eqref{eq4.12} is at least
		$c_0(R_0/c)^{N/2}\varepsilon^{\frac{N-\theta}{2}}$.  Since the opposite-center
		term is $O(\varepsilon^{\frac{N-\theta}{2}})$ with a fixed constant, the main
		term dominates it as soon as $R_0$ is large enough.  Thus the main term
		dominates the opposite-center term, which proves the pointwise bound for
		$v_1-v_1^*$.
		
		For the gradient, let $y$ satisfy $R_0+1\le|y-\widetilde a_1|$ and
		$|y-\widetilde\xi_1|\le c\,\varepsilon^{-1/2}$, so that
		$B_1(y)\subset\widetilde\Omega_\varepsilon\setminus(M_0\cup Q_0)$.  Put
		$\phi_1=v_1-v_1^*$; by \eqref{eq4.14} it solves $-\Delta\phi_1=f$ on
		$B_1(y)$, where $f$ is the right-hand side displayed below.  The interior
		gradient estimate \cite[Theorem 3.9]{GT01} gives
		\[
		|\nabla\phi_1(y)|\le C\Big(\sup_{B_1(y)}|\phi_1|
		+\sup_{B_1(y)}\Big|\mu_1|v_1|^{2^*-2}v_1-\mu_1^{-\frac{1}{2^*-2}}U_{d_1,\tilde\xi_1}^{2^*-1}
		+\beta|v_1|^{\frac{2^*}2-2}v_1|v_2|^{\frac{2^*}2}\Big|\Big).
		\]
		For $z\in B_1(y)$ one has
		$d_1+|z-\widetilde\xi_1|\ge\frac12(d_1+|y-\widetilde\xi_1|)$, so
		\eqref{eq4.15} and \eqref{eq4.17} give
		\[
		\begin{aligned}
			\sup_{B_1(y)}\Big|\mu_1|v_1|^{2^*-2}v_1
			&-\mu_1^{-\frac{1}{2^*-2}}U_{d_1,\tilde\xi_1}^{2^*-1}
			+\beta|v_1|^{\frac{2^*}2-2}v_1|v_2|^{\frac{2^*}2}\Big|\\
			&\le C R_0^{-4}\sup_{B_1(y)}|\phi_1|
			+C\varepsilon^{\frac N4}(d_1+|y-\tilde\xi_1|)^{-\frac N2-2}.
		\end{aligned}
		\]
		Inserting the pointwise bound on $B_1(y)$ and using
		$\varepsilon^{\frac N4}(d_1+|y-\widetilde\xi_1|)^{-2}
		\le\varepsilon^{\frac\theta2}R_0^{N/2}
		\varepsilon^{\frac{N-2\theta}{4}}$, we obtain
		the estimate for 
		$|\nabla(v_1-v_1^*)(y)|$.
		
		For the boundary values, on
		$\{3R<|y-\widetilde a_1|<3R+1\}$ one has
		$|y-\widetilde\xi_1|=3R+O(1)$, so the pointwise and gradient bounds above give
		$|v_1-v_1^*|+|\nabla(v_1-v_1^*)|\le CR^{-N/2}\varepsilon^{\frac{N-2\theta}{4}}$.  Evaluating the pointwise bound on $\partial B_{3R}(\widetilde a_1)$
		gives the $L^\infty$ bound for $v_1-v_1^*$ on $\partial M_3$; integrating the
		gradient bound over the collar $\{3R<|y-\widetilde a_1|<3R+1\}$, of volume at most $CR^{N-1}$, and applying the trace
		theorem on $\partial B_{3R}(\widetilde a_1)$ gives
		$C(R^{N-1})^{1/2}R^{-N/2}=CR^{-1/2}$ times
		$\varepsilon^{\frac{N-2\theta}{4}}$, the $H^{1/2}$ bound.
	\end{proof}

	\subsection{Dead-core formation}
	By the zero-extension convention above, the first component is already zero
	in the artificial core around $\widetilde a_2$, and the second is zero
	in the artificial core around $\widetilde a_1$.  The next lemma quantifies
	how far these zero sets fill the balls centered at the hole points; it is
	based on the dead-core criterion \cite[Proposition 3.1]{TVY}.
	\begin{Lem}\label{lem4.6}
		Let $R_0$ be fixed as above, let $(\phi_{1,0},\phi_{2,0})\in\Lambda_\varepsilon$,
		and let $(v_1,v_2)$ satisfy \eqref{eq4.6}, \eqref{eq4.7}, \eqref{eq4.9}
		and \eqref{eq4.8}.  Then there exists $c_1>0$, depending only on
		$N,\mu_1,\mu_2,\beta,\gamma_*$ and $R_0$, but not on $\varepsilon$, $R$ or
		$(\phi_{1,0},\phi_{2,0})$, such that for every $\rho$ with
		\[
		2R_0\le\rho\le c_1\varepsilon^{-\frac{N-4}{2(N-2)}}
		\]
		one has
		\[
		v_1=0\ \text{a.e.\ in }B_\rho(\widetilde a_2),
		\qquad
		v_2=0\ \text{a.e.\ in }B_\rho(\widetilde a_1).
		\]
	\end{Lem}
	\begin{proof}
		We prove the assertion for $v_2$; the case of $v_1$ is identical.  We show
		the existence of a dead core with radius
		$\rho=c_1\varepsilon^{-\frac{N-4}{2(N-2)}}$, where $c_1$ will be fixed later;
		the statement for smaller $\rho$ then follows by monotonicity in $\rho$.
		Since $\frac{N-4}{2(N-2)}<\frac12$ and $|\widetilde\xi_1-\widetilde a_1|<1$,
		every $y\in B_{2\rho}(\widetilde a_1)\setminus\overline{M_0}$ belongs to
		$\widetilde\Omega_\varepsilon\setminus(M_0\cup Q_0)$ and satisfies
		$|y-\widetilde\xi_1|\le2\rho+1\le c\,\varepsilon^{-1/2}$ for
		small $\varepsilon$, where $c$ is the constant of Corollary~\ref{cor4.4}; hence
		that corollary is available on the whole ball considered below.
		Since the two centers satisfy
		$c\varepsilon^{-1/2}\le|\widetilde a_1-\widetilde a_2|\le C\varepsilon^{-1/2}$,
		\[
		|y-\widetilde\xi_2|\ge c\,\varepsilon^{-\frac12},
		\qquad y\in B_{2\rho}(\widetilde a_1).
		\]
		
		By Lemma~\ref{lem4.2} and \eqref{eq4.9},
		\[
		|v_2(y)|\le A_2d_2^{\frac{N-2}{2}}(d_2+|y-\widetilde\xi_2|)^{2-N}
		\le C\varepsilon^{\frac{N-2}{2}},
		\qquad y\in B_{2\rho}(\widetilde a_1)\setminus\overline{M_0},
		\]
		with $v_2=0$ in $\overline{M_0}$; the same bound holds on
		$\partial B_{2\rho}(\widetilde a_1)$ in the trace sense.
		
		For $y\in B_{2\rho}(\widetilde a_1)\setminus\overline{M_0}$ one has
		$\gamma_*\le d_1+|y-\widetilde\xi_1|\le5\rho$, and the projection defect is
		negligible at that scale: by Corollary~\ref{corB.2},
		\[
		v_1^*(y)=\mu_1^{-\frac{1}{2^*-2}}\widetilde{P}_\varepsilon
		U_{d_1,\widetilde\xi_1}(y)\ge c\,(d_1+|y-\widetilde\xi_1|)^{2-N},
		\]
		because $|y-\widetilde a_1|\ge R_0$ on this set, so that
		$|U_{d_1,\widetilde\xi_1}-\widetilde{P}_\varepsilon
		U_{d_1,\widetilde\xi_1}|\le C\varepsilon^{\frac{N-2}{2}}
		\le\frac c2(d_1+|y-\widetilde\xi_1|)^{2-N}$, while
		$(d_1+|y-\widetilde\xi_1|)^{2-N}\ge(5\rho)^{2-N}
		\ge c\varepsilon^{\frac{N-4}{2}}\gg\varepsilon^{\frac{N-2}{2}}$.
		On the other hand, Corollary \ref{cor4.4} yields
		\begin{align*}
			|v_1-v_1^*|
			&\le C R_0^{\frac N2}\varepsilon^{\frac{N-2\theta}{4}}
			(d_1+|y-\widetilde\xi_1|)^{-\frac N2} \\
			&\le C\varepsilon^{\frac{N-2\theta}{4}}\rho^{\frac{N-4}{2}}
			(d_1+|y-\widetilde\xi_1|)^{2-N}
			, \quad y\in B_{2\rho}(\widetilde a_1)\setminus\overline{M_0},
		\end{align*}
		where we used $d_1+|y-\widetilde\xi_1|\le5\rho$.
		Since
		\[
		\frac{N-2\theta}{4}-\frac{(N-4)^2}{4(N-2)}
		=\frac14\left[(N-2\theta)-\frac{(N-4)^2}{N-2}\right]>0,
		\]
		the factor $\varepsilon^{\frac{N-2\theta}{4}}\rho^{\frac{N-4}{2}}$ tends to
		zero, and consequently
		\[
		v_1\ge c\,(d_1+|y-\widetilde\xi_1|)^{2-N}\ge c\,\rho^{2-N},
		\qquad y\in B_{2\rho}(\widetilde a_1)\setminus\overline{M_0}.
		\]
		
		On $B_{2\rho}(\widetilde a_1)\setminus\overline{M_0}$, Kato's inequality,
		$\beta<0$ and the two displays above give
		\[
		-\Delta|v_2|
		\le \mu_2|v_2|^{2^*-1}-c|v_1|^{\frac N{N-2}}|v_2|^{\frac2{N-2}}
		\le \mu_2|v_2|^{2^*-1}-c\,\rho^{-N}|v_2|^{\frac2{N-2}}.
		\]
		The first term is harmless.  Indeed
		$|v_2|\le C\varepsilon^{\frac{N-2}{2}}$ gives
		$\mu_2|v_2|^{\frac N{N-2}}\le C\varepsilon^{\frac N2}$, while
		$\rho^{-N}\ge c_1^{-N}\varepsilon^{\frac{N(N-4)}{2(N-2)}}$, and
		$\frac{N(N-4)}{2(N-2)}<\frac N2$; hence
		\[
		\mu_2|v_2|^{2^*-1}
		=\mu_2|v_2|^{\frac N{N-2}}|v_2|^{\frac2{N-2}}
		\le\frac c2\,\rho^{-N}|v_2|^{\frac2{N-2}}
		\]
		for $\varepsilon$ small.  Therefore
		\[
		-\Delta|v_2|\le-\frac c2\,\rho^{-N}|v_2|^{\frac2{N-2}}
		\qquad\text{in }B_{2\rho}(\widetilde a_1)\setminus\overline{M_0}.
		\]
		Since $|v_2|=0$ in $\overline{M_0}$, both sides extend by zero to
		$\overline{M_0}$, and $|v_2|\in H^1(B_{2\rho}(\widetilde a_1))$; thus the
		inequality holds weakly in the whole ball,
		\[
		-\Delta|v_2|\le-\frac c2\,\rho^{-N}|v_2|^{\frac2{N-2}}
		\qquad\text{weakly in }B_{2\rho}(\widetilde a_1).
		\]
		
		Apply \cite[Proposition 3.1]{TVY} on $B_{2\rho}(\widetilde a_1)$ with
		exponent $p=\frac2{N-2}$ and radii $\rho$ and $2\rho$.  Its hypotheses are the
		two displays just established, and its numeric hypothesis, obtained by
		applying the power $1-p=\frac{N-4}{N-2}>0$ to the bound
		$C\varepsilon^{\frac{N-2}{2}}$ at radius $2\rho$, reads
		\[
		c\,\rho^{-N}\ \ge\ C\left(\varepsilon^{\frac{N-2}{2}}
		\right)^{\frac{N-4}{N-2}}\rho^{-2}
		=C\,\varepsilon^{\frac{N-4}{2}}\rho^{-2}.
		\]
		With $\rho=c_1\varepsilon^{-\frac{N-4}{2(N-2)}}$ both sides carry the same
		power of $\varepsilon$, because
		$\frac{N(N-4)}{2(N-2)}=\frac{N-4}{2}+\frac{N-4}{N-2}$, and the inequality
		reduces to $c\,c_1^{-N}\ge Cc_1^{-2}$, i.e.\ to $c_1^{N-2}\le c/C$.  Hence it
		holds for $c_1$ small, in terms of $N,\mu_1,\mu_2,\beta,\gamma_*$, $R_0$ and
		$\theta$ only, and
		\[
		v_2=0\quad\text{a.e.\ in }B_\rho(\widetilde a_1).
		\]
	\end{proof}
	\begin{Rem}\label{rem4.7}
		Since $R>4R_0$ is fixed and $\rho=4R$ obeys the displayed
		condition once $\varepsilon$ is small enough, and since
		$M_3\subset B_{4R}(\widetilde a_1)$ and $Q_3\subset B_{4R}(\widetilde a_2)$,
		\[
		v_2=0\ \text{ in }\ M_3,
		\qquad
		v_1=0\ \text{ in }\ Q_3.
		\] 
		As noted after \eqref{eq4.8}, the vanishing $v_2=0$ in $M_3$ and
		$v_1=0$ in $Q_3$, which follows from Lemma~\ref{lem4.6}, promotes any
		solution of the intermediate problem \eqref{eq4.8} to a solution of the
		target problem \eqref{eq4.5}.
	\end{Rem}

	Let $(\phi_{1,0},\phi_{2,0}),(\bar\phi_{1,0},\bar\phi_{2,0})\in\Lambda_\varepsilon$,
	and let $(v_1,v_2)$ and $(\bar v_1,\bar v_2)$ be the corresponding
	zero-extended outer solutions of \eqref{eq4.8} satisfying
	\eqref{eq4.6} and \eqref{eq4.7}. 
	\begin{Lem}\label{lem4.8}
		Set $(w,z) = (v_1,v_2) - (\bar v_1,\bar v_2)$ and
		$(w_0,z_0) = (\phi_{1,0},\phi_{2,0}) - (\bar{\phi}_{1,0},\bar{\phi}_{2,0})$.
		Upon increasing $R_0$ if necessary, still independently of $\varepsilon$,
		there exist constants $C>0$, independent of $\varepsilon$ and $R_0$, such that for any $y\in \widetilde{\Omega}_{\varepsilon}\setminus(M_0\cup Q_0)$
		\begin{align*}
			|w(y)|\leq &C R_0^{N-2} \bigl(\| w_0\|_{L^{\infty}(\partial M_0)} + \| z_0\|_{L^{\infty}(\partial Q_0)}\bigr) 
			\left(\frac{d_1^{\frac{N-2}{2}}}{(d_1 + |y - \widetilde{\xi}_1|)^{N-2}}+\frac{d_2^{\frac{N-2}{2}}}{(d_2 + |y - \widetilde{\xi}_2|)^{N-2}}\right), 
			\\
			|z(y)|\leq &C R_0^{N-2} \bigl(\| w_0\|_{L^{\infty}(\partial M_0)} + \| z_0\|_{L^{\infty}(\partial Q_0)}\bigr) 
			\left(\frac{d_1^{\frac{N-2}{2}}}{(d_1 + |y - \widetilde{\xi}_1|)^{N-2}}+\frac{d_2^{\frac{N-2}{2}}}{(d_2 + |y - \widetilde{\xi}_2|)^{N-2}}\right).
		\end{align*}
	\end{Lem}
	
	\begin{proof}
		By Lemma~\ref{lem4.6} the zero extensions satisfy the equations across   $\partial B_{R_0}(\widetilde a_i)$, $i=1,2$, so the equations
		for $(v_1,v_2)$ and $(\bar v_1,\bar v_2)$ may be subtracted on
		$\widetilde{\Omega}_{\varepsilon}\setminus(M_0\cup Q_0)$.  The difference
		$(w,z)$ then satisfies
		\begin{equation}\label{eq4.22}
			\begin{aligned}
				-\Delta w ={}& \mu_1\bigl(|v_1|^{2^*-2}v_1 -|\bar v_1|^{2^*-2}\bar v_1\bigr)
				+ \beta\big(\lvert v_2\rvert^{\frac{2^{*}}{2}}\lvert v_1\rvert^{\frac{2^{*}}{2}-2}v_1-\lvert \bar v_2\rvert^{\frac{2^{*}}{2}}\lvert \bar v_1\rvert^{\frac{2^{*}}{2}-2}\bar v_1\big),\\
				-\Delta z ={}& \mu_2\bigl(|v_2|^{2^*-2}v_2 -|\bar v_2|^{2^*-2}\bar v_2\bigr)
				+ \beta\big(\lvert v_1\rvert^{\frac{2^{*}}{2}}\lvert v_2\rvert^{\frac{2^{*}}{2}-2}v_2-\lvert \bar v_1\rvert^{\frac{2^{*}}{2}}\lvert \bar v_2\rvert^{\frac{2^{*}}{2}-2}\bar v_2\big),
			\end{aligned}
		\end{equation}
		with $w=w_0$, $z=0$ in $M_0$ and $z=z_0$, $w=0$ in $Q_0$.  Note that
		$\varpi^{(i)}_4=d_i^{2}(d_i+|y-\widetilde\xi_i|)^{-4}
		=(\varpi^{(i)}_{N-2})^{2^*-2}$.
		
		By Lemma~\ref{lem4.2},
		\begin{equation}\label{eq4.23}
			\begin{aligned}
				|v_1|^{2^*-2}v_1 - |\bar v_1|^{2^*-2}\bar v_1
				&\le (2^*-1)\bigl(|v_1|^{2^*-2}+|\bar v_1|^{2^*-2}\bigr)|w|
				\le 2(2^*-1)A_1^{2^*-2}\varpi^{(1)}_4|w|,\\
				|v_2|^{\frac{2^*}{2}}|v_1|^{\frac{2^*}{2}-2}v_1
				-|\bar v_2|^{\frac{2^*}{2}}|\bar v_1|^{\frac{2^*}{2}-2}\bar v_1
				&=\bigl(|v_2|^{\frac{2^*}{2}}-|\bar v_2|^{\frac{2^*}{2}}\bigr)|\bar v_1|^{\frac{2^*}{2}-2}\bar v_1
				+|v_2|^{\frac{2^*}{2}}\bigl(|v_1|^{\frac{2^*}{2}-2}v_1-|\bar v_1|^{\frac{2^*}{2}-2}\bar v_1\bigr)\\
				&\le |v_2|^{\frac{2^*}{2}}\bigl(|v_1|^{\frac{2^*}{2}-2}v_1-|\bar v_1|^{\frac{2^*}{2}-2}\bar v_1\bigr)
				+C\sqrt{\varpi^{(1)}_4\varpi^{(2)}_4}\,|z|,
			\end{aligned}
		\end{equation}
		the last step by Lemma~\ref{lem4.2} applied to $v_2$, $\bar v_2$ and
		$\bar v_1$.  Since $t\mapsto|t|^{\frac{2^*}{2}-2}t$ is increasing and
		$\beta<0$, the first term on the right of \eqref{eq4.23}, multiplied by the
		sign of $w$, is non-positive; the same holds for the corresponding term in
		the second equation of \eqref{eq4.22}.  Substituting \eqref{eq4.23} into
		\eqref{eq4.22}, Kato's inequality gives
		\[
		\begin{aligned}
			-\Delta|w|&\le 2(2^*-1)A_1^{2^*-2}\mu_1\varpi^{(1)}_4|w|
			+C\sqrt{\varpi^{(1)}_4\varpi^{(2)}_4}\,|z|,\\
			-\Delta|z|&\le 2(2^*-1)A_2^{2^*-2}\mu_2\varpi^{(2)}_4|z|
			+C\sqrt{\varpi^{(1)}_4\varpi^{(2)}_4}\,|w|,
		\end{aligned}
		\qquad y\in\widetilde{\Omega}_{\varepsilon}\setminus(M_0\cup Q_0).
		\]
		Adding the two inequalities and using
		$2\sqrt{\varpi^{(1)}_4\varpi^{(2)}_4}\le\varpi^{(1)}_4+\varpi^{(2)}_4$, we obtain
		\begin{equation}\label{eq4.24}
			-\Delta(|w|+|z|)\le C(\varpi^{(1)}_4+\varpi^{(2)}_4)(|w|+|z|)
			\qquad\text{in }\widetilde{\Omega}_\varepsilon\setminus(\overline{M_0}\cup\overline{Q_0}).
		\end{equation}
		
		Set
		\[
		A_0:=\|w_0\|_{L^\infty(\partial M_0)}+\|z_0\|_{L^\infty(\partial Q_0)},
		\qquad
		W_2:=C R_0^{N-2}A_0\bigl(\varpi^{(1)}_{N-2}+\varpi^{(2)}_{N-2}\bigr),
		\]
		where $C$ is chosen so that $W_2\ge|w|+|z|$ on the boundary of
		$\widetilde\Omega_\varepsilon\setminus(\overline{M_0}\cup\overline{Q_0})$.
		This is possible.  Indeed, on $\partial M_0$ one has $z=0$ and
		$|w|=|w_0|\le A_0$, while $\varpi^{(1)}_{N-2}\ge cR_0^{-(N-2)}$ because
		$|y-\widetilde\xi_1|\le2R_0$; on $\partial Q_0$ the roles of $w$ and $z$
		are exchanged; and both vanish on $\partial\widetilde\Omega_\varepsilon$.
		By \eqref{eq4.1}, applied with $z=\widetilde\xi_i$, $d=d_i$ and
		$\alpha=N-2$,
		\[
		-\Delta\varpi^{(i)}_{N-2}\ge c_0R_0\,{\varpi^{(i)}_{N-2}}^{2^*-1}
		=c_0R_0\,\varpi^{(i)}_4\varpi^{(i)}_{N-2}
		\qquad\text{in }\widetilde\Omega_\varepsilon\setminus(\overline{M_0}\cup\overline{Q_0}).
		\]
		Since $\varpi^{(i)}_4={\varpi^{(i)}_{N-2}}^{2^*-2}$ and $t\mapsto t^{2^*-2}$ is
		increasing, $(\varpi^{(1)}_4-\varpi^{(2)}_4)(\varpi^{(1)}_{N-2}-\varpi^{(2)}_{N-2})\ge0$;
		\eqref{eq4.19} with $a_i=\varpi^{(i)}_4$ and $q_i=\varpi^{(i)}_{N-2}$ gives
		$\varpi^{(1)}_4\varpi^{(1)}_{N-2}+\varpi^{(2)}_4\varpi^{(2)}_{N-2}
		\ge\frac12(\varpi^{(1)}_4+\varpi^{(2)}_4)(\varpi^{(1)}_{N-2}+\varpi^{(2)}_{N-2})$.
		Consequently
		\begin{equation}\label{eq4.25}
			-\Delta W_2\ge C'\bigl(\varpi^{(1)}_4+\varpi^{(2)}_4\bigr)W_2
			\quad\text{in }\widetilde\Omega_\varepsilon\setminus(\overline{M_0}\cup\overline{Q_0}),
			\qquad
			C':=\tfrac{c_0R_0}{2}.
		\end{equation}
		
		On $\widetilde\Omega_\varepsilon\setminus(\overline{M_0}\cup\overline{Q_0})$
		one has $|y-\widetilde\xi_i|\ge R_0/2$, so \eqref{eq4.3} gives
		$\|C(\varpi^{(1)}_4+\varpi^{(2)}_4)\|_{L^{N/2}}
		\le CR_0^{-2}<\mathcal S$ for $R_0$ large.  Since $C'=c_0R_0/2\to\infty$,
		we may also assume $C'\ge C$, the constant in \eqref{eq4.24}.  By
		\eqref{eq4.24}, \eqref{eq4.25} and the choice of $C$,
		\[
		-\Delta\bigl(|w|+|z|-W_2\bigr)
		\le C(\varpi^{(1)}_4+\varpi^{(2)}_4)\bigl(|w|+|z|-W_2\bigr)
		\]
		in $\widetilde\Omega_\varepsilon\setminus(\overline{M_0}\cup\overline{Q_0})$,
		while $|w|+|z|-W_2\le0$ on its boundary.
		Lemma~\ref{lem3.1}, with $v_1=W_2\ge0$ and $v_2=|w|+|z|$, therefore yields
		$|w|+|z|\le W_2$ there, which is the asserted pair of pointwise bounds.
		
	\end{proof}
	\medskip
	
	We next prove the existence of a solution $(v_1,v_2)$ of~\eqref{eq4.8} satisfying~\eqref{eq4.6} and~\eqref{eq4.7}. 
	
	\subsection{Existence of an outer minimizer} 
	
	Let $A_1,A_2>0$ be the constants given in Lemma \ref{lem4.2}. We consider solutions $(v_{1},v_{2})$ to problem \eqref{eq4.8} satisfying \eqref{eq4.6} and
	\begin{equation}\label{eq4.26}
		|v_1(y)|\le 2A_1\frac{d_1^{(N-2)/2}}
		{(d_1+|y-\widetilde\xi_1|)^{N-2}},\quad
		|v_2(y)|\le 2A_2\frac{d_2^{(N-2)/2}}
		{(d_2+|y-\widetilde\xi_2|)^{N-2}},\quad
		y\in\widetilde{\Omega}_{\varepsilon}\setminus(M_0\cup Q_0).
	\end{equation}
	We note that, once $R_0$ is chosen large enough,
	condition \eqref{eq4.26} implies \eqref{eq4.7}.
	
	\medskip
	For any $(\phi_{1,0},\phi_{2,0})\in \Lambda_{\varepsilon}$, denote 
	\begin{equation*}
		\mathbb{S}(\phi_{1,0},\phi_{2,0}):=\Big\{ (v_{1},v_{2})
		\in H_{0}^{1}(\widetilde{\Omega}_{\varepsilon}\setminus\overline{Q_0})\times H_{0}^{1}(\widetilde{\Omega}_{\varepsilon}\setminus\overline{M_0})  \mid (v_{1}, v_{2}) \text{ satisfies } \eqref{eq4.6} \text{ and } \eqref{eq4.26} 
		\Big\}.
	\end{equation*}
	This set is a closed subset of $H_{0}^{1}(\widetilde{\Omega}_{\varepsilon}\setminus\overline{Q_0})\times H_{0}^{1}(\widetilde{\Omega}_{\varepsilon}\setminus\overline{M_0})$ under the induced product norm
	\[
	\|(v_1,v_2)\|^2_{H_{0}^{1}(\widetilde{\Omega}_{\varepsilon}\setminus\overline{Q_0})\times H_{0}^{1}(\widetilde{\Omega}_{\varepsilon}\setminus\overline{M_0})}=\|\nabla v_1\|_2^2+\|\nabla v_2\|_2^2.
	\]
	Moreover, $\mathbb{S}(\phi_{1,0},\phi_{2,0})$ is nonempty, convex, and weakly sequentially closed.  Indeed, the prescribed values on $M_0$ and $Q_0$ admit $H_0^1$ extensions $E_1,E_2$ to the corresponding auxiliary domains; truncate each $E_i$ by the two barriers
	in \eqref{eq4.26}, which changes $E_i$ only on
	$\widetilde\Omega_\varepsilon\setminus(M_0\cup Q_0)$ and leaves the prescribed functions unchanged in
	$M_0$ and $Q_0$.  The boundary traces of $\phi_{i,0}$ already satisfy the smaller barrier
	bounds, so the truncated pair still belongs to $H^1$ and gives an admissible pair.  Convexity is immediate.  If an admissible sequence converges weakly in the product energy space, then, after passing to a subsequence, it converges almost everywhere; the fixed values and the pointwise bounds therefore pass to the weak limit.

	Define a functional $\widetilde{I}\colon\mathbb{S}(\phi_{1,0},\phi_{2,0})\rightarrow \mathbb{R}$ as
	\[
	\begin{aligned}
		\widetilde{I}(v_1,v_2) :=& \frac{1}{2}\int_{\widetilde{\Omega}_{\varepsilon}\setminus\overline{Q_0}} \lvert\nabla v_{1}\rvert^{2} dy+\frac{1}{2}\int_{\widetilde{\Omega}_{\varepsilon}\setminus\overline{M_0}}\lvert\nabla v_{2}\rvert^{2} dy 
		-\frac{\mu_1}{2^{*}}\int_{\widetilde{\Omega}_{\varepsilon}\setminus\overline{Q_0}}\lvert v_1\rvert^{2^{*}}dy -\frac{\mu_2}{2^{*}}\int_{\widetilde{\Omega}_{\varepsilon}\setminus\overline{M_0}}\lvert v_2\rvert^{2^{*}} dy \\
		&-\frac{2\beta}{2^{*}}\int_{\widetilde{\Omega}_{\varepsilon}\setminus\overline{M_0\cup Q_0}}\lvert v_1\rvert^{\frac{2^{*}}{2}}\lvert v_2\rvert^{\frac{2^{*}}{2}}dy,~~(v_{1},v_{2}) \in \mathbb{S}(\phi_{1,0},\phi_{2,0}).
	\end{aligned}
	\]

	The solution to \eqref{eq4.8} can be obtained by solving the following minimization problem
	\begin{equation}\label{eq4.27}
		\min\big\{\widetilde{I}(v_1,v_2) \mid (v_{1},v_{2}) \in \mathbb{S}(\phi_{1,0},\phi_{2,0})\big\}.
	\end{equation}
	
	\begin{Lem}\label{lem4.9}
		Upon increasing $R_0$ if necessary, still independently of $\varepsilon$,
		the minimization problem \eqref{eq4.27} is attained by a solution to
		\eqref{eq4.8} in $\mathbb{S}(\phi_{1,0},\phi_{2,0})$.
	\end{Lem}
	
	\begin{proof}
		Write $\mathbb S:=\mathbb S(\phi_{1,0},\phi_{2,0})$.  Since $\mathbb S$
		is weakly sequentially closed, it suffices to prove that
		$\widetilde I$ is coercive and weakly lower semicontinuous on
		$\mathbb S$.
		
		(\romannumeral1) \emph{Coercivity.}
		By \eqref{eq4.26} the nonlinear integrals over
		$\widetilde\Omega_\varepsilon\setminus(M_0\cup Q_0)$ are bounded,
		while by \eqref{eq4.6} those over $M_0$ and $Q_0$ are fixed
		constants.  Since $\beta<0$, there is $C>0$ such that
		\[
		\widetilde I(v_1,v_2)\ge
		\frac12\|\nabla v_1\|^2_{L^2(\widetilde\Omega_\varepsilon\setminus\overline{Q_0})}
		+\frac12\|\nabla v_2\|^2_{L^2(\widetilde\Omega_\varepsilon\setminus\overline{M_0})}
		-C,
		\]
		so that $\widetilde I(v_1,v_2)\to+\infty$ whenever
		$\|(v_1,v_2)\|_{H_0^1(\widetilde\Omega_\varepsilon\setminus\overline{Q_0})
			\times H_0^1(\widetilde\Omega_\varepsilon\setminus\overline{M_0})}\to\infty$.
		
		(\romannumeral2) \emph{Weak lower semicontinuity.}
		Since $\beta<0$, Fatou's lemma gives the weak lower semicontinuity
		on $\mathbb S$ of $-\frac{2\beta}{2^*}\int|v_1|^{\frac{2^*}2}
		|v_2|^{\frac{2^*}2}$, the integral being taken over
		$\widetilde\Omega_\varepsilon\setminus(M_0\cup Q_0)$.  The
		complementary functional
		$\widetilde J:=\widetilde I+\frac{2\beta}{2^*}\int|v_1|^{\frac{2^*}2}
		|v_2|^{\frac{2^*}2}$ decouples and is of class $C^2$; along the tangent
		space of $\mathbb S$ its second variation equals
		\[
		\widetilde J''[(\nu,\eta),(\nu,\eta)]
		=\|\nabla\nu\|_2^2+\|\nabla\eta\|_2^2
		-(2^*-1)\mu_1\int|v_1|^{2^*-2}\nu^2
		-(2^*-1)\mu_2\int|v_2|^{2^*-2}\eta^2,
		\]
		hence, by H\"older's and Sobolev's inequalities,
		\[
		\widetilde J''[(\nu,\eta),(\nu,\eta)]
		\ge\Big[1-\frac{(2^*-1)\mu_1}{\mathcal S}
		\|v_1\|_{2^*}^{2^*-2}
		-\frac{(2^*-1)\mu_2}{\mathcal S}
		\|v_2\|_{2^*}^{2^*-2}\Big]\|(\nu,\eta)\|^2\ge0,
		\]
		where all integrals and $L^{2^*}$ norms are taken over
		$\widetilde\Omega_\varepsilon\setminus(M_0\cup Q_0)$ and
		$\|(\nu,\eta)\|$ is the product energy norm; the last inequality
		follows from \eqref{eq4.7}, which \eqref{eq4.26} implies once $R_0$
		is large.
		Thus $\widetilde J$ is convex and, being continuous, weakly lower
		semicontinuous; the same then holds for $\widetilde I$, and
		$(\romannumeral1)$ together with the weak sequential closedness of
		$\mathbb S$ yields the minimum in \eqref{eq4.27}.
		
		(\romannumeral3) \emph{The minimizer solves \eqref{eq4.8}.}
		Let $(v_1,v_2)\in\mathbb S$ attain the minimum.  The energy depends on
		$v_i$ only through $|v_i|$, and by \eqref{eq4.6} $v_1>0$ on $M_0$ and
		$v_2>0$ on $Q_0$, the bubble dominating the datum $\phi_{i,0}$; hence
		$(|v_1|,|v_2|)\in\mathbb S$ has the same energy and we may assume
		$v_i\ge0$.
		For $\xi\in C_c^\infty(\widetilde\Omega_\varepsilon\setminus(M_0\cup Q_0))$
		with $\xi\ge0$ and small $t>0$ the pair $(v_1-t\xi,v_2)$ is
		admissible, so $t\mapsto\widetilde I(v_1-t\xi,v_2)$ is minimal at
		$t=0$.  Differentiating at $0+$ and using $\beta<0$ together with
		$\xi\ge0$ gives
		\[
		\int\nabla v_1\cdot\nabla\xi\le\mu_1\int v_1^{2^*-1}\xi,
		\]
		that is, $-\Delta v_1\le\mu_1v_1^{2^*-1}$ weakly in
		$\widetilde\Omega_\varepsilon\setminus(M_0\cup Q_0)$, and likewise
		for $v_2$.  Since $v_i\ge0$, this is \eqref{eq4.10}, so
		Lemma~\ref{lem4.2} applies with \eqref{eq4.6}, \eqref{eq4.7} and
		\eqref{eq4.9} and yields
		\[
		|v_i(y)|\le A_i\varpi^{(i)}_{N-2}(y)
		=A_i\,d_i^{\frac{N-2}{2}}\,(d_i+|y-\widetilde\xi_i|)^{-(N-2)},
		\qquad y\in\widetilde\Omega_\varepsilon\setminus(M_0\cup Q_0),
		\]
		which leaves a margin of a full factor $2$ in \eqref{eq4.26}.  Since
		$\varpi^{(i)}_{N-2}$ is bounded below on compact subsets of
		$\widetilde\Omega_\varepsilon\setminus(M_0\cup Q_0)$, a two-sided
		variation $v_i+s\xi$ with $\xi\in C_c^\infty$ supported there keeps
		$(v_1,v_2)$ in $\mathbb S$ for $|s|$ small; hence the Euler--Lagrange
		equation of $\widetilde I$ at $(v_1,v_2)$, written with $v_i\ge0$, is
		\eqref{eq4.8}.
	\end{proof}
	
	This allows us to introduce the outer solution operator.
	
	For every $(\phi_{1,0},\phi_{2,0})\in\Lambda_\varepsilon$, the outer
	solution is unique.  Indeed, existence follows from the minimization
	argument; if $(v_1,v_2)$ and $(\bar v_1,\bar v_2)$ are two solutions with the
	same prescribed data, then $w_0=z_0=0$, and Lemma~\ref{lem4.8} gives
	$w=z=0$ in $\widetilde\Omega_\varepsilon\setminus(M_0\cup Q_0)$, while both
	differences vanish identically in $M_0\cup Q_0$, where the data coincide.  In
	what follows, the outer components are always understood through the
	zero-extension convention fixed above.

	\begin{Def}\label{def4.10}
		For $(\phi_{1,0},\phi_{2,0})\in\Lambda_\varepsilon$, let
		$(v_1,v_2)$ be the unique zero-extended outer solution of
		\eqref{eq4.8} satisfying \eqref{eq4.6} and \eqref{eq4.7}.  As noted after
		\eqref{eq4.8}, this same pair solves the outer problem \eqref{eq4.5} with
		the data \eqref{eq4.6}, and in what follows $S$ is understood as the
		solution map of \eqref{eq4.5}.  Define
		\[
		S:\Lambda_\varepsilon\longrightarrow
		\mathbf K_{\mathbf d,\boldsymbol\tau,\varepsilon}^{\perp},
		\qquad
		S(\phi_{1,0},\phi_{2,0})
		:=(v_1-v_1^*,v_2-v_2^*).
		\]
	\end{Def}
	
	The map is well defined.  Indeed, by \eqref{eq4.6},
	$v_1-v_1^*=\phi_{1,0}$ in $M_0$ and
	$v_2-v_2^*=\phi_{2,0}$ in $Q_0$, while by \eqref{eq4.9} both components are
	extended by zero, so that $S(\phi_{1,0},\phi_{2,0})\in
	H_0^1(\widetilde\Omega_\varepsilon)^2$.  Since
	$\operatorname{supp}Z^h\subset M_0$ and
	$\operatorname{supp}Y^h\subset Q_0$, one has
	\[
	\int_{\widetilde\Omega_\varepsilon}
	\nabla(v_1-v_1^*)\cdot\nabla Z^h
	=\int_{M_0}\nabla\phi_{1,0}\cdot\nabla Z^h=0,
	\]
	with the analogous identity for the second component.  Thus
	$S(\phi_{1,0},\phi_{2,0})$ belongs to the stated global orthogonal space.
	Moreover, $S$ depends only on the restrictions of its inputs to $M_0$ and
	$Q_0$.

	\section{The inner problem}\label{sec5}

	We focus on the first component; the construction in $Q_3$ is identical.
	Recall that
	\[
	v_1^*:=\mu_1^{-\frac1{2^*-2}}
	\widetilde P_\varepsilon
	U_{d_1,\widetilde\xi_1}.
	\]
	
	For the first component, let $\varphi_{1,0}$ be a prescribed datum on
	$\widetilde\Omega_\varepsilon\setminus M_3$; its trace on $\partial M_3$ is
	the boundary datum, and the same symbol denotes that trace in boundary
	formulas.
	We seek a correction $\varphi_1$ satisfying
	\begin{equation}\label{eq5.1}
		\begin{cases}
			- \Delta (v_1^* + \varphi_1) = \mu_1 |v_1^* + \varphi_1|^{2^*- 2}(v_1^* + \varphi_1)+ \displaystyle\sum_{h=0}^N c_h Z^h , &~\text{in}~ M_3,\\[3mm]
			\varphi_1=\varphi_{1,0}, &~\text{in}~ \widetilde\Omega_{\varepsilon} \setminus M_3,
		\end{cases}
	\end{equation}
	subject to
	\begin{equation}\label{eq5.2}
		\int_{M_3} \nabla \varphi_1 \cdot \nabla Z^h \, \mathrm{d}y = 0, \quad h=0,1,\dots,N.
	\end{equation}
	We regard $\varphi_{1,0}$ as a boundary datum in
	$H^{1/2}(\partial M_3)\cap L^\infty(\partial M_3)$, vanishing on
	$\partial B_{r_1\sqrt\varepsilon}(\widetilde a_1)$; the admissible class
	$\mathcal D_{1,\varepsilon}$ of such data is introduced just before
	Proposition~\ref{prop5.7} below.

	\subsection{Solvability of the inner problem}
	We solve the boundary value problem \eqref{eq5.1} in three steps: we
	first record the linearized operator and its projected inverse, then
	remove the non-zero boundary condition by a harmonic lifting, and finally
	apply the contraction mapping principle in the orthogonal complement
	space.
	
	For $q:=\frac{2N}{N-4}$, we work in the Banach space
	\[
	X_1:=H_0^1(M_3)\cap L^q(M_3),
	\]
	equipped with the norm
	\[
	\|v_1\|_{X_1}
	:=\|\nabla v_1\|_{L^2(M_3)}
	+R^{-1}\|v_1\|_{L^q(M_3)}.
	\]
	
	The exponent $q$ is determined by
	\[
	\frac1q=\frac12-\frac2N,
	\qquad q=\frac{2N}{N-4}.
	\]
	Thus an $L^2$ right-hand side gains exactly the $L^q$ integrability under the
	Dirichlet inverse.  This choice closes the nonlinear scheme in the form
	\[
	X_1\longrightarrow L^2(M_3)\longrightarrow X_1.
	\]

	The inner problem is solved in the subspace of $X_1$ orthogonal to $Z^h$,
	$h=0,\ldots,N$, with the moments taken over $M_3$, rather than over
	$\widetilde\Omega_\varepsilon$ as in $K_1^\perp$,
	\[
	K_{M_3}^\perp
	:=\left\{\varphi\in X_1:
	\int_{M_3}\nabla\varphi\cdot\nabla Z^h\,dy=0,
	\quad h=0,\ldots,N\right\}.
	\]
	The second component is treated in exactly the same way on $Q_3$, with $Y^h$
	in place of $Z^h$; we write
	\[
	X_2:=H_0^1(Q_3)\cap L^q(Q_3),
	\qquad
	K_{Q_3}^\perp
	:=\left\{\varphi\in X_2:
	\int_{Q_3}\nabla\varphi\cdot\nabla Y^h\,dy=0,
	\quad h=0,\ldots,N\right\},
	\]
	and use the analogous notation for all the objects below without further
	comment.
	
	\medskip

	\subsubsection{The linearized inner operator}
	
	Set
	\[
	V_1:=(2^*-1)\bigl(\widetilde P_{\varepsilon}U_{d_1,\widetilde\xi_1}\bigr)^{2^*-2}.
	\]
	Since \(\widetilde P_{\varepsilon}U_{d_1,\widetilde\xi_1}\in L^\infty(M_3)\),
	one has \(V_1\in L^{N/2}(M_3)\). Hölder's inequality and the Sobolev embedding
	\(H_0^1(M_3)\hookrightarrow L^{2^*}(M_3)\) imply that the bilinear form
	\[
	\int_{M_3}
	\Big(
	\nabla v_1\cdot\nabla\varphi_1
	-
	V_1v_1\varphi_1
	\Big)\,dy
	\]
	is bounded on \(H_0^1(M_3)\times H_0^1(M_3)\).
	Let \(L_1\) be the corresponding operator, formally defined by
	\begin{equation}\label{eq5.3}
		L_1v:=-\Delta v-V_1v
		\quad\text{in }M_3.
	\end{equation}
	This is the localized version of the whole-space operator
	\(-\Delta - (2^*-1) U_{d_1,\widetilde\xi_1}^{2^*-2}\), with
	\(U_{d_1,\widetilde\xi_1}\) replaced by
	\(\widetilde P_{\varepsilon}U_{d_1,\widetilde\xi_1}\) and the domain restricted to \(M_3\).
	For the inner problem we always use the homogeneous Dirichlet realization of \(L_1\) on \(X_1\subset H_0^1(M_3)\), together with the local moment conditions defining \(K_{M_3}^\perp\).

	For $v_1 \in X_1$, we have $-\Delta v_1 \in H^{-1}(M_3)$. For the potential term, since $v_1 \in L^q(M_3)$ with $q = \frac{2N}{N-4}$ and $V_1 \in L^{N/2}(M_3)$, H\"older's inequality with $\frac{1}{2} = \frac{2}{N} + \frac{1}{q}$ yields $V_1 v_1 \in L^2(M_3)$. Since $L^2(M_3) \hookrightarrow H^{-1}(M_3)$, we conclude that $L_1 v_1 \in H^{-1}(M_3)$, and the operator $L_1: X_1 \to H^{-1}(M_3)$ is bounded.

	\begin{Lem}\label{lem5.1}
		There exist $R_*>R_0$, $\varepsilon_*>0$, and $C>0$, independent of
		$R\ge R_*$, $0<\varepsilon<\varepsilon_*$, and the finite-dimensional
		parameters in $\mathcal P_{\gamma_*}$, with the following property.  For every
		$f\in L^2(M_3)$ there are a unique $v_1\in K_{M_3}^\perp$ and
		unique $c_0,\ldots,c_N\in\mathbb R$ such that
		\begin{equation}\label{eq5.4}
			L_1v_1=f+\sum_{h=0}^Nc_hZ^h\quad\hbox{in }M_3.
		\end{equation}
		Moreover,
		\begin{equation}\label{eq5.5}
			\|\nabla v_1\|_{L^2(M_3)}+\sum_{h=0}^N|c_h|
			\le C\|f\|_{H^{-1}(M_3)},
		\end{equation}
		and, with $q=2N/(N-4)$,
		\begin{equation}\label{eq5.6}
			\|v_1\|_{X_1}\le CR\|f\|_{L^2(M_3)}.
		\end{equation}
		Existence, uniqueness and \eqref{eq5.5} remain valid for $f\in H^{-1}(M_3)$,
		with $v_1\in H_0^1(M_3)$ satisfying
		$\int_{M_3}\nabla v_1\cdot\nabla Z^h=0$, $h=0,\ldots,N$, in place of
		$K_{M_3}^\perp$; the space $X_1$ and the exponent $q$ enter only
		\eqref{eq5.6}.
	\end{Lem}
	
	\begin{proof} 
		We first prove the augmented estimate \eqref{eq5.5}.  If it
		were false, there would exist
		$R_n\to\infty$, $\varepsilon_n\to0$,
		$d_{1,n}\to d_1^*\in
		[\gamma_*,\gamma_*^{-1}]$, functions
		$v_{1,n}\in K_{M_{3,n}}^\perp$, coefficients $c_{h,n}$, and
		$f_n\in H^{-1}(M_{3,n})$ satisfying
		\begin{equation}\label{eq5.7}
			\|\nabla v_{1,n}\|_{L^2(M_{3,n})}
			+\sum_{h=0}^N|c_{h,n}|=1,\qquad
			\|f_n\|_{H^{-1}(M_{3,n})}\to0,
		\end{equation}
		and
		\begin{equation}\label{eq5.8}
			L_{1,n}v_{1,n}=f_n+\sum_{h=0}^Nc_{h,n}Z_{1,n}^h
			\quad\hbox{in }M_{3,n}.
		\end{equation}
		After passing to a subsequence, $c_{h,n}\to\bar c_h$.  All centers in
		the following rescaling are $\widetilde\xi_{1,n}$; they vary with the
		parameters and with $\varepsilon_n$.  Throughout this rescaling we write
		$U:=U_{1,0}$ and $\psi^h:=\psi^h_{1,0}$, $h=0,\ldots,N$, for the bubble and
		the kernel derivatives at $d=1$, $\widetilde\xi=0$; with this normalization
		\begin{equation}\label{eq5.9}
			d^{\frac{N-2}{2}}\psi^h_{d,\widetilde\xi}
			\bigl(\widetilde\xi+dz\bigr)=d^{-1}\psi^h_{1,0}(z).
		\end{equation}
		
		Define 
		\[
		w_n(z) = d_{1,n}^{\frac{N-2}{2}} v_{1,n}\bigl(d_{1,n} z + \widetilde\xi_{1,n}\bigr), \quad
		z \in \widetilde M_{3,n} := d_{1,n}^{-1}\bigl(M_{3,n} - \widetilde\xi_{1,n}\bigr)
		\]
		and
		\[
		\widetilde M_{3,n} = B_{\rho_n}(-\tau_{1,n})
		\setminus \overline{B_{r_n}(-\tau_{1,n})},
		\]
		where
		\[
		\rho_n = \frac{3R_n}{d_{1,n}}, \qquad
		r_n = \frac{r_1\sqrt{\varepsilon_n}}{d_{1,n}}.
		\]
		The two balls are concentric, since
		$M_{3,n}=B_{3R_n}(\widetilde a_{1,n})\setminus
		\overline{B_{r_1\sqrt{\varepsilon_n}}(\widetilde a_{1,n})}$ and both are
		centered at $\widetilde a_{1,n}$, whose image is $-\tau_{1,n}$ with
		$|\tau_{1,n}|\le\gamma_*$.  Since $R_n \to \infty$, $\varepsilon_n \to 0$,
		and $d_{1,n} \to d_1^* > 0$, we get
		\[
		\rho_n \to \infty, \quad r_n \to 0.
		\]
		
		By scale invariance,
		\[
		\|\nabla w_n\|_{L^2(\widetilde M_{3,n})}
		=\|\nabla v_{1,n}\|_{L^2(M_{3,n})}\le1.
		\]
		Extend $w_n$ by zero to $\mathbb R^N$.  The energy bound and Sobolev's
		inequality give, after taking a subsequence, weak convergence in $D^{1,2}$
		and strong convergence in $L^2$ on every fixed ball to a function $w$.
		The shrinking hole has zero limiting $D^{1,2}$-capacity.  More precisely,
		if $\eta_n$ is zero on $B_{r_n}(-\tau_{1,n})$, equals one outside the
		concentric ball of twice the radius, and
		$|\nabla\eta_n|\le Cr_n^{-1}$, then
		\[
		\int_{\mathbb R^N}|\nabla\eta_n|^2
		\le C r_n^{N-2}\to0.
		\]
		Hence every compactly supported test function may be replaced by
		$\eta_n$ times that test function, and the limiting equation is posed on all
		of $\mathbb R^N$.
		Passing to a subsequence gives
		\begin{equation}\label{eq5.10}
			-\Delta w-(2^*-1)U^{2^*-2}w
			=\sum_{h=0}^N\bar c_h\,\chi_0(d_1^*z)\psi^h
			\quad\hbox{in }\mathbb R^N.
		\end{equation}
		Testing against $\psi^k$ gives
		\[
		0=\sum_{h=0}^N\bar c_h
		\int_{\mathbb R^N}\chi_0(d_1^*z)\psi^h\psi^k\,dz
		=\bar c_k\int_{\mathbb R^N}\chi_0(d_1^*z) |\psi^k|^2\,dz,
		\qquad k=0,\ldots,N,
		\]
		because $\chi_0$ is radial, $\psi^0$ is radial, and $\psi^h$ is odd in
		$z_h$ for $h\ge1$, so all integrals with $h\ne k$ vanish by symmetry while
		the diagonal one is positive.  Hence $\bar c_k=0$ for
		every $k$, and therefore $c_{h,n}\to0$.
		
		The orthogonality condition   means
		\[
		\int_{M_{3,n}} \nabla v_{1,n} \cdot \nabla Z_{1,n}^h \, \mathrm{d}y = 0, \quad h = 0,1,\dots,N.
		\]
		
		Rescaling $Z^h_{1,n}$ exactly as $w_n$, that is, setting
		$\widehat Z^h_n(z):=d_{1,n}^{\frac{N-2}{2}}
		Z^h_{1,n}\bigl(d_{1,n}z+\widetilde\xi_{1,n}\bigr)$, one has
		$\widehat Z^h_n\to d_1^{*\,-1}\chi_0(d_1^*z)\psi^h$ strongly in
		$H^1(\mathbb R^N)$ by \eqref{eq5.9}, since $Z^h_{1,n}$ are supported in the
		fixed ball $B_{2/\gamma_*}$.
		By $\bar c_h=0$, the limiting equation is homogeneous, so the bubble
		nondegeneracy gives $w=\sum_{h=0}^Na_h\psi^h$.  The moment conditions read
		$\int_{\widetilde M_{3,n}}\nabla w_n\cdot\nabla\widehat Z^h_n\,dz=0$.  Since
		$\widehat Z^h_n$ converge strongly in $H^1$ while $\nabla w_n\rightharpoonup\nabla w$
		weakly in $L^2$, passing to the limit
		gives $\int_{\mathbb R^N}\nabla w\cdot\nabla\bigl(\chi_0(d_1^*z)\psi^k\bigr)
		dz=0$ for every $k$, the constant factor $d_1^{*\,-1}$ being nonzero.
		Testing the homogeneous equation \eqref{eq5.10} against
		$\chi_0(d_1^*z)\psi^k$ and using this identity yields, for
		$w=\sum_{h=0}^Na_h\psi^h$,
		\[
		(2^*-1)\sum_{h=0}^Na_h
		\int_{\mathbb R^N}U^{2^*-2}
		\chi_0(d_1^*z)\psi^h\psi^k\,dz
		=(2^*-1)a_k\int_{\mathbb R^N}U^{2^*-2}
		\chi_0(d_1^*z)|\psi^k|^2\,dz=0,
		\qquad k=0,\ldots,N,
		\]
		the sum collapsing by the symmetry above.  Hence $a_k=0$ for every $k$, and
		$w\equiv0$.
		
		Since $c_{h,n}\to0$, after rescaling we may absorb the finite-dimensional
		terms into a new remainder $g_n$ satisfying
		$\|g_n\|_{H^{-1}(\widetilde M_{3,n})}\to0$.  Put
		\[
		U_n(z)=d_{1,n}^{(N-2)/2}
		(\widetilde P_{\varepsilon_n}U_{
			d_{1,n},\widetilde\xi_{1,n}})
		(d_{1,n}z+\widetilde\xi_{1,n}).
		\]
		It satisfies $0\le U_n\le U_{1,0}$ and converges to $U_{1,0}$ away from
		the shrinking hole.  The rescaled equation is therefore
		\begin{equation}\label{eq5.11}
			-\Delta w_n-(2^*-1)U_n^{2^*-2}w_n=g_n
			\quad\hbox{in }\widetilde M_{3,n}.
		\end{equation}
		
		It remains to show $\|\nabla w_n\|_{L^2(\widetilde M_{3,n})}\to0$; together with
		$c_{h,n}\to0$, this will contradict
		\eqref{eq5.7}.  Testing
		\eqref{eq5.11} with $w_n$ gives
		\begin{equation}\label{eq5.12}
			\int_{\widetilde M_{3,n}} |\nabla w_n|^2 \, \mathrm{d}z
			= (2^{*}-1) \int_{\widetilde M_{3,n}} U_n^{2^{*}-2} w_n^2 \, \mathrm{d}z
			+ \int_{\widetilde M_{3,n}} g_n w_n \, \mathrm{d}z .
		\end{equation}
		The last term is understood as the $H^{-1}$--$H_0^1$ pairing, with
		\[
		|\langle g_n,w_n\rangle|
		\le \|g_n\|_{H^{-1}(\widetilde M_{3,n})}\|\nabla w_n\|_2\longrightarrow0.
		\]
		
		For the potential term, fix $\rho > 0$ and decompose $\widetilde M_{3,n} = (\widetilde M_{3,n} \cap B_\rho(0)) \cup (\widetilde M_{3,n} \setminus B_\rho(0))$.
		On $\widetilde M_{3,n} \cap B_\rho(0)$, $U_n$ is uniformly bounded and $w_n \to 0$ strongly in $L^2(B_\rho(0))$, hence
		\[
		\int_{\widetilde M_{3,n} \cap B_\rho(0)} U_n^{2^{*}-2} w_n^2 \, \mathrm{d}z \to 0 \quad \text{as } n \to \infty.
		\]
		On the complement, since $(2^{*}-2)\frac N2=2^*$ and $0\le U_n\le U_{1,0}$, the
		H\"older's inequality together with Sobolev's inequality
		gives, for every $\rho>0$,
		\[
		\begin{aligned}
			\int_{\widetilde M_{3,n} \setminus B_\rho(0)} U_n^{2^{*}-2} w_n^2 \, \mathrm{d}z
			\le& \|U_n^{2^{*}-2}\|_{L^{N/2}(|z|>\rho)}
			\|w_n\|_{L^{2^*}(\mathbb R^N)}^2\\
			\le& C\|U_{1,0}\|^{2^{*}-2}_{L^{2^*}(|z|>\rho)}
			\|\nabla w_n\|_{L^2(\widetilde M_{3,n})}^2\\
			\le& C\|U_{1,0}\|^{2^{*}-2}_{L^{2^*}(|z|>\rho)}.
		\end{aligned}
		\]  This bound is independent of
		$n$ and tends to $0$ as $\rho\to\infty$.  Hence
		\[
		\limsup_{n\to\infty}\int_{\widetilde M_{3,n}} U_n^{2^{*}-2} w_n^2 \, \mathrm{d}z
		\le C\|U_{1,0}\|^{2^{*}-2}_{L^{2^*}(|z|>\rho)}
		\]
		for every $\rho>0$, and letting $\rho\to\infty$ we obtain
		$\int_{\widetilde M_{3,n}} U_n^{2^{*}-2} w_n^2 \to 0$.
		Consequently, \eqref{eq5.12} forces $\|w_n\|_{H^1_0(\widetilde M_{3,n})}^2 \to 0$.
		Together with $c_{h,n}\to0$, this contradicts
		\eqref{eq5.7} and proves
		\eqref{eq5.5}.  In particular, for the datum $f=L_1v_1$ the tuple is
		$(v_1,0,\ldots,0)$ by uniqueness, so \eqref{eq5.5} reads
		\begin{equation}\label{eq5.13}
			\|\nabla v_1\|_{L^2(M_3)}
			\le C\|L_1v_1\|_{H^{-1}(M_3)},
			\qquad v_1\in K_{M_3}^\perp.
		\end{equation}

		The preceding energy estimate extends by density to $H_0^1$ functions
		with the same zero moments: approximate by smooth functions and correct
		the finitely many moments using $Z^h$.
		For each fixed $R$ and $\varepsilon$, the augmented operator
		\[
		(v,c)\longmapsto
		\left(L_1v-\sum_{h=0}^Nc_hZ^h,
		\left(\int_{M_3}\nabla v\cdot\nabla Z^k\,dy\right)_{k=0}^N\right)
		\]
		is a Fredholm operator of index zero from
		$H_0^1(M_3)\times\mathbb R^{N+1}$ to
		$H^{-1}(M_3)\times\mathbb R^{N+1}$.  Indeed, multiplication by $V_1$ is
		compact relative to the Dirichlet Laplacian, and the remaining terms are
		finite dimensional.  Estimate \eqref{eq5.5} gives
		injectivity, hence bijectivity.  This proves the asserted existence and
		uniqueness.
		
		Now let $f \in L^2(M_3)$ and let $v_1 \in K_{M_3}^\perp$ be the unique $H^1_0$--solution given by the augmented estimate \eqref{eq5.5}.
		From the augmented estimate \eqref{eq5.5} and the continuous embedding $L^2(M_3) \hookrightarrow H^{-1}(M_3)$,
		\begin{equation}\label{eq5.14}
			\|\nabla v_1\|_2+|c|\le C\|f\|_{H^{-1}(M_3)}.
		\end{equation}
		Extending $f$ by zero to $B_{3R}(\widetilde a_1)$ and using that every
		$\varphi\in H_0^1(M_3)$ extends by zero to $H_0^1(B_{3R}(\widetilde a_1))$ with
		$\|\varphi\|_{L^2}\le CR\|\nabla\varphi\|_{L^2}$ (Poincar\'e on a ball of
		radius $3R$), we get
		\[
		\|f\|_{H^{-1}(M_3)}
		=\sup_{\|\nabla\varphi\|_{L^2(M_3)}\le1}\int_{M_3}f\varphi
		\le\|f\|_{L^2(M_3)}\sup_{\|\nabla\varphi\|_{L^2}\le1}\|\varphi\|_{L^2}
		\le CR\|f\|_{L^2(M_3)}.
		\]
		Consequently
		\begin{equation}\label{eq5.15}
			\|\nabla v_1\|_2+|c|\le CR\|f\|_2,
			\qquad \|v_1\|_2\le CR^2\|f\|_2.
		\end{equation}
		Writing $-\Delta v_1=F:=f+V_1v_1+\sum_hc_hZ^h$, the three terms are
		estimated separately: $\|f\|_2$;
		$\|V_1v_1\|_2\le\|V_1\|_\infty\|v_1\|_2\le CR^2\|f\|_2$ by \eqref{eq5.15};
		and $\|\sum_hc_hZ^h\|_2\le C\sum_h|c_h|\le CR\|f\|_2$, the functions $Z^h$ being
		uniformly bounded in $L^2$.    Hence $\|F\|_2\le CR^2\|f\|_2$.
		The regular part of the Dirichlet Green function of $M_3$ is harmonic with positive
		boundary values $|x-y|^{2-N}$, hence positive by the maximum principle, so
		\[
		0\le G_{M_3}(x,y)\le I_2(x-y),
		\]
		uniformly in $R$ and $\varepsilon$, where $I_2$ is the Newtonian potential.  Extending $F$ by zero
		and applying the Hardy--Littlewood--Sobolev inequality to this kernel, which
		belongs to $L^{N/(N-2),\infty}$, gives
		\[
		\|v_1\|_{L^q(M_3)}\le \|I_2\ast F\|_{L^q(\mathbb R^N)}\le C\|F\|_2\le CR^2\|f\|_2,
		\qquad \frac1q=\frac12-\frac2N.
		\]
		Together with \eqref{eq5.15}, this proves
		$\|v_1\|_{X_1}\le CR\|f\|_2$, as asserted in
		\eqref{eq5.6}.
	\end{proof}

	In the sequel $L_1^{-1}f$ denotes the first component of the augmented
	solution of Lemma~\ref{lem5.1}: if $\mathcal{A}$ is the augmented operator
	displayed in the proof above, then $L_1^{-1}f:=[\mathcal{A}^{-1}(f,0)]_1$,
	that is,
	\[
	L_1^{-1}f
	=\bigl(L_1|_{K_{M_3}^\perp}\bigr)^{-1}
	\Bigl(f+\sum_{h=0}^Nc_hZ^h\Bigr),
	\]
	where $c_0,\ldots,c_N$ are the multipliers attached to $f$ by that lemma.
	Thus $L_1^{-1}$ is the projected Dirichlet inverse on
	$K_{M_3}^\perp$, not an inverse of the bare differential
	expression: the restriction alone is invertible only on $\mathrm{Ran}(L_1)$, a
	closed subspace of codimension $N+1$, and the multipliers are exactly what
	carry $f$ onto that range.
	
	\subsubsection{Reduction to a zero boundary condition problem}
	
	Let $\varphi_{1,0}^*$ be the harmonic lifting of the prescribed data, namely the
	unique solution of
	\[
	\begin{cases}
		\Delta \varphi_{1,0}^*=0,&\text{in }M_3,\\[2mm]
		\varphi_{1,0}^*=\varphi_{1,0},&\text{on }\partial M_3.
	\end{cases}
	\]
	
	We seek a solution of the form
	\[
	\varphi_1 = \varphi_{1,0}^* + \tilde{\varphi}_1.
	\]
	
	The weak harmonic equation can be tested by $Z^h\in H_0^1(M_3)$, so
	\begin{equation}\label{eq5.16}
		\int_{M_3}\nabla\varphi_{1,0}^*\cdot\nabla Z^h=0,
		\qquad
		\int_{M_3}\nabla\varphi_1\cdot\nabla Z^h
		=\int_{M_3}\nabla\widetilde\varphi_1\cdot\nabla Z^h.
	\end{equation}
	Thus the full correction and its homogeneous part have the same moments.
	No normal derivative of the $H^1$ lifting is required.
	
	The problem to be solved is the following.
	\begin{equation}\label{eq5.17}
		\begin{cases}
			- \Delta (v_1^*  + \varphi_{1,0}^* + \tilde{\varphi}_1)  = \mu_1 |v_1^* + \varphi_{1,0}^* + \tilde{\varphi}_1|^{2^*- 2}(v_1^*  + \varphi_{1,0}^* + \tilde{\varphi}_1) + \sum_{h=0}^N c_h Z^h, &~\text{in}~ M_3,
			\\[3mm]
			\tilde{\varphi}_1=0, &~\text{on}~ \partial M_3,
			\\[3mm]
			\tilde{\varphi}_1 \in K_{M_3}^\perp.
		\end{cases}
	\end{equation}
	
	Equivalently, writing $L_1=-\Delta-V_1$ and using that the harmonic lifting $\varphi_{1,0}^*$ satisfies $L_1\varphi_{1,0}^*=-V_1\varphi_{1,0}^*$, problem~\eqref{eq5.17} reads
	\begin{equation}\label{eq5.18}
		\begin{cases}
			L_1 \tilde{\varphi}_1 = \widetilde{l}_1 + \mathcal N_1(\varphi_1) + V_1 \varphi_{1,0}^*+\sum_{h=0}^{N}c_h Z^h, & \text{in } M_3, \\[3mm]
			\tilde{\varphi}_1 = 0, & \text{on } \partial M_3,\\[3mm]
			\tilde{\varphi}_1\in K_{M_3}^\perp,
		\end{cases}
	\end{equation}
	where
	\begin{equation}\label{eq5.19}
		\begin{aligned}
			\widetilde{l}_1 :=  \mu_1^{-\frac{1}{2^*-2}}\Big[\big(\widetilde{P}_{\varepsilon} U_{d_1,\widetilde{\xi}_1}\big)^{2^*-1} - U_{d_1,\widetilde{\xi}_1}^{2^*-1}\Big],
		\end{aligned}
	\end{equation}
	and,  
	\begin{equation}\label{eq5.20}
		\begin{aligned}
			\mathcal N_1(\varphi_1) := {}&\mu_1 f\big(v_1^* + \varphi_1\big)
			- \mu_1 f\big(v_1^*\big)
			- \mu_1 f'(v_1^*)\varphi_1,
		\end{aligned}
	\end{equation} 
	with $f(t):=|t|^{2^*-2}t$.

	\begin{Lem}\label{lem5.2}
		Let $\varphi_{1,0}\in H^{1/2}(\partial M_3)\cap L^\infty(\partial M_3)$
		have zero trace on $\partial B_{r_1\sqrt\varepsilon}(\widetilde a_1)$, and let
		$\varphi_{1,0}^*$ be its harmonic lifting, as defined above.  There exists
		$C>0$, independent of $R$ and $\varepsilon$, such that
		\begin{align}
			\|\varphi_{1,0}^*\|_{L^p(M_3)}
			&\le CR^{N/p}\|\varphi_{1,0}\|_{L^\infty(\partial M_3)},
			\qquad 1\le p\le\infty,
			\label{eq5.21}\\
			\|\nabla \varphi_{1,0}^*\|_{L^2(M_3)}
			&\le C\|\varphi_{1,0}\|_{H^{1/2}(\partial M_3)}.
			\label{eq5.22}
		\end{align}
	\end{Lem}
	
	\begin{proof}
		The $H^{1/2}$ assumption gives the unique $H^1$ harmonic lifting with the
		stated boundary values.  By the maximum principle and $|M_3|\le CR^N$,
		Hölder's inequality gives, for every $1\le p\le\infty$,
		\[
		\|\varphi_{1,0}^*\|_{L^p(M_3)}
		\le |M_3|^{1/p}\|\varphi_{1,0}^*\|_{L^\infty(M_3)}
		\le CR^{N/p}\|\varphi_{1,0}\|_{L^\infty(\partial M_3)}.
		\]
		It remains to prove \eqref{eq5.22} with a uniform constant.  After the dilation
		$y=\widetilde a_1+Rz$, the outer boundary is $\partial B_3$ and, for fixed
		large $R$ and small $\varepsilon$, the shrinking hole lies in $B_1$.  The
		trace map of $B_3\setminus\overline{B_2}$ onto
		$H^{1/2}(\partial B_3)\oplus H^{1/2}(\partial B_2)$ is onto, so there is
		$F\in H^1(B_3\setminus\overline{B_2})$ which equals the rescaled datum on
		$\partial B_3$ and vanishes on $\partial B_2$; extending $F$ by zero inside
		$B_2$ keeps it in $H^1(B_3)$, and scaling back gives a competitor $E$ with the
		same boundary values as $\varphi_{1,0}^*$ and
		\[
		\|\nabla E\|_{L^2(M_3)}
		\le C\|\varphi_{1,0}\|_{H^{1/2}(\partial M_3)}.
		\]
		Here $C$ is independent of $R$.  Indeed, under $y=\widetilde a_1+Rz$ the
		rescaled function $F$ lives on the fixed annulus $B_3\setminus\overline{B_2}$.  Scaling back multiplies the Dirichlet integral
		by $R^{N-2}$.  The homogeneous $H^{1/2}$ seminorm on
		$\partial B_{3R}(\widetilde a_1)$ scales in the same way, whereas the $L^2$
		part of the $H^{1/2}$ norm scales as $R^{N-1}$; expressed in the original
		variables, this leaves only the harmless factor $R^{-1}\le1$.  Hence the extension constant is uniform
		in $R$.  Finally, the Dirichlet principle
		yields $\|\nabla \varphi_{1,0}^*\|_{L^2(M_3)}\le\|\nabla E\|_{L^2(M_3)}$, which
		proves \eqref{eq5.22}.
	\end{proof}
	
	\subsubsection{Contraction argument for the inner problem}
	
	Define the mapping $\mathcal{T}: K_{M_3}^\perp \to K_{M_3}^\perp$ by
	\[
	\mathcal{T}(\tilde{\varphi}_1) := L_1^{-1} \Big[ \widetilde{l}_1 + \mathcal N_1(\varphi_{1,0}^* + \tilde{\varphi}_1) + V_1 \varphi_{1,0}^* \Big],
	\]
	where $L_1^{-1}$ is as defined after Lemma~\ref{lem5.1}.  The multipliers
	$c_h$ produced in the same step do
	not enter the definition of $\mathcal{T}$; they are recovered at the fixed
	point, where they are precisely the residuals that the reduction has to
	annihilate.

	We next prove that $\mathcal{T}$ has a unique fixed point in $K_{M_3}^\perp$ for sufficiently small $\varepsilon$.  We first record the estimates used in the contraction argument.

	\begin{Lem}\label{lem5.3}
		For each fixed $\theta\in\bigl(0,\frac{N-4}{N-2}\bigr)$,
		\[
		\|\widetilde l_1\|_{L^2(M_3)}\le C_R\,\varepsilon^{\frac{N-\theta}{4}} .
		\]
	\end{Lem}
	\begin{proof}
		By Lemma~\ref{lem3.2},
		$|\widetilde l_1|\le C\varepsilon^{\frac{N-2}{2}}(G_1+G_2)$ with
		$G_1=(d_1+|y-\widetilde\xi_1|)^{-4}$ and
		$G_2=|y-\widetilde a_1|^{-(N-2)}(d_1+|y-\widetilde\xi_1|)^{-4}$.
		Since $d_1\ge\gamma_*$ and $|y-\widetilde a_1|\ge r_1\sqrt\varepsilon$, for
		$\theta>0$
		\[
		\|G_2\|_{L^2(M_3)}^2
		\le C\varepsilon^{-\frac{N-4+\theta}{2}}
		\int_{M_3}\frac{\mathrm dy}{|y-\widetilde a_1|^{N-\theta}}
		\le C_\theta R^{\theta}\varepsilon^{-\frac{N-4+\theta}{2}},
		\]
		the last integral being $O(R^\theta)$; and $\|G_1\|_{L^2(M_3)}\le C_R$.
		Hence
		\[
		\|\widetilde l_1\|_{L^2(M_3)}
		\le C_R\varepsilon^{\frac{N-2}{2}}
		\Bigl(1+\varepsilon^{-\frac{N-4+\theta}{4}}\Bigr)
		\le C_R\varepsilon^{\frac{N-\theta}{4}},
		\]
		because $\frac{N-2}{2}>\frac{N-\theta}{4}$.
	\end{proof}
	
	\begin{Lem}\label{lem5.4}
		Let $\varphi_1$ satisfy $\|\varphi_1\|_{L^q(M_3)}\le1$.  There exists $C_R>0$,
		depending only on $R$ and $N$, such that
		\[
		\|\mathcal N_1(\varphi_1)\|_{L^2(M_3)}
		\le C_R\begin{cases}
			\|\varphi_{1,0}\|_{L^\infty(\partial M_3)}^2+\|\tilde\varphi_1\|_{L^q(M_3)}^2, & N=5,\\[1mm]
			\|\varphi_{1,0}\|_{L^\infty(\partial M_3)}^{2^*-1}+\|\tilde\varphi_1\|_{L^q(M_3)}^{2^*-1}, & N\ge6.
		\end{cases}
		\]
	\end{Lem}
	
	\begin{proof}
		We use two elementary consequences.  Since $|M_3|\le\omega_N(3R)^N$, for every $1\le t\le q$,
		\begin{equation}\label{eq5.23}
			\|\varphi_1\|_{L^t(M_3)} \le |M_3|^{\frac1t-\frac1q}\|\varphi_1\|_{L^q(M_3)}
			\le C R^{N(\frac1t-\frac1q)}\|\varphi_1\|_{L^q(M_3)},
		\end{equation}
		and, since $\varphi_1=\varphi_{1,0}^*+\tilde\varphi_1$, by \eqref{eq5.21},
		\begin{equation}\label{eq5.24}
			\|\varphi_1\|_{L^q(M_3)}
			\le CR^{\frac{N-4}{2}}\|\varphi_{1,0}\|_{L^\infty(\partial M_3)}
			+\|\tilde\varphi_1\|_{L^q(M_3)}.
		\end{equation}
		Since pointwise remainder obeys
		$|\mathcal N_1(\varphi_1)|\le C\bigl(|\varphi_1|^2+|\varphi_1|^{2^*-1}\bigr)$,  applying \eqref{eq5.23} with $t=4$ and with
		$t=2(2^*-1)=\frac{2N+4}{N-2}<q$,
		\[
		\|\mathcal N_1(\varphi_1)\|_{L^2(M_3)}
		\le C\|\varphi_1\|_{L^4}^2+C\|\varphi_1\|_{L^{2(2^*-1)}}^{2^*-1}
		\le C_R\Bigl(\|\varphi_1\|_{L^q}^2+\|\varphi_1\|_{L^q}^{2^*-1}\Bigr).
		\]
		Since $\|\varphi_1\|_{L^q(M_3)}\le1$, the second term dominates for $N\ge6$ and the
		first for $N=5$; inserting \eqref{eq5.24} therefore gives
		the assertion.
	\end{proof}

	We next establish the Lipschitz estimate for the nonlinear remainder.

	\begin{Lem}\label{lem5.5} 
		There is $C_R>0$
		such that whenever $\varphi,\psi\in L^q(M_3)$ satisfy
		\[
		\max\{\|\varphi\|_{L^q(M_3)}, \|\psi\|_{L^q(M_3)}\} \le1,
		\] there holds 
		\[
		\|\mathcal N_1(\varphi)-\mathcal N_1(\psi)\|_{L^2(M_3)}
		\le
		\begin{cases}
			C_R\bigl(\|\varphi\|_{L^q(M_3)}+ \|\psi\|_{L^q(M_3)}\bigr)
			\|\varphi-\psi\|_{L^q(M_3)},&N=5,\\[1mm]
			C_R (\|\varphi\|_{L^q(M_3)}+\|\psi\|_{L^q(M_3)})^{\frac4{N-2}}
			\|\varphi-\psi\|_{L^q(M_3)},&N\ge6.
		\end{cases}
		\] 
	\end{Lem}
	\begin{proof}
		Throughout this proof all norms are taken over $M_3$.  Denote 
		$r=\|\varphi\|_{L^q}+\|\psi\|_{L^q}$. We first note that
		\begin{equation}\label{eq5.25}
			\|\varphi\|_{L^p}+\|\psi\|_{L^p}\le C_R r,\qquad 1\le p\le q.
		\end{equation}

		\emph{Case $N=5$.}  Since $2^*-3=\frac13>0$ and $v_1^*$ is bounded,  we have  
		\[
		|\mathcal N_1(\varphi)-\mathcal N_1(\psi)|
		\le C\bigl(|\varphi|+|\psi|+|\varphi|^{4/3}+|\psi|^{4/3}\bigr)
		|\varphi-\psi| .
		\]
		Since $\frac25+\frac1{10}=\frac12$ and
		$\||u|^{4/3}\|_{L^{5/2}}=\|u\|_{L^{10/3}}^{4/3}$, Hölder's inequality and
		\eqref{eq5.25} give
		\[
		\|\mathcal N_1(\varphi)-\mathcal N_1(\psi)\|_{L^2(M_3)}
		\le C_R\bigl(r+r^{4/3}\bigr)\|\varphi-\psi\|_{L^{10}(M_3)}.
		\]
		
		\emph{Case $N\ge6$.}  Here $2^*-3\le0$ and $v_1^*\ge0$, so with $a=v_1^*$ the factor
		$(a+|u|+|v|)^{2^*-3}$ is at most $(|u|+|v|)^{2^*-3}$ and the estimate gives
		\[
		|\mathcal N_1(\varphi)-\mathcal N_1(\psi)|
		\le C\bigl(|\varphi|+|\psi|\bigr)^{2^*-2}|\varphi-\psi| .
		\]
		Since $ 
		\frac{N}2(2^*-2)=2^*<q$, H\"older's inequality
		and \eqref{eq5.25} give
		\[
		\|\mathcal N_1(\varphi)-\mathcal N_1(\psi)\|_{L^2(M_3)}
		\le C_R r^{2^*-2}\|\varphi-\psi\|_{L^q(M_3)}
		=C_R r^{\frac4{N-2}}\|\varphi-\psi\|_{L^q(M_3)}.
		\]
	\end{proof}

	\begin{Lem}\label{lem5.6}
		There exists $C>0$, independent of $R$
		and $\varepsilon$, such that
		\[
		\|V_1\varphi_{1,0}^*\|_{L^2(M_3)}\le C(1+R^{\frac{N-8}{2}}(\ln R)^{4/N})\|\varphi_{1,0}\|_{L^\infty(\partial M_3)}.
		\]
	\end{Lem}
	\begin{proof}
		Take $m=1$ for $N\le7$ and $m=N/8$ for $N\ge8$, and set
		$m':=m/(m-1)\in(1,\infty]$.  Since $V_1\leq c\,U_{d_1,\widetilde\xi_1}^{2^*-2}$, 
		\[
		\int_{M_3}V_1^{2m}\le C+ C\int_1^{3R}r^{N-1-8m}\,dr
		\le\begin{cases}C,&N\le7,\\ C\ln R,&N\ge8,\end{cases}
		\]
		while \eqref{eq5.21} with $p=2m'$ gives
		$\|\varphi_{1,0}^*\|_{L^{2m'}}\le CR^{N/(2m')}\|\varphi_{1,0}\|_{L^\infty(\partial M_3)}$.
		Hölder's inequality therefore yields, for $N\ge8$,
		\[
		\|V_1\varphi_{1,0}^*\|_{L^2(M_3)}^2
		\le\|V_1\|_{L^{2m}}^2\|\varphi_{1,0}^*\|_{L^{2m'}}^2
		\le C(\ln R)^{8/N}R^{N-8}\|\varphi_{1,0}\|_{L^\infty(\partial M_3)}^2.
		\]
		For $N\le7$ the same inequality with $m=1$
		and $\int V_1^2\le C$ gives
		$\|V_1\varphi_{1,0}^*\|_{L^2}\le C\|\varphi_{1,0}\|_{L^\infty(\partial M_3)}$.
	\end{proof}
	
	We consider boundary data $\varphi_{1,0}$ in
	$H^{1/2}(\partial M_3)\cap L^\infty(\partial M_3)$ vanishing on
	$\partial B_{r_1\sqrt\varepsilon}(\widetilde a_1)$ and satisfying
	\begin{equation}\label{eq5.26}
		\|\varphi_{1,0}\|_{L^\infty(\partial M_3)}
		\le C_{\mathrm b}R^{-N/2}\varepsilon^{(N-2\theta)/4},
		\qquad
		\|\varphi_{1,0}\|_{H^{1/2}(\partial M_3)}
		\le C_{\mathrm b}R^{-1/2}\varepsilon^{(N-2\theta)/4}.
	\end{equation}
	Here $C_{\mathrm b}>0$ is a fixed constant, larger than twice the constant
	$C$ of Corollary~\ref{cor4.4}, and independent of $R$ and $\varepsilon$.  We write
	\[
	\mathcal D_{1,\varepsilon}:=
	\Bigl\{\varphi_{1,0}\in K_1^\perp:\ \varphi_{1,0}|_{\partial M_3}\in
	H^{1/2}(\partial M_3)\cap L^\infty(\partial M_3),\quad
	\varphi_{1,0}=0\ \text{on }\partial B_{r_1\sqrt\varepsilon}(\widetilde a_1),\quad
	\eqref{eq5.26}\ \text{holds}\Bigr\}
	\]
	for the corresponding class of admissible data, and let
	$\mathcal D_{2,\varepsilon}$ be defined analogously, with $Q_3$ and
	$\varphi_{2,0}$ in place of $M_3$ and $\varphi_{1,0}$.
	
	\begin{Prop}\label{prop5.7}
		Let $\varphi_{1,0}\in\mathcal D_{1,\varepsilon}$.  For $R$
		sufficiently large and then $\varepsilon$ sufficiently small,
		problem~\eqref{eq5.18} has a unique solution
		\(\widetilde\varphi_1\) in the class
		\[
		\left\{
		\xi\in K_{M_3}^\perp:
		\|\xi\|_{L^q(M_3)}\le a_R
		\right\},
		\]
		for some \(a_R>0\) depending only on \(R\) and \(N\).
		Moreover, the solution satisfies
		\begin{equation}\label{eq5.27}
			\|\tilde{\varphi}_1\|_{X_1} 
			\le C R^{-3/2} \varepsilon^{\frac{N-2\theta}{4}},
		\end{equation}
		where $C>0$ is independent of $R$ and $\varepsilon$.
	\end{Prop}
	
	\begin{proof}
		For $\varphi=\varphi_{1,0}^*+\tilde \varphi$ with $\tilde\varphi\in K_{M_3}^\perp$,  
		assume first that $\| {\varphi}_1\|_{L^q(M_3)} \leq 1$.
		By Lemmas \ref{lem5.1}, \ref{lem5.3} and \ref{lem5.4} and the bounds defining $\mathcal D_{1,\varepsilon}$, we have
		
		\begin{equation}\label{eq5.28}
			\begin{aligned}
				\|\mathcal{T}(\tilde{\varphi})\|_{X_1} &\le C R \Bigl( \|\widetilde{l}_1\|_{L^2(M_3)} + \|\mathcal N_1(\varphi_{1,0}^* + \tilde{\varphi})\|_{L^2(M_3)} + \|V_1 \varphi_{1,0}^*\|_{L^2(M_3)} \Bigr)\\
				&\le    C_R \bigl(\varepsilon^{\frac{N-\theta}{4}} +\|\varphi_{1,0}\|_{L^\infty(\partial M_3)}^{\beta} + \|\tilde{\varphi}\|_{L^q(M_3)}^\beta\bigr) + C(1+R^{\frac{N-8}{2}}(\ln R)^{4/N})\|\varphi_{1,0}\|_{L^\infty(\partial M_3)}\\
				&\le   C_R \|\tilde{\varphi}\|_{L^q(M_3)}^\beta +C (R^{1-N/2}+R^{-3}(\ln R)^{4/N})\varepsilon^{\frac{N-2\theta}{4}}\\
				&\le C_R \|\tilde{\varphi}\|_{L^q(M_3)}^\beta + C R^{-3/2} \varepsilon^{\frac{N-2\theta}{4}},
			\end{aligned}
		\end{equation}
		where $\beta=\min\set{2, 2^*-1}>1$, $C_R>1$ is a constant depending only on $R$ and $N$, and $C>0$ is a constant independent of $R$ and $\varepsilon$.   
		Consider the set 
		\[
		\mathcal B= \Set{\xi\in K_{M_3}^\perp | \|\xi\|_{L^q(M_3)}\leq (2R^3C_R)^{\frac{1}{1-\beta}}},
		\]
		where the constant $C_R$ is the largest  one appearing in   Lemmas \ref{lem5.3}, \ref{lem5.4}, \ref{lem5.5} and \eqref{eq5.28}.   
		
		For $\tilde\varphi\in\mathcal B$, by \eqref{eq5.24} and $\varphi_{1,0}\in\mathcal D_{1,\varepsilon}$,
		\[
		\|\varphi\|_{L^q(M_3)}\leq CR^{-2}\varepsilon^{\frac{N-2\theta}{4}}+(2R^3C_R)^{\frac{1}{1-\beta}}.
		\]
		Hence $\|\varphi\|_{L^q(M_3)}<(R^3C_R)^{\frac{1}{1-\beta}}\le1$ for $R$ large and
		$\varepsilon$ small, and
		\[
		\|\mathcal T(\tilde\varphi)\|_{L^q(M_3)}\leq \frac1{2R^2} (2R^3C_R)^\frac{1}{1-\beta} +C R^{-1/2} \varepsilon^{\frac{N-2\theta}{4}}.
		\]
		Therefore, for sufficiently small $\varepsilon$, $\mathcal T(\tilde\varphi)\in \mathcal B$.
		
		Now let $\tilde\varphi,\tilde\psi\in\mathcal B$ and denote
		\[
		\varphi=\varphi_{1,0}^*+\tilde\varphi,
		\qquad
		\psi=\varphi_{1,0}^*+\tilde\psi.
		\]
		Lemmas
		\ref{lem5.5} and \ref{lem5.1} give
		\[\begin{aligned}
			\|\mathcal T(\tilde\varphi)-\mathcal T(\tilde\psi)\|_{X_1}
			\le & CRC_R(\|\varphi\|_{L^q(M_3)}+\|\psi\|_{L^q(M_3)})^{\beta-1}
			\|\tilde\varphi-\tilde\psi\|_{L^q(M_3)}\\
			\le & CR^{-2} 
			\|\tilde\varphi-\tilde\psi\|_{L^q(M_3)} \leq CR^{-1} \|\tilde\varphi-\tilde\psi\|_{X_1}.
		\end{aligned}
		\]
		Thus $\mathcal T$ is a contraction on $\mathcal B$ for large $R$.
		
		At the fixed point $\tilde\varphi_1\in\mathcal B$, \eqref{eq5.28} gives
		\[
		\|\tilde\varphi_1\|_{X_1}
		\le  \frac1{2R^2} \|\tilde\varphi_1\|_{L^q(M_3)} +C R^{-1/2} \varepsilon^{\frac{N-2\theta}{4}}\le  \frac1{2R} \|\tilde\varphi_1\|_{X_1} +C R^{-1/2} \varepsilon^{\frac{N-2\theta}{4}}.
		\]
		This completes the proof.
	\end{proof}
	
	\medskip
	
	\subsection{The inner solution operator}\label{subsec5.2}
	
	For $\varphi_{1,0}\in\mathcal D_{1,\varepsilon}$, let $\varphi_{1,0}^*$ be its
	harmonic lifting in $M_3$ and let $\tilde\varphi_1\in K_{M_3}^\perp$ be the unique
	solution of \eqref{eq5.18} given by Proposition~\ref{prop5.7}. We define \(T_1(\varphi_{1,0})\) by gluing
	\(\varphi_{1,0}\) outside \(M_3\) with
	\(\varphi_{1,0}^*+\widetilde\varphi_1\) inside \(M_3\). The two pieces agree
	on \(\partial M_3\), so \(T_1(\varphi_{1,0})\) is precisely the solution
	\(\varphi_1\) of \eqref{eq5.1}. Since the local moment conditions imply the
	global orthogonality conditions, one obtains
	\[
	T_1:\mathcal D_{1,\varepsilon}\longrightarrow K_1^\perp.
	\]
	\begin{Rem}\label{rem5.8}
		By elliptic regularity applied to \eqref{eq5.1}, each such
		$\varphi_1$ is of class $C^1$ in the interior of $M_3$, so its restrictions
		to spheres such as $\partial B_{R_0}(\widetilde a_1)$ are ordinary
		restrictions of a continuous function, and the $L^\infty$ norms used
		below on such spheres are meaningful.
	\end{Rem}
	Define \(T_2\) analogously and put \(T:=T_1\times T_2\).
	
	\subsection{Estimates for the inner problem}
	
	\begin{Lem}\label{lem5.9}
		Let $R$ be sufficiently large and $\varepsilon$ sufficiently small, and
		let $\varphi_1=T_1(\varphi_{1,0})$ and $\bar\varphi_1=T_1(\bar\varphi_{1,0})$
		be inner solutions with data
		$\varphi_{1,0},\bar\varphi_{1,0}\in\mathcal D_{1,\varepsilon}$.  Then
		\begin{gather}
			\|\nabla\varphi_1\|_{L^2(M_3)}+\|\varphi_1\|_{L^q(M_3)}
			\le CR^{-\frac12}\varepsilon^{\frac{N-2\theta}{4}},
			\label{eq5.29}\\
			\|\varphi_1-\bar\varphi_1\|_{L^q(M_3)}
			\le C(R^2+R^{\frac{N-4}{2}}(\ln R)^{4/N})
			\|\varphi_{1,0}-\bar\varphi_{1,0}\|_{L^\infty(\partial M_3)},
			\label{eq5.30}
		\end{gather}
		where $C>0$ is independent of $R$ and $\varepsilon$.
	\end{Lem}
	\begin{proof}
		We write $\varphi_1=\varphi_{1,0}^*+\tilde\varphi_1$ and
		$\bar\varphi_1=\bar\varphi_{1,0}^*+\tilde\psi_1$, where
		$\varphi_{1,0}^*,\bar\varphi_{1,0}^*$ are the harmonic liftings of the
		prescribed data and $\tilde\varphi_1,\tilde\psi_1\in K_{M_3}^\perp$ are
		the parts with zero boundary data.
		Since $\varphi_{1,0}^*$ is harmonic and $\tilde\varphi_1\in H_0^1(M_3)$,
		the cross term vanishes,
		$\int_{M_3}\nabla\varphi_{1,0}^*\cdot\nabla\tilde\varphi_1\,dy=0$, and
		therefore the Dirichlet energies split exactly,
		\[
		\|\nabla\varphi_1\|_{L^2(M_3)}^2
		=\|\nabla\varphi_{1,0}^*\|_{L^2(M_3)}^2
		+\|\nabla\tilde\varphi_1\|_{L^2(M_3)}^2 .
		\]
		Inserting \eqref{eq5.22}, \eqref{eq5.27} and the bounds defining $\mathcal D_{1,\varepsilon}$ and taking
		square roots gives
		$\|\nabla\varphi_1\|_{L^2(M_3)}\le CR^{-\frac12}\varepsilon^{\frac{N-2\theta}{4}}$.
		
		By the triangle inequality, \eqref{eq5.21} with $p=q=\frac{2N}{N-4}$ and
		$\|u\|_{L^q}\le R\|u\|_{X_1}$,
		\[
		\|\varphi_1\|_{L^q(M_3)}
		\le\|\varphi_{1,0}^*\|_{L^q(M_3)}+\|\tilde\varphi_1\|_{L^q(M_3)}
		\le CR^{\frac{N-4}{2}}\|\varphi_{1,0}\|_{L^\infty(\partial M_3)}
		+R\|\tilde\varphi_1\|_{X_1},
		\]
		\eqref{eq5.27} bounds the last term by $CR^{-\frac12}\varepsilon^{\frac{N-2\theta}{4}}$,
		while the bound for $\|\varphi_{1,0}\|_{L^\infty(\partial M_3)}$ in
		$\mathcal D_{1,\varepsilon}$ turns the first term into
		$CR^{-2}\varepsilon^{\frac{N-2\theta}{4}}$, which is absorbed.
		Together with the preceding estimate this is \eqref{eq5.29}.

		Since
		\[\tilde\varphi_1-\tilde\psi_1=
		\mathcal T(\tilde\varphi_1)-\mathcal T(\tilde\psi_1)
		=L_1^{-1}(\mathcal N_1(\varphi_1)-\mathcal N_1(\bar\varphi_1)
		+V_1\bigl(\varphi_{1,0}^*-\bar\varphi_{1,0}^*\bigr)),
		\]
		and  by Lemmas~\ref{lem5.1}, \ref{lem5.5}, \ref{lem5.6}, and \eqref{eq5.29}, we have
		\[\begin{aligned}
			\|\tilde\varphi_1-\tilde\psi_1\|_{L^q(M_3)}
			\le& CR^2\bigl\|\mathcal N_1(\varphi_1)
			-\mathcal N_1(\bar\varphi_1)\bigr\|_{L^2(M_3)}
			+CR^2\bigl\|V_1\bigl(\varphi_{1,0}^*-\bar\varphi_{1,0}^*\bigr)
			\bigr\|_{L^2(M_3)} \\
			\le & C_R\varepsilon^{(\beta-1)\frac{N-2\theta}{4}}\|\tilde\varphi_1-\tilde\psi_1\|_{L^q(M_3)}
			+C(R^2+R^{\frac{N-4}{2}}(\ln R)^{4/N})\|\varphi_{1,0}-\bar\varphi_{1,0}\|_{L^\infty(\partial M_3)}.
		\end{aligned}
		\]
		Since $\varphi_{1,0}\mapsto\varphi_{1,0}^*$ is linear, 
		$\varphi_{1,0}^*-\bar\varphi_{1,0}^*=(\varphi_{1,0}-\bar\varphi_{1,0})^*$.
		Applying  Lemma~\ref{lem5.2} to  
		$\varphi_{1,0}-\bar\varphi_{1,0}$ gives,
		\[ 
		\bigl\|\varphi_{1,0}^*-\bar\varphi_{1,0}^*\bigr\|_{L^q(M_3)}
		\le CR^{\frac{N-4}{2}}\|\varphi_{1,0}-\bar\varphi_{1,0}\|_{L^\infty(\partial M_3)} .
		\]
		Thus, for sufficiently small $\varepsilon$,
		\begin{align*}
			\|\varphi_1-\bar\varphi_1\|_{L^q(M_3)}
			&\le\bigl\|\varphi_{1,0}^*-\bar\varphi_{1,0}^*\bigr\|_{L^q(M_3)}
			+\|\tilde\varphi_1-\tilde\psi_1\|_{L^q(M_3)}\\
			&\le \frac12 \|\varphi_1-\bar\varphi_1\|_{L^q(M_3)}+ C(R^2+R^{\frac{N-4}{2}}(\ln R)^{4/N})\|\varphi_{1,0}-\bar\varphi_{1,0}\|_{L^\infty(\partial M_3)} ,
		\end{align*}
		which implies \eqref{eq5.30}.
	\end{proof}
	
	\begin{Lem}\label{lem5.10}
		Let $\varphi_1\in T_1(\mathcal D_{1,\varepsilon})$.  Then
		\begin{equation}\label{eq5.31}
			\|\varphi_1\|_{L^\infty(\partial B_{R_0}(\widetilde a_1))}
			\le CR^{-\frac12}\varepsilon^{\frac{N-2\theta}{4}}.
		\end{equation}
		Here $C$ is independent of $\varepsilon$ and $R$.
	\end{Lem}
	\begin{proof}
		For $y_0\in\partial B_{R_0}(\widetilde a_1)$ put $B:=B_1(y_0)$.  Since
		$R_0+1<3R$ and $B$ is far from the small hole, $B\subset M_3$; and since
		$\operatorname{supp}Z^h\subset B_3(\widetilde a_1)\subset B_{R_0-1}(\widetilde a_1)$,
		the multipliers vanish on $B$.  By the mean value theorem applied to
		$f(t)=|t|^{2^*-2}t$, the nonlinear term can be written as $a\varphi_1$ with
		$a=\mu_1\int_0^1f'(v_1^*+t\varphi_1)\,\mathrm dt$; since
		$f'(u)=(2^*-1)|u|^{4/(N-2)}$ and $v_1^*$ is bounded, this gives
		$|a|\le C(1+|\varphi_1|^{4/(N-2)})$.  Hence $-\Delta\varphi_1=a\varphi_1+\widetilde l_1$
		in $B$, and
		$\|\varphi_1\|_{L^q(M_3)}\le CR^{-1/2}\varepsilon^{\frac{N-2\theta}{4}}\le1$
		by \eqref{eq5.29}, so $\|a\|_{L^s(B)}\le C$ for $s=q(N-2)/4>N/2$.  By Lemma~\ref{lem3.2},
		$|\widetilde l_1|\le C\varepsilon^{\frac{N-2}{2}}(G_1+G_2)$ with $G_1,G_2$ as in the proof of
		Lemma~\ref{lem5.3}; since $|y-\widetilde\xi_1|\simeq R_0$ on $B$, both are bounded
		there, and hence $|\widetilde l_1|\le C\varepsilon^{\frac{N-2}{2}}$ on $B$.  The interior boundedness
		estimate \cite[Theorem 8.17]{GT01} on $B$ therefore gives
		$|\varphi_1(y_0)|\le C(\|\varphi_1\|_{L^q(M_3)}+\varepsilon^{\frac{N-2}{2}})$
		with $C$ independent of $y_0$.  Taking the supremum over
		$y_0\in\partial B_{R_0}(\widetilde a_1)$ and using \eqref{eq5.29} together with
		$\varepsilon^{\frac{N-2}{2}}\le R^{-1/2}\varepsilon^{\frac{N-2\theta}{4}}$ for
		small $\varepsilon$ gives \eqref{eq5.31}.
	\end{proof}

	\begin{Lem}\label{lem5.11}
		Let $\varphi_1,\bar\varphi_1\in T_1(\mathcal D_{1,\varepsilon})$.  Then
		\begin{equation}\label{eq5.32}
			\|\varphi_1-\bar\varphi_1\|_{L^\infty(
				\partial B_{R_0}(\widetilde a_1))}
			\le C(R^2+R^{\frac{N-4}{2}}(\ln R)^{4/N})
			\|\varphi_{1}-\bar\varphi_{1}\|_{L^\infty(\partial M_3)},
		\end{equation}
	\end{Lem}
	\begin{proof}
		For $y_0\in\partial B_{R_0}(\widetilde a_1)$ and $B:=B_1(y_0)$ as in the proof of
		Lemma~\ref{lem5.10}, the common inhomogeneity $\widetilde l_1$ cancels in the
		difference and, since $\operatorname{supp}Z^h\cap B=\varnothing$, the
		multipliers vanish.  Thus
		\begin{equation}\label{eq5.33}
			-\Delta(\varphi_1-\bar\varphi_1)=a(\varphi_1-\bar\varphi_1),
			\qquad |a|\le C(1+|\varphi_1|^{4/(N-2)}
			+|\bar\varphi_1|^{4/(N-2)}).
		\end{equation}
		Here $a$ is the average of $f'$ along the segment joining
		$v_1^*+\varphi_1$ and $v_1^*+\bar\varphi_1$.
		By \eqref{eq5.29} both $L^q$ norms are at most
		$CR^{-1/2}\varepsilon^{\frac{N-2\theta}{4}}\le1$, so the coefficient has a
		uniform $L^s$ bound for $s=q(N-2)/4>N/2$.  The same local boundedness
		estimate gives
		$\|\varphi_1-\bar\varphi_1\|_{L^\infty(\partial B_{R_0}(\widetilde a_1))}
		\le C\|\varphi_1-\bar\varphi_1\|_{L^q(M_3)}$, and \eqref{eq5.30} bounds the
		right-hand side by $C(R^2+R^{\frac{N-4}{2}}(\ln R)^{4/N})
		\|\varphi_1-\bar\varphi_1\|_{L^\infty(\partial M_3)}$, which is \eqref{eq5.32}.
	\end{proof}
	
	\section{Matching the inner and outer solutions}\label{sec6}

	The outer correction $\varphi_1=v_1-v_1^*$, restricted to the outer
	component of $\partial M_3$, supplies the boundary datum of the inner
	problem in $M_3$, and similarly for $\partial Q_3$ and $Q_3$.  Conversely,
	the restriction of the inner solution $\varphi_1$ to $M_0$ supplies the new
	prescribed datum of the outer problem, and similarly for $Q_0$.  At a fixed point
	the two prescriptions are compatible, and the inner and outer solutions
	agree on the overlap regions $M_3\setminus M_0$ and $Q_3\setminus Q_0$.

	For $\Phi=(\varphi_1,\varphi_2)$ and
	$\bar\Phi=(\bar\varphi_1,\bar\varphi_2)$ in $S(\Lambda_\varepsilon)$, define
	\[
	\mathrm{dist}(\Phi,\bar\Phi):=
	\|\varphi_1-\bar\varphi_1\|_{L^q(M_0)}
	+\|\varphi_1-\bar\varphi_1\|_{L^\infty(\partial M_0)}
	+\|\varphi_2-\bar\varphi_2\|_{L^q(Q_0)}
	+\|\varphi_2-\bar\varphi_2\|_{L^\infty(\partial Q_0)} .
	\]
	By Remark~\ref{rem5.8}, the boundary norms in the definition of
	$\mathrm{dist}$ are ordinary suprema of continuous traces.  The function
	$\mathrm{dist}$ so defined is a metric on $S(\Lambda_\varepsilon)$.  Since for every
	$H=(H_1,H_2)\in\Lambda_\varepsilon$ the first component of $S(H)$ equals
	$H_1$ on $M_0$ and the second equals $H_2$ on $Q_0$, writing $\Phi=S(H)$ and
	$\bar\Phi=S(\bar H)$ with $H,\bar H\in\Lambda_\varepsilon$ and using the
	bounds in the definition of $\Lambda_\varepsilon$, both the $L^q$ and the
	$L^\infty$ norms above are at most $2\varepsilon^{\frac{N-2\theta}{4}}$
	(up to the factor $|M_0|^{1/q}$ in the $L^q$ case), so the distance is
	finite.  Moreover $\mathrm{dist}(\Phi,\bar\Phi)=0$ forces
	$\varphi_1=\bar\varphi_1$ almost everywhere on $M_0$ and
	$\varphi_2=\bar\varphi_2$ almost everywhere on $Q_0$, that is,
	$H_1=\bar H_1$ and $H_2=\bar H_2$ on these regions; by outer uniqueness
	this forces $\Phi=\bar\Phi$.
	
	The space $S(\Lambda_\varepsilon)$ is complete under this distance.  Let $\Phi_n$ be a
	$\mathrm{dist}$-Cauchy sequence.  Its components converge in $L^q$ on $M_0$ and $Q_0$ and uniformly on
	$\partial M_0$ and $\partial Q_0$, to limits $\varphi_1$ and $\varphi_2$.  The components of $\Phi_n$
	obey the bounds defining $\Lambda_{\varepsilon,1}$ and $\Lambda_{\varepsilon,2}$, and these pass to the limits: the $L^\infty$ bound on an almost-everywhere
	convergent subsequence, and the gradient bounds by weak closedness of closed sets in
	$H^1$, the $L^q$ convergence
	ensuring that the weak limits are $\varphi_1$ and $\varphi_2$. Because
	$\operatorname{supp}Z^h\subset M_0$ and $\operatorname{supp}Y^h\subset
	Q_0$, the moment conditions involve only $\varphi_1|_{M_0}$ and $\varphi_2|_{Q_0}$, so these conditions
	pass to the limit.  Extending $\varphi_1$ and $\varphi_2$ outside
	$M_0$ and $Q_0$ by standard $H^1$ extensions vanishing on the small-hole
	spheres gives an element $H\in\Lambda_\varepsilon$, and since the outer
	solution $S(H)$ restricts to $\varphi_1$ on $M_0$ and to $\varphi_2$ on
	$Q_0$, $\Phi_n\to S(H)$ in $\mathrm{dist}$.
	
	We apply the contraction principle to the composition
	$\Phi\mapsto S(T(\Phi))$ on $S(\Lambda_\varepsilon)$.
	
	\begin{Lem}\label{lem6.1}
		For $R$ sufficiently large and $\varepsilon$ sufficiently small, every
		$\Phi=(\varphi_1,\varphi_2)\in S(\Lambda_\varepsilon)$ satisfies
		$\varphi_1|_{\partial M_3}\in\mathcal D_{1,\varepsilon}$ and
		$\varphi_2|_{\partial Q_3}\in\mathcal D_{2,\varepsilon}$, and
		$T(\Phi)\in\Lambda_\varepsilon$.  Consequently the composition
		$\Phi\mapsto S(T(\Phi))$ maps $S(\Lambda_\varepsilon)$ into itself.
	\end{Lem}
	
	\begin{proof}
		Let $\Phi=(\varphi_1,\varphi_2)\in S(\Lambda_\varepsilon)$ and let
		$\varphi_{1,0}$ be the trace of $\varphi_1$ on $\partial M_3$, which
		vanishes on $\partial B_{r_1\sqrt\varepsilon}(\widetilde a_1)$.  By the ``Moreover'' part
		of Corollary~\ref{cor4.4}, $\varphi_{1,0}$ satisfies both bounds in
		\eqref{eq5.26} with the constant $C$ of that corollary; since
		$C_{\mathrm b}$ dominates $C$, this gives
		$\varphi_{1,0}\in\mathcal D_{1,\varepsilon}$, and the second component gives
		$\varphi_2|_{\partial Q_3}\in\mathcal D_{2,\varepsilon}$ in the same way.
		
		Let $(\psi_1,\psi_2):=T(\Phi)$.  Since $T$ produces global $H_0^1$
		functions and $Z^h$, $Y^h$ are supported in $M_0\cup Q_0$, the local moment
		conditions give membership in $K_i^\perp$.  By \eqref{eq5.29},
		$\|\nabla\psi_1\|_{L^2(M_0)}\le\|\nabla\psi_1\|_{L^2(M_3)}
		\le CR^{-\frac12}\varepsilon^{\frac{N-2\theta}{4}}$, and
		Lemma~\ref{lem5.10} together with
		$\varepsilon^{\frac{N-2}{2}}\le R^{-\frac12}\varepsilon^{\frac{N-2\theta}{4}}$
		and $\psi_1=0$ on $\partial B_{r_1\sqrt\varepsilon}(\widetilde a_1)$ gives
		\begin{equation}\label{eq6.1}
			\|\nabla\psi_1\|_{L^2(M_0)}+\|\psi_1\|_{L^\infty(\partial M_0)}
			\le CR^{-\frac12}\varepsilon^{\frac{N-2\theta}{4}} .
		\end{equation}
		The same holds for the second component.  Choosing $R$ large so that
		$CR^{-1/2}\le1$ therefore gives $T(\Phi)\in\Lambda_\varepsilon$ and hence
		$S(T(\Phi))\in S(\Lambda_\varepsilon)$.
	\end{proof}
	
	\begin{Lem}\label{lem6.2}
		The composition $\Phi\mapsto S(T(\Phi))$ is a contraction on
		$S(\Lambda_\varepsilon)$, namely
		\begin{equation*}
			\mathrm{dist}\big(S(T(\Phi)), S(T(\bar\Phi))\big)\leq \frac{1}{2}\,\mathrm{dist}(\Phi,\bar\Phi).
		\end{equation*}
	\end{Lem}
	
	\begin{proof}
		Let $\Phi=(\varphi_1,\varphi_2),\bar\Phi=(\bar\varphi_1,\bar\varphi_2)
		\in S(\Lambda_\varepsilon)$ and set $w=\varphi_1-\bar\varphi_1$.
		
		When $w$ appears in a boundary norm it denotes the trace of $\varphi_1-\bar\varphi_1$, which vanishes on
		$\partial B_{r_1\sqrt\varepsilon}(\widetilde a_1)$.  By Lemma~\ref{lem4.8},
		\begin{equation}\label{eq6.2}
			\|w\|_{L^\infty(\partial B_{3R}(\widetilde a_1))}
			\le CR^{-(N-2)}\mathrm{dist}(\Phi,\bar\Phi).
		\end{equation}
		Let $\psi_1=T_1(\varphi_1)$ and $\bar\psi_1=T_1(\bar\varphi_1)$, so that
		$\psi_1,\bar\psi_1\in T_1(\mathcal D_{1,\varepsilon})$, and put
		$u:=\psi_1-\bar\psi_1$.  By \eqref{eq5.32}, \eqref{eq5.30} and \eqref{eq6.2},
		\[
		\|u\|_{L^\infty(\partial B_{R_0}(\widetilde a_1))}+\|u\|_{L^q(M_3)}
		\le C\Bigl(R^2+R^{\frac{N-4}{2}}(\ln R)^{4/N}\Bigr)\|w\|_{L^\infty(\partial M_3)}
		\le CR^{-1}\mathrm{dist}(\Phi,\bar\Phi).
		\]
		Since $u=0$ on $\partial B_{r_1\sqrt\varepsilon}(\widetilde a_1)$ and
		$M_0\subset M_3$, this controls both norms entering $\mathrm{dist}$:
		\begin{equation}\label{eq6.3}
			\|u\|_{L^q(M_0)}+\|u\|_{L^\infty(\partial M_0)}
			\le CR^{-1}\mathrm{dist}(\Phi,\bar\Phi).
		\end{equation}
		The same estimates hold for the second component on $Q_0$.  Since $S$
		preserves the prescribed data on $M_0$ and $Q_0$, \eqref{eq6.3} and its
		analogue for the second component give
		\[
		\mathrm{dist}\bigl(S(T(\Phi)),S(T(\bar\Phi))\bigr)
		\le CR^{-1}\mathrm{dist}(\Phi,\bar\Phi).
		\]
		Choosing $R$ sufficiently large so that $CR^{-1}\le1/2$ proves the
		assertion.
	\end{proof}

	\smallskip
	Lemma \ref{lem6.2} implies that $\Phi\mapsto S(T(\Phi))$ has a unique fixed point
	in $S(\Lambda_\varepsilon)$. To see it is a solution, we consider the following problem 
	\begin{equation}\label{eq6.4}
		\left\{
		\begin{aligned}
			-\Delta(v_1^*+\phi_1)
			={}&\mu_1|v_1^*+\phi_1|^{2^*-2}(v_1^*+\phi_1)\\
			&+\beta |v_1^*+\phi_1|^{\frac{2^*}{2}-2}(v_1^*+\phi_1)
			|v_2^*+\phi_2|^{\frac{2^*}{2}}
			+\sum_{h=0}^N c_hZ^h
			&&\text{in }\widetilde\Omega_\varepsilon,\\[1mm]
			-\Delta(v_2^*+\phi_2)
			={}&\mu_2|v_2^*+\phi_2|^{2^*-2}(v_2^*+\phi_2)\\
			&+\beta |v_2^*+\phi_2|^{\frac{2^*}{2}-2}(v_2^*+\phi_2)
			|v_1^*+\phi_1|^{\frac{2^*}{2}}
			+\sum_{h=0}^N\lambda_hY^h
			&&\text{in }\widetilde\Omega_\varepsilon,\\[1mm]
			&(\phi_1,\phi_2)\in
			\mathbf K_{\mathbf d,\boldsymbol\tau,\varepsilon}^\perp.
		\end{aligned}
		\right.
	\end{equation}

	\begin{Lem}\label{lem6.3}
		A pair $\Phi=(\phi_1,\phi_2)\in S(\Lambda_\varepsilon)$ satisfies
		$S(T(\Phi))=\Phi$ if and only if $\Phi\in\Lambda_\varepsilon$ and $\Phi$
		solves~\eqref{eq6.4} with
		\[
		(2^*-1)\mu_1\|\phi_1+v_1^*\|_{L^{2^*}(\widetilde{\Omega}_{\varepsilon}\setminus(M_0\cup Q_0))}^{2^*-2}<\frac{\mathcal S}{2},
		\qquad
		(2^*-1)\mu_2\|\phi_2+v_2^*\|_{L^{2^*}(\widetilde{\Omega}_{\varepsilon}\setminus(M_0\cup Q_0))}^{2^*-2}<\frac{\mathcal S}{2}.
		\]
		In particular, the fixed point belongs to
		$S(\Lambda_\varepsilon)\cap\Lambda_\varepsilon$.
	\end{Lem}
	
	\begin{proof}
		(i)	Assume $S(T(\phi_1,\phi_2))=(\phi_1,\phi_2)$.  We denote
		\begin{align*} 
			\zeta_{1}:=T_1 \phi_1,\quad
			\zeta_{2}:=T_2 \phi_2
			. 
		\end{align*}
		
		By definition of the global replacement operator,
		$\zeta_1=\phi_1$ in
		$\widetilde\Omega_\varepsilon\setminus M_3$.  On the other hand, the fixed
		point identity $\Phi=S(T(\Phi))$ and the fact that $S$ preserves its input on
		$M_0$ give
		\begin{equation}\label{eq6.5}
			\zeta_1=\phi_1\quad\text{in }M_0,
			\qquad
			\zeta_1=\phi_1\quad\text{in }
			\widetilde\Omega_\varepsilon\setminus M_3.
		\end{equation}
		and, similarly,
		\[
		\zeta_2=\phi_2\quad\text{in }Q_0,
		\qquad
		\zeta_2=\phi_2\quad\text{in }
		\widetilde\Omega_\varepsilon\setminus Q_3.
		\]
		By the definition of $T_1$, we know that
		\begin{equation}\label{eq6.6}
			- \Delta (v_1^* + \zeta_{1}) = \mu_1 |v_1^* + \zeta_{1}|^{2^*- 2}(v_1^* + \zeta_{1})+ \sum_{h=0}^N c_h Z^h, \quad~\text{in}~ M_3.
		\end{equation}
		Note that
		\begin{equation}\label{eq6.7}
			- \Delta (v_1^* + \phi_1) = \mu_1 |v_1^* +\phi_1|^{2^*- 2}(v_1^* + \phi_1)+ \beta |v_1^* + \phi_1|^{\frac{2^*}{2}-2}(v_1^* + \phi_1) |v_2^* + \phi_2|^{\frac{2^*}{2}},\quad~\text{in}~ \widetilde{\Omega}_{\varepsilon}\setminus(M_0\cup Q_0).
		\end{equation}
		By Lemma~\ref{lem4.6}, the coupling vanishes in $M_3$, and
		$\operatorname{supp}Z^h\subset M_0$ makes the multiplier term vanish on
		$M_3\setminus M_0$.  Hence \eqref{eq6.5}, \eqref{eq6.6}, and \eqref{eq6.7} give
		\begin{equation*}
			-\Delta (\zeta_{1} -\phi_1)=\mu_1 |v_1^* + \zeta_{1}|^{2^*- 2}(v_1^* + \zeta_{1})-\mu_1 |v_1^* +\phi_1|^{2^*- 2}(v_1^* + \phi_1),~\text{in}~M_3 \setminus M_0,
		\end{equation*}
		and $\zeta_{1} -\phi_1\in H_0^{1}(M_3 \setminus M_0)$.
		
		By the integral mean value formula,
		\[
		-\Delta(\zeta_1-\phi_1)=W_4(\zeta_1-\phi_1),\qquad
		W_4=\mu_1(2^*-1)\int_0^1
		|v_1^*+\phi_1+t(\zeta_1-\phi_1)|^{2^*-2}\,dt.
		\]
		Multiplying by $(\zeta_1 - \phi_1)$ and integrating by parts over $M_3 \setminus M_0$, we obtain
		\begin{equation*}
			\int_{M_3 \setminus M_0} |\nabla (\zeta_1 - \phi_1)|^2 \, \mathrm{d}y = \int_{M_3 \setminus M_0} W_4 (\zeta_1 - \phi_1)^2 \, \mathrm{d}y.
		\end{equation*}
		By H\"older's inequality and the Sobolev embedding $H_0^1 \hookrightarrow L^{2^*}$,
		\begin{equation*}
			\int_{M_3 \setminus M_0} W_4 (\zeta_1 - \phi_1)^2 \, \mathrm{d}y \leq \|W_4\|_{L^{N/2}} \|\zeta_1 - \phi_1\|_{L^{2^*}}^2 \leq C \|W_4\|_{L^{N/2}} \int_{M_3 \setminus M_0} |\nabla (\zeta_1 - \phi_1)|^2 \, \mathrm{d}y.
		\end{equation*}
		Using $|a+b|^{2^*-2}\le C(|a|^{2^*-2}+|b|^{2^*-2})$ and the Sobolev
		embedding, we obtain
		\[
		\|W_4\|_{L^{N/2}(M_3\setminus M_0)}
		\le C_N\|v_1^*\|_{L^{2^*}(M_3\setminus M_0)}^{2^*-2}
		+C\|\phi_1\|_{L^{2^*}(M_3\setminus M_0)}^{2^*-2}
		+C\|\zeta_1\|_{L^{2^*}(M_3\setminus M_0)}^{2^*-2}\leq C_NR_0^{-2}+C\varepsilon^{\frac{N-2\theta}{4}(2^*-2)}.
		\]
		Enlarging $R_0$ if necessary, we obtain $C\|W_4\|_{L^{N/2}(M_3\setminus M_0)}<1$ for
		$\varepsilon$ small.  
		The preceding energy identity then yields
		$\zeta_1=\phi_1$ in $M_3\setminus M_0$.  The same conclusion holds for
		the second component, so $T(\Phi)=\Phi$ globally.  Lemma~\ref{lem6.1}
		therefore gives $\Phi\in\Lambda_\varepsilon$, and the inner and outer
		equations together give \eqref{eq6.4}.  The two smallness inequalities are
		inherited from the outer solution.
		
		(ii) Conversely, keep the hypothesis $\Phi\in S(\Lambda_\varepsilon)$ of
		this lemma, and suppose also that $\Phi\in\Lambda_\varepsilon$ solves
		\eqref{eq6.4} with the stated smallness bounds.  As an element of the actual
		outer image, $\Phi$ already has the exterior boundary conditions and the
		dead-core property.  Since $S$ depends only on the local prescribed data
		and preserves them, outer uniqueness gives $S(\Phi)=\Phi$.  It therefore
		remains to prove $T(\Phi)=\Phi$, and then $S(T(\Phi))=S(\Phi)=\Phi$.
		
		On $M_3$ the function $\phi_1$ solves \eqref{eq5.1} with datum
		$\phi_1|_{\partial M_3}$: the coupling term vanishes by Lemma~\ref{lem4.6},
		and its moments vanish because $\operatorname{supp}Z^h\subset M_0$ and
		$\Phi\in\mathbf K_{\mathbf d,\boldsymbol\tau,\varepsilon}^\perp$.  Let
		$\varphi_{1,0}^*$ be the harmonic lifting of that datum; it has zero moments,
		so the zero-boundary part $\phi_1-\varphi_{1,0}^*$ lies in $K_{M_3}^\perp$.
		By Lemma~\ref{lem6.1} the datum satisfies \eqref{eq5.26}, and by the local
		$H^1$ bounds, \eqref{eq5.21} and the exterior estimates,
		\[
		\|\phi_1-\varphi_{1,0}^*\|_{L^q(M_3)}
		\le\|\phi_1\|_{L^q(M_3)}+\|\varphi_{1,0}^*\|_{L^q(M_3)}=o(1),
		\]
		which is at most the radius $a_R$ of Proposition~\ref{prop5.7} once
		$\varepsilon$ is small, since $a_R$ depends only on $R$ and $N$.  Both
		$\phi_1$ and $T_1(\phi_1)$ solve \eqref{eq5.18} with this datum, so the
		uniqueness part of Proposition~\ref{prop5.7} gives $\phi_1=T_1(\phi_1)$ on
		$M_3$.  The same argument on $Q_3$ gives $\phi_2=T_2(\phi_2)$.  Thus
		$T(\Phi)=\Phi$, hence $S(T(\Phi))=S(\Phi)=\Phi$.
	\end{proof}

	\begin{Lem}\label{lem6.4}
		For the matched solution constructed above,
		\begin{equation}\label{eq6.8}
			\|\nabla\phi_1\|_{L^2(\widetilde\Omega_\varepsilon)}
			+\|\nabla\phi_2\|_{L^2(\widetilde\Omega_\varepsilon)}
			\le C\varepsilon^{\frac{N-2\theta}{4}}
			\Bigl(\log\frac1\varepsilon\Bigr)^{1/2}.
		\end{equation}
	\end{Lem}
	\begin{proof}
		By Lemma~\ref{lem6.3} the matched pair satisfies $T(\Phi)=\Phi$, so
		$\phi_1|_{M_3}=T_1(\phi_1|_{\partial M_3})$ with datum in
		$\mathcal D_{1,\varepsilon}$; by Remark~\ref{rem4.7} the coupling vanishes on $M_3$.
		Lemma~\ref{lem5.9} and the $H^{1/2}$ bound of Corollary~\ref{cor4.4}
		therefore give
		\[
		\|\nabla\phi_1\|_{L^2(M_3)}
		\le C\|\phi_1\|_{H^{1/2}(\partial M_3)}
		+CR^{-\frac32}\varepsilon^{\frac{N-2\theta}{4}}
		\le CR^{-\frac12}\varepsilon^{\frac{N-2\theta}{4}}.
		\]
		On $\widetilde\Omega_\varepsilon\setminus M_3$ one has
		$|y-\widetilde a_1|\ge R_0+1$ and $|y-\widetilde\xi_1|\le c\varepsilon^{-1/2}$,
		so the gradient form of Lemma~\ref{lem4.3} given in
		Corollary~\ref{cor4.4} yields
		\[
		|\nabla\phi_1(y)|
		\le CR_0^{N/2}\varepsilon^{\frac{N-2\theta}{4}}
		(d_1+|y-\widetilde\xi_1|)^{-N/2},
		\]
		and integrating in polar coordinates,
		\[
		\int_{\widetilde\Omega_\varepsilon\setminus M_3}|\nabla\phi_1|^2
		\le C\varepsilon^{\frac{N-2\theta}{2}}
		\int_{3R}^{\varepsilon^{-1/2}}\frac{\mathrm dr}{r}
		\le C\varepsilon^{\frac{N-2\theta}{2}}\log\frac1\varepsilon .
		\]
		Adding the two contributions and repeating the argument on $Q_3$ for
		$\phi_2$ gives \eqref{eq6.8}.  
	\end{proof}
	
	\section{Finite-dimensional reduction}\label{sec7}
	Section~\ref{sec5} gives
	\begin{equation}\label{eq7.1}
		\begin{cases}
			-\Delta v_1-\mu_1|v_1|^{2^*-2}v_1
			-\beta |v_1|^{\frac{2^*}{2}-2}v_1|v_2|^{\frac{2^*}{2}}
			=\displaystyle\sum_{h=0}^Nc_hZ^h,\\[3mm]
			-\Delta v_2-\mu_2|v_2|^{2^*-2}v_2
			-\beta |v_2|^{\frac{2^*}{2}-2}v_2|v_1|^{\frac{2^*}{2}}
			=\displaystyle\sum_{h=0}^N\lambda_hY^h.
		\end{cases}
	\end{equation}
	It remains to choose $(\mathbf d,\boldsymbol\tau)$ so that
	$c_h=\lambda_h=0$ for $h=0,\ldots,N$.  Recall that
	$v_i=v_i^*+\phi_i$.
	
	Let $R_\varepsilon:=\varepsilon^{-\alpha}$, where $\alpha$ is a fixed number with
	\[
	\frac{N-4+2\theta}{2N}<\alpha<\frac{N-4}{2(N-2)}<\frac12 .
	\]
	This interval is nonempty precisely when $\theta<\frac{N-4}{N-2}$.  With such a choice,
	and for $\varepsilon$ small,
	\[
	3R<R_\varepsilon\ll\varepsilon^{-\frac{N-4}{2(N-2)}}<\varepsilon^{-1/2},
	\]
	because $R$ is fixed.
	Choose radial cutoff functions $\chi_i$ centered at $\widetilde\xi_i$ such
	that
	\[
	\chi_i=1\quad\hbox{on }B_{R_\varepsilon}(\widetilde\xi_i),
	\qquad
	\chi_i=0\quad\hbox{outside }B_{2R_\varepsilon}(\widetilde\xi_i),
	\]
	and
	\[
	0\le\chi_i\le1,
	\qquad |\nabla\chi_i|\le \frac{C}{R_\varepsilon},
	\qquad |\Delta\chi_i|\le \frac{C}{(R_\varepsilon)^2}.
	\]
	For small $\varepsilon$, $\chi_i=1$ on the support of every $Z_i^h$,
	$h=0,\ldots,N$.
	
	Since $\alpha<\frac{N-4}{2(N-2)}$ and
	$|\widetilde\xi_i-\widetilde a_i|\le\gamma_*^{-1}$, the supports of
	$\chi_1,\chi_2$ lie in the dead-core balls $B_\rho(\widetilde a_1)$,
	$B_\rho(\widetilde a_2)$ of Lemma~\ref{lem4.6} with
	$\rho=c_1\varepsilon^{-\frac{N-4}{2(N-2)}}$; hence $v_2=0$ on
	$\operatorname{supp}\chi_1$ and $v_1=0$ on $\operatorname{supp}\chi_2$, and
	every coupling term below vanishes on the support of the corresponding
	cutoff.  Hence each integrand below is the projection of the $i$-th
	left-hand side of \eqref{eq7.1} onto the corresponding kernel direction.
	
	Define the reduced residuals by
	\begin{align*}
		\widetilde c_0&=\int_{\widetilde\Omega_\varepsilon}
		\bigl(-\Delta v_1-\mu_1|v_1|^{2^*-2}v_1\bigr)
		\chi_1\partial_{d_1}(\widetilde P_\varepsilon U_{d_1,\widetilde\xi_1}),\\
		\widetilde c_k&=\int_{\widetilde\Omega_\varepsilon}
		\bigl(-\Delta v_1-\mu_1|v_1|^{2^*-2}v_1\bigr)
		\chi_1\partial_{\widetilde\xi_{1,k}}(\widetilde P_\varepsilon U_{d_1,\widetilde\xi_1}),
		\qquad k=1,\ldots,N,\\
		\widetilde\lambda_0&=\int_{\widetilde\Omega_\varepsilon}
		\bigl(-\Delta v_2-\mu_2|v_2|^{2^*-2}v_2\bigr)
		\chi_2\partial_{d_2}(\widetilde P_\varepsilon U_{d_2,\widetilde\xi_2}),\\
		\widetilde\lambda_k&=\int_{\widetilde\Omega_\varepsilon}
		\bigl(-\Delta v_2-\mu_2|v_2|^{2^*-2}v_2\bigr)
		\chi_2\partial_{\widetilde\xi_{2,k}}(\widetilde P_\varepsilon U_{d_2,\widetilde\xi_2}),
		\qquad k=1,\ldots,N .
	\end{align*}
	
	Put $\eta_{i,0}=d_i$ and
	$\eta_{i,k}=\widetilde\xi_{i,k}$ for $k=1,\ldots,N$, the directions of the
	kernel; the free parameters of the reduction are $(d_i,\tau_i)$, and
	$\eta_{i,k}$ is the geometric label for the same direction, related to
	$\tau_{i,k}$ by $\partial_{\tau_{i,h}}=d_i\,\partial_{\widetilde\xi_{i,h}}$
	at fixed $d_i$ (this constant factor is absorbed in the normalization of the
	reduced map below).  Then
	\[
	\widetilde c_k=\sum_{h=0}^NM^{(1)}_{kh}c_h,
	\qquad
	\widetilde\lambda_k=\sum_{h=0}^NM^{(2)}_{kh}\lambda_h,
	\]
	where
	\[
	M^{(i)}_{kh}
	=\int_{\widetilde\Omega_\varepsilon}
	Z_i^h\,\chi_i\,
	\partial_{\eta_{i,k}}
	(\widetilde P_\varepsilon U_{d_i,\widetilde\xi_i}).
	\]
	Because $\chi_i=1$ on $\operatorname{supp}Z_i^h$,
	\[
	M^{(i)}_{kh}=\int_{\widetilde\Omega_\varepsilon}
	\chi_0(y-\widetilde\xi_i)
	(\widetilde P_\varepsilon\psi_i^h)
	(\widetilde P_\varepsilon\psi_i^k).
	\]
	After translating the center, its limit is the same integral with
	$\psi_{d_i,0}^h,\psi_{d_i,0}^k$ in place of the projected functions, i.e.\ the
	Gram matrix $M^{(i)}_\infty=\bigl(\int_{\mathbb R^N}\chi_0\,
	\partial_{\eta_{i,k}}U\,\partial_{\eta_{i,h}}U\bigr)_{k,h}$ of the derivatives of the limit bubble.  By the nondegeneracy of the limit kernel, $M^{(i)}_\infty$ is positive definite; by
	continuity in the center and compactness of the scale interval, $M^{(i)}$ is
	uniformly invertible for small $\varepsilon$.  Hence
	\[
	(c_0,\ldots,c_N)=0\iff
	(\widetilde c_0,\ldots,\widetilde c_N)=0,
	\]
	with the analogous equivalence for $(\lambda_0,\ldots,\lambda_N)$.
	
	We now estimate these reduced residuals.
	The two auxiliary estimates used below --- the quadratic estimate for
	$\mathcal N_{v_i^*}$ and the projection-defect bound \eqref{eqB.proj} --- are
	collected, together with their proofs, in Appendix~\ref{secB}.
	
	The coefficient expansions in Lemmas~\ref{lem7.4} and~\ref{lem7.5}
	differ only in the parameter derivative that is paired with the reduced
	residual.  The estimates that are common to both derivatives are collected
	in Lemma~\ref{lem7.res} below, and the expansion of the main term in
	Lemma~\ref{lem7.mainterm}.
	
	\begin{Lem}\label{lem7.res}
		Let $\eta$ be any of the parameter derivatives $\partial_{d_i}$ or
		$\partial_{\widetilde\xi_{i,h}}$, $h=1,\ldots,N$.  Define
		\begin{equation}\label{eq7.2}
			\begin{aligned}
				F_1&:=-\mu_i\int_{\widetilde\Omega_\varepsilon}
				\left[|v_i|^{2^*-2}v_i-(v_i^*)^{2^*-1}
				-(2^*-1)(v_i^*)^{2^*-2}\phi_i\right]\chi_i\,\partial_\eta\bigl(\widetilde P_\varepsilon U_{d_i,\widetilde\xi_i}\bigr),\\
				F_2&:=\int_{\widetilde\Omega_\varepsilon}
				\phi_i\,\chi_i\left[\partial_\eta\bigl(U_{d_i,\widetilde\xi_i}^{2^*-1}\bigr)
				-\partial_\eta\bigl(\widetilde P_\varepsilon
				U_{d_i,\widetilde\xi_i}\bigr)^{2^*-1}\right],\\
				F_3&:=-\int_{\widetilde\Omega_\varepsilon}
				\phi_i\left(\partial_\eta\bigl(\widetilde P_\varepsilon
				U_{d_i,\widetilde\xi_i}\bigr)\Delta\chi_i
				+2\nabla\chi_i\cdot\nabla\partial_\eta\bigl(\widetilde P_\varepsilon
				U_{d_i,\widetilde\xi_i}\bigr)\right).
			\end{aligned}
		\end{equation}
		Then
		\begin{equation}\label{eq7.3}
			|F_j|=o\left(\varepsilon^{\frac{N-2}{2}}\right),
			\qquad j=1,2,3,
		\end{equation}
	\end{Lem}
	\begin{proof}
		The bracket in the definition of $F_1$ is precisely
		$\mathcal N_{v_i^*}(\phi_i)$, and the weight
		$\chi_i\partial_\eta(\widetilde P_\varepsilon U_{d_i,\widetilde\xi_i})$ is
		bounded by $Cv_i^*$ on its support, so \eqref{eqB.quad} gives
		$|F_1|=o(\varepsilon^{(N-2)/2})$.
		
		For $F_2$, let $\Pi_i$ be the projection defect of
		Corollary~\ref{corB.2}.  For any parameter derivative
		$\partial_\eta$,
		\[
		|\partial_\eta(U_{d_i,\widetilde\xi_i}^{2^*-1}
		-(\widetilde P_\varepsilon U_{d_i,\widetilde\xi_i})^{2^*-1})|
		\le C U_{d_i,\widetilde\xi_i}^{2^*-2}\Pi_i.
		\]
		Indeed, $0\le\Pi_i\le U_{d_i,\widetilde\xi_i}$,
		$|\partial_\eta U_{d_i,\widetilde\xi_i}|\le CU_{d_i,\widetilde\xi_i}$,
		and $|\partial_\eta\Pi_i|\le C\Pi_i$; also
		$|U_{d_i,\widetilde\xi_i}^b-(\widetilde P_\varepsilon
		U_{d_i,\widetilde\xi_i})^b|\le CU_{d_i,\widetilde\xi_i}^{b-1}\Pi_i$
		for every fixed $b>0$.  Weighted Cauchy--Schwarz, \eqref{eqB.proj} and
		\eqref{eq6.8} therefore give
		\begin{equation}\label{eq7.4}
			\begin{aligned}
				|F_2|
				&\le C\left(\int\chi_iU_{d_i,\widetilde\xi_i}^{2^*-2}\Pi_i^2\right)^{1/2}
				\left(\int\chi_iU_{d_i,\widetilde\xi_i}^{2^*-2}\phi_i^2\right)^{1/2}\\
				&\le C\varepsilon^{(N-2\theta)/4}\Bigl(\log\frac1\varepsilon\Bigr)^{1/2}
				\left[\varepsilon^{(N-2)/2}(R_\varepsilon)^{(N-4)/2}
				+\varepsilon^{N/4}\right]
				=o(\varepsilon^{(N-2)/2}).
			\end{aligned}
		\end{equation}
		For the first term the excess exponent is
		$(N-2\theta)/4-\alpha(N-4)/2>1-\theta/2>0$; for the second it is
		$1-\theta/2$.
		
		For $F_3$, on $R_\varepsilon\le r=|y-\widetilde\xi_i|\le2R_\varepsilon$,
		Appendix~\ref{secB} and $\alpha<1/2$ imply
		\[
		|\partial_{\widetilde\xi_{i,h}}(\widetilde P_\varepsilon
		U_{d_i,\widetilde\xi_i})|\le Cr^{1-N},
		\qquad
		|\partial_{d_i}(\widetilde P_\varepsilon
		U_{d_i,\widetilde\xi_i})|\le Cr^{2-N}.
		\]
		Their gradients gain one more factor $r^{-1}$: apply the interior
		gradient estimate on balls of radius comparable to $r$ to the
		differentiated Poisson equations.  These balls avoid both holes and the
		outer boundary.  On this collar $R_\varepsilon\le r=|y-\widetilde\xi_i|\le2R_\varepsilon\le c\,\varepsilon^{-1/2}$
		for small $\varepsilon$, so Corollary~\ref{cor4.4} gives
		$|\phi_i|\le C\varepsilon^{(N-2\theta)/4}r^{-N/2}$.
		Hence
		\begin{equation}\label{eq7.5}
			|F_3|\le
			C\varepsilon^{(N-2\theta)/4}
			\begin{cases}
				(R_\varepsilon)^{-N/2-1},&\eta=\widetilde\xi_{i,h},\\[1mm]
				(R_\varepsilon)^{-N/2},&\eta=d_i,
			\end{cases}
			=o(\varepsilon^{(N-2)/2}),
		\end{equation}
		the last step following from the choice of $\alpha$: the interval for
		$\alpha$ is nonempty by the choice of $\theta$ in Section~\ref{sec1}, and
		$\alpha>(N-4+2\theta)/(2N)$ covers both the scale derivative and, a fortiori,
		the translations.
	\end{proof}
	
	\begin{Lem}\label{lem7.mainterm}
		Let $\eta$ be any of the parameter derivatives $\partial_{d_i}$ or
		$\partial_{\widetilde\xi_{i,h}}$, $h=1,\ldots,N$.  Then the main term
		admits the expansion
		\begin{equation}\label{eq7.mainterm}
			\begin{aligned}
				&\mu_i^{-\frac{1}{2^*-2}}
				\int_{\widetilde\Omega_{\varepsilon}}
				\left[U_{d_i,\widetilde\xi_i}^{2^*-1}
				-\bigl(\widetilde{P}_{\varepsilon}U_{d_i,\widetilde\xi_i}\bigr)^{2^*-1}\right]
				\chi_i\,\partial_\eta\bigl(\widetilde P_\varepsilon U_{d_i,\widetilde\xi_i}\bigr)\\
				={}&\mu_i^{-\frac{1}{2^*-2}}(2^*-1)\alpha_N
				\Bigl[d_i^{\frac{N-2}{2}}
				\widetilde H_\varepsilon(\widetilde\xi_i,\widetilde\xi_i)
				\int_{\mathbb R^N}
				U_{d_i,\tilde{\xi}_i}^{2^*-2}\,
				\partial_\eta U_{d_i,\tilde{\xi}_i}\,dy\\
				&\quad+\frac{\varepsilon^{\frac{N-2}{2}}r_i^{N-2}}
				{d_i^{\frac{N-2}{2}}(1+|\tau_i|^2)^{\frac{N-2}{2}}}
				\int_{\widetilde\Omega_{\varepsilon}}
				\frac{\chi_iU_{d_i,\tilde{\xi}_i}^{2^*-2}
					\partial_\eta U_{d_i,\tilde{\xi}_i}}
				{|y-\tilde a_i|^{N-2}}\,dy\Bigr]+o\left(\varepsilon^{\frac{N-2}{2}}\right).
			\end{aligned}
		\end{equation}
		where $\widetilde H_\varepsilon$ is the regular part of the projected bubble
		(Lemma~\ref{lemB.1}).  For every translation derivative
		$\eta=\widetilde\xi_{i,h}$ the $\widetilde H_\varepsilon$-term, i.e.\ the
		first term in the bracket of \eqref{eq7.mainterm}, is $o(\varepsilon^{(N-2)/2})$,
		by oddness of the parameter kernel.   
	\end{Lem}
	\begin{proof}
		Lemma~\ref{lemB.1} gives, for every parameter derivative $\partial_\eta$,
		\begin{equation}\label{eq7.7}
			\begin{aligned}
				&\mu_i^{-\frac{1}{2^*-2}}
				\int_{\widetilde\Omega_{\varepsilon}}
				\left[U_{d_i,\tilde{\xi}_i}^{2^*-1}
				-\bigl(\widetilde{P}_{\varepsilon}U_{d_i,\tilde{\xi}_i}\bigr)^{2^*-1}\right]
				\chi_i\,\partial_\eta\bigl(\widetilde P_\varepsilon U_{d_i,\tilde{\xi}_i}\bigr)\\
				={}&\mu_i^{-\frac{1}{2^*-2}}(2^*-1)
				\int_{\widetilde\Omega_{\varepsilon}}
				U_{d_i,\tilde{\xi}_i}^{2^*-2}\Pi_i\,\chi_i\,\partial_\eta\bigl(\widetilde P_\varepsilon U_{d_i,\tilde{\xi}_i}\bigr)+O\left(\int_{\widetilde\Omega_{\varepsilon}},
				U_{d_i,\tilde{\xi}_i}^{2^*-2}\Pi_i^2\right)\\
				={}&\mu_i^{-\frac{1}{2^*-2}}(2^*-1)
				\int_{\widetilde\Omega_{\varepsilon}}
				U_{d_i,\tilde{\xi}_i}^{2^*-2}\Pi_i\,\chi_i
				\frac{\partial U_{d_i,\tilde{\xi}_i}}{\partial\eta}
				-\mu_i^{-\frac{1}{2^*-2}}(2^*-1)
				\int_{\widetilde\Omega_{\varepsilon}}
				U_{d_i,\tilde{\xi}_i}^{2^*-2}\Pi_i\,\chi_i
				\frac{\partial\Pi_i}{\partial\eta}\\
				&+O\left(\int_{\widetilde\Omega_{\varepsilon}}
				\chi_iU_{d_i,\tilde{\xi}_i}^{2^*-2}\Pi_i^2\right)\\
				={}&\mu_i^{-\frac{1}{2^*-2}}(2^*-1)
				\int_{\widetilde\Omega_{\varepsilon}}
				U_{d_i,\tilde{\xi}_i}^{2^*-2}\Pi_i\,\chi_i
				\frac{\partial U_{d_i,\tilde{\xi}_i}}{\partial\eta}
				+o\left(\varepsilon^{\frac{N-2}{2}}\right),
			\end{aligned}
		\end{equation}
		where the two remainder terms are absorbed as follows: the first
		because $|\partial_\eta(\widetilde P_\varepsilon
		U_{d_i,\widetilde\xi_i})|\le C\widetilde P_\varepsilon
		U_{d_i,\widetilde\xi_i}$, the second because
		$|\partial_\eta\Pi_i|\le C\Pi_i$ and \eqref{eqB.proj} makes the
		weighted square integrable to $o(\varepsilon^{(N-2)/2})$.  Inserting next the
		explicit form of $\Pi_i$ from Lemma~\ref{lemB.1}, the last integral splits
		into three pieces: the piece containing
		$\widetilde H_\varepsilon(y,\widetilde\xi_i)$, the piece containing
		$|y-\widetilde a_i|^{-(N-2)}$, and the piece containing the remainder
		$\mathcal R_i$.  The second piece is precisely the second term of the bracket
		in \eqref{eq7.mainterm}, while the first piece is the first term of that
		bracket, as we now verify.
		For the piece containing $\widetilde H_\varepsilon$, replace
		$\widetilde H_\varepsilon(y,\widetilde\xi_i)$ first by
		$\widetilde H_\varepsilon(\widetilde\xi_i,\widetilde\xi_i)$ and then by
		their difference.  The constant part contributes
		$\mu_i^{-\frac{1}{2^*-2}}(2^*-1)\alpha_Nd_i^{(N-2)/2}
		\widetilde H_\varepsilon(\widetilde\xi_i,\widetilde\xi_i)
		\int_{\mathbb R^N}U_{d_i,\widetilde\xi_i}^{2^*-2}
		\partial_\eta U_{d_i,\widetilde\xi_i}\,dy$, up to the
		error made by removing the own-hole ball, which is
		$O(\varepsilon^{N/2})$, and by the radial cutoff, which is smaller
		still; both are $o(\varepsilon^{(N-2)/2})$ after the
		prefactor $\varepsilon^{(N-2)/2}$ of
		$\widetilde H_\varepsilon(\widetilde\xi_i,\widetilde\xi_i)$.  For the
		fluctuation, with $\xi_i=\sqrt\varepsilon\widetilde\xi_i$,
		\[
		|\widetilde H_\varepsilon(y,\widetilde\xi_i)
		-\widetilde H_\varepsilon(\widetilde\xi_i,\widetilde\xi_i)|
		=\varepsilon^{(N-2)/2}|H(\sqrt\varepsilon y,\xi_i)-H(\xi_i,\xi_i)|
		\le C\varepsilon^{(N-1)/2}|y-\widetilde\xi_i|,
		\]
		and the resulting weighted integral is finite uniformly in the
		parameters, so the fluctuation is $O(\varepsilon^{(N-1)/2})
		=o(\varepsilon^{(N-2)/2})$.  This identifies the first term of the bracket in \eqref{eq7.mainterm}; for a
		translation derivative the full-space integral of
		$U_{d_i,\widetilde\xi_i}^{2^*-2}\partial_{\widetilde\xi_{i,h}}
		U_{d_i,\widetilde\xi_i}$ vanishes by oddness, so that term is $o(\varepsilon^{(N-2)/2})$.
		
		Finally, on the support of $\chi_i$, $|y-\widetilde a_l|\ge
		c\varepsilon^{-1/2}$ for $l\ne i$.  The weight
		$U_{d_i,\widetilde\xi_i}^{2^*-2}|\partial_\eta
		U_{d_i,\widetilde\xi_i}|$ is bounded by
		$C(1+|y-\widetilde\xi_i|)^{-N-2}$ for the scale derivative and has one
		additional decay power for translations.  Both weights have bounded
		integrals after multiplication by $1+|y-\widetilde a_i|^{2-N}$,
		uniformly for the allowed center shifts.  Therefore \eqref{eqB.2} bounds the
		remaining piece by
		$C\bigl(\varepsilon^{N/2}+\varepsilon^{(N-1)/2}+\varepsilon^{3(N-2)/2}\bigr)
		=o(\varepsilon^{(N-2)/2})$.
		This completes the proof of Lemma~\ref{lem7.mainterm}.
	\end{proof}

	\begin{Lem}\label{lem7.hole}
		Let $\eta$ be one of the derivatives $\widetilde\xi_{1,i}$ $(i=1,\ldots,N)$
		and $d_1$, and let $\mathcal K_\eta$ denote $\partial_{z_i}U$ if
		$\eta=\widetilde\xi_{1,i}$ and $\Lambda U:=\frac{N-2}{2}U+z\cdot\nabla U$ if
		$\eta=d_1$.  Then
		\begin{equation}\label{eq7.hole}
			\int_{\widetilde\Omega_\varepsilon}
			\frac{\chi_1U_{d_1,\widetilde\xi_1}^{2^*-2}\,\partial_\eta U_{d_1,\widetilde\xi_1}}
			{|y-\widetilde a_1|^{N-2}}\,dy
			=-d_1^{-\frac N2}\left[\frac{N(N-2)\omega_N}{2^*-1}\,
			\mathcal K_\eta(-\tau_1)+O\bigl(\varepsilon+(R_\varepsilon)^{-N-1}\bigr)\right].
		\end{equation}
	\end{Lem}
	\begin{proof}
		With $y=\widetilde\xi_1+d_1z$ the own hole becomes
		$B_{r_1\sqrt\varepsilon/d_1}(-\tau_1)$ and the other hole lies outside the
		cutoff support; moreover $\mathrm dy=d_1^N\mathrm dz$,
		$|y-\widetilde a_1|^{-(N-2)}=d_1^{-(N-2)}|z+\tau_1|^{-(N-2)}$ and
		$U^{2^*-2}_{d_1,\widetilde\xi_1}\partial_\eta U_{d_1,\widetilde\xi_1}
		=-d_1^{-\frac N2-2}U^{2^*-2}\mathcal K_\eta$, since
		$\partial_{\widetilde\xi_{1,i}}U_{d_1,\widetilde\xi_1}
		=-d_1^{-\frac N2}\partial_{z_i}U$ and
		$\partial_{d_1}U_{d_1,\widetilde\xi_1}=-d_1^{-\frac N2}\Lambda U$; the powers
		of $d_1$ therefore combine to $-d_1^{-\frac N2}$.  Near $-\tau_1$ both bubble
		factors are bounded, so omitting the hole costs at most
		$C\int_0^{r_1\sqrt\varepsilon/d_1}r\,dr=O(\varepsilon)$; on the support of
		$\chi_1-1$, which begins at radius $R_\varepsilon/d_1$, the integrand is
		bounded by $C|z|^{-2N-1}$ for a translation derivative and by
		$C|z|^{-2N}$ for the scale derivative, so replacing $\chi_1$ by $1$ costs
		$O\bigl((R_\varepsilon)^{-N}\bigr)$ in either case.  Finally
		$-\Delta_z|z+\tau_1|^{-(N-2)}=N(N-2)\omega_N\delta_{-\tau_1}$ and
		$-\Delta\mathcal K_\eta=(2^*-1)U^{2^*-2}\mathcal K_\eta$ (for the derivatives
		of $U$ and for $\Lambda U$), so integrating by parts
		\[
		\int_{\mathbb R^N}\frac{U^{2^*-2}\mathcal K_\eta}{|z+\tau_1|^{N-2}}\,dz
		=\frac{1}{2^*-1}\int_{\mathbb R^N}
		\Bigl(-\Delta_z\frac{1}{|z+\tau_1|^{N-2}}\Bigr)\mathcal K_\eta\,dz
		=\frac{N(N-2)\omega_N}{2^*-1}\,\mathcal K_\eta(-\tau_1),
		\]
		which is \eqref{eq7.hole}.
	\end{proof}
	
	\begin{Lem}\label{lem7.4} The coefficients admit the expansions
		\begin{align*}
			\widetilde{c}_i&
			=-B_1 \varepsilon^{\frac{N-2}{2}} d_1^{1-N} r_1^{N-2} \frac{\tau_{1,i}}{(1+|\tau_1|^2)^{N-1}} + o\left(\varepsilon^{\frac{N-2}{2}}\right) 
			,
			\\
			\widetilde{\lambda}_i&
			=-B_2 \varepsilon^{\frac{N-2}{2}} d_2^{1-N} r_2^{N-2} \frac{\tau_{2,i}}{(1+|\tau_2|^2)^{N-1}} + o\left(\varepsilon^{\frac{N-2}{2}}\right),
		\end{align*}
		where $B_l = \mu_l^{-\frac{1}{2^*-2}} \alpha_N^{2}N (N-2)^2 \omega_N>0 $, $l=1,2$ and $i=1,\dots,N$.
	\end{Lem}
	
	\begin{proof} 
		By the definition of $\widetilde c_i$ in terms of
		$-\Delta v_1-\mu_1|v_1|^{2^*-2}v_1$, we have
		\begin{equation}\label{eq7.9}
			\begin{aligned}
				&\int_{\widetilde\Omega_{\varepsilon}}\Big(-\Delta v_1 - \mu_1 |v_1|^{2^*-2} v_1 \Big)\chi_1\frac{\partial (\widetilde{P}_{\varepsilon}U_{d_1, \tilde{\xi}_1})}{\partial \tilde{\xi}_{1,i}}\\
				=
				&\mu_1^{-\frac{1}{2^*-2}}\int_{\widetilde\Omega_{\varepsilon}}\left[U_{d_1,\tilde{\xi}_1}^{2^*-1}-\big(\widetilde{P}_{\varepsilon} U_{d_1,\widetilde{\xi}_1}\big)^{2^*-1}\right]\chi_1\frac{\partial \bigl(\widetilde{P}_{\varepsilon}U_{d_1,\tilde{\xi}_1}\bigr)}{\partial \tilde{\xi}_{1,i}} 
				+F_1+F_2+F_3
			\end{aligned}
		\end{equation}
		where $F_1,F_2,F_3$ are the three terms $F_j$ of
		\eqref{eq7.2} with $i=1$ and $\eta=\widetilde\xi_{1,i}$; by
		Lemma~\ref{lem7.res} each of them is $o(\varepsilon^{(N-2)/2})$.
		For the main term, \eqref{eq7.mainterm} with $i=1$ and
		$\eta=\widetilde\xi_{1,i}$ reduces, since the $\widetilde H_\varepsilon$-term
		of the bracket in \eqref{eq7.mainterm} vanishes for translation derivatives
		by oddness, to the hole term, whose integral is
		\eqref{eq7.hole}.  Using
		$\partial_{z_i}U(-\tau_1)=(N-2)\alpha_N\tau_{1,i}(1+|\tau_1|^2)^{-\frac N2}$
		and $2^*-1=\frac{N+2}{N-2}$, that second term equals
		\begin{equation}\label{eq7.10}
			-B_1\varepsilon^{\frac{N-2}{2}}d_1^{1-N}r_1^{N-2}
			\frac{\tau_{1,i}}{(1+|\tau_1|^2)^{N-1}}+o\bigl(\varepsilon^{\frac{N-2}{2}}\bigr),
			\qquad B_1:=\mu_1^{-\frac{1}{2^*-2}}\alpha_N^2N(N-2)^2\omega_N>0 .
		\end{equation}
		Combining \eqref{eq7.9} and \eqref{eq7.10} with the definition of $\widetilde c_i$ completes the proof.
	\end{proof}
	
	\begin{Lem}\label{lem7.5} The coefficients admit the expansions
		\begin{align*}
			\widetilde{c}_0&
			=\tfrac12B_1\,\varepsilon^{\frac{N-2}{2}}d_1^{N-3}H(a_1,a_1)
			-\tfrac12B_1\,\varepsilon^{\frac{N-2}{2}}
			d_1^{1-N}r_1^{N-2}
			\frac{1-|\tau_1|^2}{(1+|\tau_1|^2)^{N-1}}
			+o\left(\varepsilon^{\frac{N-2}{2}}\right),
			\\
			\widetilde{\lambda}_0&
			=\tfrac12B_2\,\varepsilon^{\frac{N-2}{2}}d_2^{N-3}H(a_2,a_2)
			-\tfrac12B_2\,\varepsilon^{\frac{N-2}{2}}
			d_2^{1-N}r_2^{N-2}
			\frac{1-|\tau_2|^2}{(1+|\tau_2|^2)^{N-1}}
			+o\left(\varepsilon^{\frac{N-2}{2}}\right),
		\end{align*}
		where $B_l$ is the constant of Lemma~\ref{lem7.4}.
	\end{Lem}
	
	\begin{proof} 
		By the definition of $\widetilde c_0$ in terms of
		$-\Delta v_1-\mu_1|v_1|^{2^*-2}v_1$, we have
		\begin{equation}\label{eq7.11}
			\begin{aligned}
				&\int_{\widetilde\Omega_{\varepsilon}}\Big(-\Delta v_1 - \mu_1 |v_1|^{2^*-2} v_1 \Big)\chi_1\frac{\partial (\widetilde{P}_{\varepsilon}U_{d_1, \tilde{\xi}_1})}{\partial d_1}\\
				=
				&\mu_1^{-\frac{1}{2^*-2}}\int_{\widetilde\Omega_{\varepsilon}}\left[U_{d_1,\tilde{\xi}_1}^{2^*-1}-\big(\widetilde{P}_{\varepsilon} U_{d_1,\widetilde{\xi}_1}\big)^{2^*-1}\right]\chi_1\frac{\partial \bigl(\widetilde{P}_{\varepsilon}U_{d_1,\tilde{\xi}_1}\bigr)}{\partial d_1} 
				+F_1+F_2+F_3
			\end{aligned}
		\end{equation}
		where $F_1,F_2,F_3$ are the three terms $F_j$ of
		\eqref{eq7.2} with $i=1$ and $\eta=d_1$; by
		Lemma~\ref{lem7.res} each of them is $o(\varepsilon^{(N-2)/2})$.
		For the main term, \eqref{eq7.mainterm} with $i=1$ and $\eta=d_1$ has both terms
		of the bracket in \eqref{eq7.mainterm} contributing.  Its $\widetilde H_\varepsilon$-term equals
		\begin{equation}
			\mu_1^{-\frac{1}{2^*-2}}(2^*-1)\alpha_Nd_1^{\frac{N-2}{2}}
			\widetilde H_\varepsilon(\tilde{\xi}_1,\tilde{\xi}_1)
			\int_{\mathbb R^N}U_{d_1,\widetilde\xi_1}^{2^*-2}\partial_{d_1}U_{d_1,\widetilde\xi_1}\,dy
			+o\left(\varepsilon^{\frac{N-2}{2}}\right),
		\end{equation}
		and, since \(\widetilde H_\varepsilon(\tilde{\xi}_1,\tilde{\xi}_1)
		=\varepsilon^{\frac{N-2}{2}}H(a_1,a_1)+o(\varepsilon^{\frac{N-2}{2}})\), and  
		\[
		\begin{aligned}
			\int_{\mathbb{R}^N} U_{d_1,\widetilde\xi_1}^{2^*-2} \partial_{d_1}U_{d_1,\widetilde\xi_1} \, \mathrm{d}y
			&= d_1^{\frac{N}{2}-2} \frac{N-2}{2(2^*-1)} \int_{\mathbb{R}^N} U^{2^*-1} \, \mathrm{d}z\\
			&= d_1^{\frac{N}{2}-2} \frac{N-2}{2(2^*-1)}\,N(N-2)\omega_N\alpha_N,
		\end{aligned}
		\]
		this first term equals
		\begin{equation}\label{eq7.13}
			\tfrac12B_1\varepsilon^{\frac{N-2}{2}}d_1^{N-3}H(a_1,a_1)
			+o\left(\varepsilon^{\frac{N-2}{2}}\right).
		\end{equation}

		For the hole term (the second term of the bracket in \eqref{eq7.mainterm}), since
		$U_{d_1,\widetilde\xi_1}(y)=d_1^{-\frac{N-2}{2}}U((y-\widetilde\xi_1)/d_1)$
		one has $\partial_{d_1}U_{d_1,\widetilde\xi_1}=-d_1^{-\frac N2}\Lambda U$ with
		$\Lambda U$ defined in Lemma~\ref{lem7.hole}, and
		$\Lambda U(-\tau_1)=\frac{N-2}{2}\alpha_N(1-|\tau_1|^2)(1+|\tau_1|^2)^{-\frac N2}$;
		by \eqref{eq7.hole} with $\eta=d_1$ and $2^*-1=\frac{N+2}{N-2}$ it equals
		\begin{equation}\label{eq7.14}
			-\tfrac12B_1\varepsilon^{\frac{N-2}{2}}d_1^{-(N-1)}r_1^{N-2}
			\frac{1-|\tau_1|^2}{(1+|\tau_1|^2)^{N-1}}+o\bigl(\varepsilon^{\frac{N-2}{2}}\bigr).
		\end{equation}
		
		Combining \eqref{eq7.13} and \eqref{eq7.14} with \eqref{eq7.mainterm} (applied with $\eta=d_1$), we obtain
		\begin{equation}\label{eq7.15}
			\begin{aligned}
				&\mu_1^{-\frac{1}{2^*-2}}\int_{\widetilde\Omega_{\varepsilon}}\left[U_{d_1,\tilde{\xi}_1}^{2^*-1}-\big(\widetilde{P}_{\varepsilon} U_{d_1,\widetilde{\xi}_1}\big)^{2^*-1}\right]\chi_1\frac{\partial \bigl(\widetilde{P}_{\varepsilon}U_{d_1,\tilde{\xi}_1}\bigr)}{\partial d_1}\\
				=&\tfrac12B_1 \varepsilon^{\frac{N-2}{2}}d_1^{N-3}H(a_1,a_1) -\tfrac12B_1 \; \varepsilon^{\frac{N-2}{2}} \; d_1^{-(N-1)} \; r_1^{N-2} \; \frac{1-|\tau_1|^2}{(1+|\tau_1|^2)^{N-1}} + o\bigl(\varepsilon^{\frac{N-2}{2}}\bigr). 
			\end{aligned}
		\end{equation}

		Combining \eqref{eq7.11} and \eqref{eq7.15} with the definition of $\widetilde c_0$ completes the proof.
	\end{proof}
	
	\begin{proof}[\textbf{Proof of Theorem~\ref{thm2.1}}]
		For $i=1,2$ set
		\[
		\mathcal C_{1,k}=\widetilde c_k,
		\qquad
		\mathcal C_{2,k}=\widetilde\lambda_k,
		\qquad k=0,\ldots,N.
		\]
		Define the exact normalized reduced map
		$\Psi_\varepsilon:\mathcal P_{\gamma_*}\to\mathbb R^{2N+2}$ by
		\begin{align*}
			(\Psi_\varepsilon)_{d_i}
			&=\varepsilon^{-\frac{N-2}{2}}d_i^{N-1}\mathcal C_{i,0},\\
			(\Psi_\varepsilon)_{\tau_i}
			&=-\frac{d_i^{N-1}(1+|\tau_i|^2)^{N-1}}
			{B_i r_i^{N-2}\varepsilon^{\frac{N-2}{2}}}
			(\mathcal C_{i,1},\ldots,\mathcal C_{i,N}).
		\end{align*}
		Lemma~\ref{lemC.2}, with the constants chosen on a slightly larger parameter box, gives continuity on $\mathcal P_{\gamma_*}$.  Lemmas~\ref{lem7.4} and~\ref{lem7.5}
		give, uniformly on $\mathcal P_{\gamma_*}$,
		\begin{align}
			(\Psi_\varepsilon)_{d_i}
			& =\tfrac12B_i\Bigl[d_i^{2N-4}H(a_i,a_i)
			-r_i^{N-2}\frac{1-|\tau_i|^2}{(1+|\tau_i|^2)^{N-1}}\Bigr]
			+o(1),\label{eq7.16}\\
			(\Psi_\varepsilon)_{\tau_i}&=\tau_i+o(1).
			\label{eq7.17}
		\end{align}
		
		As required in the initial choice of the parameter range, $\gamma_*>0$ is small enough that
		\[
		g_i(d):=d^{2N-4}H(a_i,a_i)-r_i^{N-2}
		\]
		is negative at $d=\gamma_*$ and positive at $d=\gamma_*^{-1}$.
		On $|\tau_i|=\gamma_*$,
		\[
		(\Psi_\varepsilon)_{\tau_i}\cdot\tau_i
		=|\tau_i|^2+o(1)>0.
		\]
		If a zero were situated on one of the two $d_i$-faces, then
		\eqref{eq7.17} would imply $\tau_1=\tau_2=o(1)$;
		\eqref{eq7.16} would then have the sign of $g_i$ on that
		face, a contradiction.  Hence the exact map has no zero on
		$\partial \mathcal P_{\gamma_*}$ for small $\varepsilon$.  The same argument applies to the straight-line
		homotopy between $\Psi_\varepsilon$ and
		\[
		\Psi_0(\mathbf d,\boldsymbol\tau)
		=\bigl(g_1(d_1),g_2(d_2),\tau_1,\tau_2\bigr).
		\]
		Therefore
		\[
		\deg(\Psi_\varepsilon,\mathcal P_{\gamma_*},0)
		=\deg(\Psi_0,\mathcal P_{\gamma_*},0)=1.
		\]
		There is consequently a zero
		$(\mathbf d_\varepsilon,\boldsymbol\tau_\varepsilon)$ of the exact reduced
		map.
		
		At this zero, \eqref{eq7.17} gives
		$\tau_{i,\varepsilon}\to0$, and the scale equations give the asymptotic
		balance
		\[
		d_{i,\varepsilon}^{2N-4}H(a_i,a_i)=r_i^{N-2}(1+o(1)).
		\]
		Thus
		\[
		d_{i,\varepsilon}
		=\left(\frac{r_i^{N-2}}{H(a_i,a_i)}\right)^{\frac1{2(N-2)}}
		(1+o(1)),
		\]
		which is the parameter limit stated in the theorem.
		
		By continuity in the center and compactness of the scale interval, the
		matrices $M^{(i)}$ are uniformly invertible for small $\varepsilon$.  Hence the vanishing of all components of
		$\Psi_\varepsilon$ implies
		$c_h=\lambda_h=0$ for $h=0,\ldots,N$, so the constructed pair solves the
		rescaled system without Lagrange multipliers.
		
		It remains to verify its sign.  Testing the $i$-th signed equation with
		$v_i^-:=\max\{-v_i,0\}$ gives
		\[
		\int_{\widetilde\Omega_\varepsilon}|\nabla v_i^-|^2
		\le \mu_i\int_{\widetilde\Omega_\varepsilon}(v_i^-)^{2^*}.
		\]
		If $v_i^-\not\equiv0$, the Sobolev inequality yields a fixed positive
		lower bound for $\|v_i^-\|_{2^*}$.  On the other hand,
		$v_i^*\ge0$ and hence $v_i^-\le|\phi_i|$; the correction estimate makes
		$\|\phi_i\|_{2^*}=o(1)$.  This contradiction proves $v_i\ge0$.  Neither
		component is identically zero, because the $L^{2^*}$-norm of its projected
		bubble is bounded away from zero whereas $\|\phi_i\|_{2^*}=o(1)$.
		The dead-core statement follows from Lemma~\ref{lem4.6} and the matching
		construction.  Scaling back by \eqref{eq2.1} completes the proof of
		Theorem~\ref{thm1.1} as well. At last, since $\theta$ is arbitrary chosen in $(0,\frac{N-4}{N-2})$,  enlarging $\theta$ in \eqref{eq6.8}, we can obtain  
		\eqref{estimate}.
	\end{proof}
	
	\appendix
	\section{Basic Estimates}
	
	We begin with an elementary inequality for products of weights centered at
	two different points, which is used in the matching argument of
	Section~\ref{sec6}.
	
	\begin{Lem}\label{lemA.1}
		Let $z_1,z_2\in\mathbb R^N$, $d_1,d_2\in[\gamma_*,\gamma_*^{-1}]$, and let
		$\alpha_1,\alpha_2>0$, $0<\tau\le\min\{\alpha_1,\alpha_2\}$.  Then there is
		$C=C(N,\gamma_*)>0$ such that
		\begin{equation}\label{eqA.1}
			\begin{aligned}
				&\frac{1}{(d_1+|y-z_1|)^{\alpha_1}}
				\frac{1}{(d_2+|y-z_2|)^{\alpha_2}}\\
				&\qquad\le\frac{C}{|z_1-z_2|^{\tau}}
				\left(
				\frac{1}{(d_1+|y-z_1|)^{\alpha_1+\alpha_2-\tau}}
				+\frac{1}{(d_2+|y-z_2|)^{\alpha_1+\alpha_2-\tau}}
				\right)
			\end{aligned}
		\end{equation}
		for every $y\in\mathbb R^N$.
	\end{Lem}
	\begin{proof}
		If $|y-z_1|\le|z_1-z_2|/2$, then $|y-z_2|\ge|z_1-z_2|/2$ and, since
		$d_2\le\gamma_*^{-1}$, one has $d_2+|y-z_2|\ge c|z_1-z_2|$.  Hence
		\[
		(d_2+|y-z_2|)^{-\alpha_2}\le C|z_1-z_2|^{-\alpha_2}
		\le C|z_1-z_2|^{-\tau}(d_1+|y-z_1|)^{-(\alpha_2-\tau)},
		\]
		because $\tau\le\alpha_2$ and $d_1+|y-z_1|\le C|z_1-z_2|$; multiplying by
		$(d_1+|y-z_1|)^{-\alpha_1}$ gives the first term on the right of
		\eqref{eqA.1}.  The case $|y-z_2|\le|z_1-z_2|/2$ is symmetric and gives
		the second term.  If $|y-z_i|>|z_1-z_2|/2$ for $i=1,2$, then both weights are
		bounded by $C|z_1-z_2|^{-\alpha_i}$, and
		\eqref{eqA.1} follows from
		$|z_1-z_2|^{-(\alpha_1+\alpha_2)}
		\le C|z_1-z_2|^{-\tau}(d_1+|y-z_1|)^{-(\alpha_1+\alpha_2-\tau)}$.
	\end{proof}
	
	In this appendix only, $P_\varepsilon$ denotes projection on the
	single-hole domain $\Omega\setminus\overline{B_{r\varepsilon}(a)}$.
	The estimates for the full two-hole domain are proved in
	Appendix~\ref{secB}.  For comparison with the usual one-hole expansion, we retain the
	following bounds (see also \cite[Lemma 3.1]{GeMussoPistoia2010}) and give
	the boundary comparison that proves them.

	\begin{Lem}
		\label{lemA.2}
		Let $a\in\Omega$, $r>0$, and let $\gamma_*,\varepsilon>0$ be small. Let $d>0$ and $\tau\in\mathbb{R}^N$ be such that $\gamma_*<d<\gamma_*^{-1}$, $|\tau|<\gamma_*^{-1}$. Finally, let $\delta=d\sqrt{\varepsilon}$ and $\xi=a+\delta\tau$. Let us define
		\[
		\mathcal R(x):= P_{\varepsilon}U_{\delta,\xi}(x)-U_{\delta,\xi}(x)+\alpha_N\delta^{\frac{N-2}{2}}H(x,\xi)+\alpha_N\frac{1}{\delta^{\frac{N-2}{2}}(1+|\tau|^2)^{\frac{N-2}{2}}}\left(\frac{r\varepsilon}{|x-a|}\right)^{N-2}.
		\]
		In the estimates below the parameters are the coordinates $(\delta,\tau)$; thus
		$\partial_{\tau_h}\mathcal R$ is taken at fixed $\delta$, and the scale derivative is
		taken at fixed $\tau$, written $\partial_\delta^{\,\tau}\mathcal R$ (the superscript
		always specifies the variable that is held fixed).  As throughout the paper,
		the scale derivative without superscript, $\partial_\delta \mathcal R$, keeps the center
		$\xi$ fixed; the two are related by
		$\partial_\delta \mathcal R=\partial_\delta^{\,\tau}\mathcal R-\delta^{-1}\tau\cdot\nabla_\tau \mathcal R$.
		Then there exists a positive constant $C>0$ depending on the fixed domain, $r$, $\gamma_*$, and the distance of $a$ from the outer boundary such that for any $x\in\Omega\setminus B_{r\varepsilon}(a)$
		\begin{align*}
			|\mathcal R(x)|&\leq C\delta^{\frac{N-2}{2}}\left[\frac{\varepsilon^{N-2}(1+\varepsilon\delta^{1-N})}{|x-a|^{N-2}}+\delta^2+\left(\frac{\varepsilon}{\delta}\right)^{N-2}\right]\\
			|\partial_{\tau_i}\mathcal R(x)|&\leq C\delta^{\frac{N}{2}}\left[\frac{\varepsilon^{N-2}(1+\varepsilon\delta^{-N})}{|x-a|^{N-2}}+\delta^2+\frac{\varepsilon^{N-2}}{\delta^{N-1}}\right]\\
			|\partial_{\delta}^{\,\tau}\mathcal R(x)|&\leq C\delta^{\frac{N-4}{2}}\left[\frac{\varepsilon^{N-2}(1+\varepsilon\delta^{1-N})}{|x-a|^{N-2}}+\delta^2+\left(\frac{\varepsilon}{\delta}\right)^{N-2}\right].
		\end{align*}
	\end{Lem}
	\begin{proof}
		The function $\mathcal R$ and its indicated parameter derivatives are harmonic.
		On $\partial\Omega$, expand the bubble in $\delta^2/|x-\xi|^2$;
		the first term is canceled by the regular part.
		On $\partial B_{r\varepsilon}(a)$, subtract $U_{\delta,\xi}(a)$ from
		$U_{\delta,\xi}(x)$ and apply Taylor's formula over a distance
		$r\varepsilon$.  These give precisely the constant boundary errors and the
		hole-boundary errors displayed above.  Comparison with a constant plus a
		multiple of $(r\varepsilon/|x-a|)^{N-2}$ proves the three bounds.
		For the derivatives the boundary remains fixed, since $\varepsilon$ is
		fixed when $\delta$ and $\tau$ are differentiated.
	\end{proof}

	\begin{Cor}
		\label{corA.3}
		Under the assumptions above, $\partial_{\xi_h}$ is taken with $\delta$ fixed,
		and $\partial_\delta$ with the center $\xi$ fixed, as throughout the paper; the
		scale derivative at fixed $\tau$ is written $\partial_\delta^{\,\tau}$. Then
		\begin{align*}
			\partial_{\xi_h}(P_{\varepsilon}U_{\delta,\xi}(x)-U_{\delta,\xi}(x))
			&=-\alpha_N\delta^{\frac{N-2}{2}}\partial_{\xi_h}H(x,\xi)
			+\frac{\alpha_N(N-2)\tau_h}{\delta^{\frac{N}{2}}(1+|\tau|^2)^{\frac{N}{2}}}
			\left(\frac{r\varepsilon}{|x-a|}\right)^{N-2}
			+\frac{1}{\delta}\partial_{\tau_h}\mathcal R(x),\\
			\partial_{\delta}\left(P_{\varepsilon}U_{\delta,\xi}(x)-U_{\delta,\xi}(x)\right)
			&=-\alpha_N\frac{N-2}{2}\delta^{\frac{N-4}{2}}H(x,\xi)\\
			&\quad+\frac{\alpha_N(N-2)}{2\delta^{\frac{N}{2}}}
			\frac{1-|\tau|^2}{(1+|\tau|^2)^{\frac N2}}
			\left(\frac{r\varepsilon}{|x-a|}\right)^{N-2}\\
			&\quad+\partial_\delta^{\,\tau}\mathcal R(x)-\frac1\delta\tau\cdot\nabla_\tau \mathcal R(x),
		\end{align*}
		and in particular there exists $C>0$ depending on the fixed domain, $r$, $\gamma_*$, and the distance of $a$ from the outer boundary such that
		\begin{align*}
			|P_{\varepsilon}U_{\delta,\xi}(x)-U_{\delta,\xi}(x)|&\leq C\delta^{\frac{N-2}{2}}+C\frac{\delta^{\frac{3}{2}N-3}}{|x-a|^{N-2}}\\
			|P_{\varepsilon}\psi_i^h(x)-\psi_i^h(x)|&\leq C\delta^{\frac{N-2}{2}}+C\frac{\delta^{\frac{3}{2}N-4}}{|x-a|^{N-2}}\\
			|P_{\varepsilon}\psi_i^0(x)-\psi_i^0(x)|&\leq C\delta^{\frac{N-4}{2}}+C\frac{\delta^{\frac{3}{2}N-4}}{|x-a|^{N-2}},
		\end{align*}
		for $h=1,\ldots,N$, for every $\varepsilon>0$ small enough, where
		$\psi_i^0=\partial_\delta U_{\delta,\xi}$ is the scale kernel at
		fixed center. In particular, if we fix a compact set $K\Subset\Omega$, then
		one can choose one constant $C>0$ uniformly for $a\in K$.
	\end{Cor}

	\begin{proof}
		The result follows directly from Lemma~\ref{lemA.2} and the boundedness of $H(x,\xi)$ and its derivatives for $x\in\Omega$ and $\xi\in K'\Subset\Omega$. We also use $\delta=d\sqrt\varepsilon$ with $d$ in a fixed compact subset of $(0,+\infty)$ and the fact that parameter differentiation commutes with the projection $P_{\varepsilon}$.
	\end{proof}
	
	\section{Basic Estimates in the Rescaled Variables}\label{secB}
	Here the projection is again taken on the full two-hole domain
	$\widetilde\Omega_\varepsilon$.  Put
	\[
	\widetilde H_\varepsilon(y,\widetilde\xi_i)
	=\varepsilon^{(N-2)/2}H(\sqrt\varepsilon y,\sqrt\varepsilon\widetilde\xi_i).
	\]
	The following coarse bounds suffice for every use of the expansion in the
	main text.  Constants are uniform on the fixed compact parameter range.
	
	\begin{Lem}\label{lemB.1}
		For $i\ne l$ define
		\[
		\begin{aligned}
			\mathcal R_i(y):={}&
			\widetilde P_\varepsilon U_{d_i,\widetilde\xi_i}(y)
			-U_{d_i,\widetilde\xi_i}(y)
			+\alpha_Nd_i^{(N-2)/2}\widetilde H_\varepsilon(y,\widetilde\xi_i)\\
			&+\frac{\alpha_N\varepsilon^{(N-2)/2}r_i^{N-2}}
			{d_i^{(N-2)/2}(1+|\tau_i|^2)^{(N-2)/2}|y-\widetilde a_i|^{N-2}}.
		\end{aligned}
		\]
		Here the superscript specifies the variable that is kept fixed: $\partial_{d_i}$
		keeps the center $\widetilde\xi_i$ fixed, as throughout the paper, whereas
		$\partial_{d_i}^{\,\tau}$ denotes the scale derivative at fixed $\tau_i$.  The
		latter is the derivative for which the coarse bound below is proved, because the
		boundary values of the remainder are differentiated in the coordinates
		$(d_i,\tau_i)$; the two are related by the chain rule
		\begin{equation}\label{eqB.1}
			\partial_{d_i}
			=\partial_{d_i}^{\,\tau}-d_i^{-1}\tau_i\cdot\nabla_{\tau_i}.
		\end{equation}
		If $\partial_\eta$ is any of $\partial_{\tau_{i,h}}$,
		$\partial_{d_i}^{\,\tau}$, or $\partial_{d_i}$, then
		\begin{equation}\label{eqB.2}
			|\mathcal R_i|+|\partial_\eta\mathcal R_i|
			\le C\left(
			\varepsilon^{N/2}
			+\frac{\varepsilon^{(N-1)/2}}{|y-\widetilde a_i|^{N-2}}
			+\frac{\varepsilon^{N-2}}{|y-\widetilde a_l|^{N-2}}\right)
			\quad\hbox{in }\widetilde\Omega_\varepsilon.
		\end{equation}
	\end{Lem}
	\begin{proof}
		The remainder is harmonic in $y$.  On the outer boundary, the normalization
		of $H$ cancels the leading far-field term of the bubble; the residual is
		$O(\varepsilon^{N/2}+\varepsilon^{N-2})$.
		On the $i$-th hole, the capacity term equals the bubble value at
		$\widetilde a_i$, so the residual is $O(\sqrt\varepsilon)$ by Taylor's
		formula.  On the other hole it is $O(\varepsilon^{(N-2)/2})$.
		The constant and the two harmonic functions
		$|y-\widetilde a_j|^{2-N}$, multiplied by the coefficients in
		\eqref{eqB.2}, dominate these three boundary errors.
		The maximum principle gives the claimed estimate.
		
		Parameter differentiation leaves the two-hole domain fixed.  The differentiated
		remainder is harmonic and has the same three coarse boundary bounds.  On the $i$-th hole one differentiates the Taylor difference
		$U(y)-U(\widetilde a_i)$, whose size is still
		$O(\sqrt\varepsilon)$.  The same barriers prove the derivative estimates.
		The chain rule gives the fixed-center formula.
	\end{proof}
	
	\begin{Cor}\label{corB.2}
		Put $\Pi_i:=U_{d_i,\widetilde\xi_i}-\widetilde P_\varepsilon
		U_{d_i,\widetilde\xi_i}$, the projection defect.  At fixed scale,
		respectively fixed center, one has
		\begin{align*}
			\partial_{\widetilde\xi_{i,h}}\Pi_i
			={}&\alpha_Nd_i^{(N-2)/2}
			\partial_{\widetilde\xi_{i,h}}\widetilde H_\varepsilon
			-\frac{\alpha_N(N-2)\tau_{i,h}\varepsilon^{(N-2)/2}r_i^{N-2}}
			{d_i^{N/2}(1+|\tau_i|^2)^{N/2}|y-\widetilde a_i|^{N-2}}
			-d_i^{-1}\partial_{\tau_{i,h}}\mathcal R_i,\\
			\partial_{d_i}\Pi_i
			={}&\alpha_N\frac{N-2}{2}d_i^{(N-4)/2}\widetilde H_\varepsilon
			-\frac{\alpha_N(N-2)(1-|\tau_i|^2)\varepsilon^{(N-2)/2}r_i^{N-2}}
			{2d_i^{N/2}(1+|\tau_i|^2)^{N/2}|y-\widetilde a_i|^{N-2}}
			-\partial_{d_i}\mathcal R_i.
		\end{align*}
		Moreover, on the full rescaled domain,
		\begin{align*}
			|\Pi_i|&\le C\varepsilon^{(N-2)/2}
			(1+|y-\widetilde a_i|^{2-N}),\\
			|\partial_{\widetilde\xi_{i,h}}\Pi_i|
			&\le C\left(\varepsilon^{(N-1)/2}
			+\varepsilon^{(N-2)/2}|y-\widetilde a_i|^{2-N}\right),\\
			|\partial_{d_i}\Pi_i|
			&\le C\varepsilon^{(N-2)/2}(1+|y-\widetilde a_i|^{2-N}).
		\end{align*}
	\end{Cor}
	\begin{proof}
		Differentiate the definition of $\mathcal R_i$ to obtain the
		identities.  For the coarse bounds, $\Pi_i$ is harmonic with boundary
		value $U_{d_i,\widetilde\xi_i}$.  This boundary value is $O(1)$ on its own small
		hole and $O(\varepsilon^{(N-2)/2})$ on the other boundary components.
		For a translation derivative the latter bound improves to
		$O(\varepsilon^{(N-1)/2})$.  Comparison with a constant plus
		$C\varepsilon^{(N-2)/2}|y-\widetilde a_i|^{2-N}$ gives the displayed
		bounds, including the second hole.  All scale factors $d_i$ are bounded
		above and below on the prescribed compact range.
	\end{proof}
	
	\begin{Lem}\label{lemB.3}
		Let $\Pi_i$ be the projection defect of Corollary~\ref{corB.2}.
		Then, uniformly on $\widetilde\Omega_\varepsilon$,
		\begin{align}
			|\partial_{\widetilde\xi_{i,h}}\Pi_i|
			&\le C\left(\varepsilon^{(N-1)/2}
			+\varepsilon^{(N-2)/2}|y-\widetilde a_i|^{2-N}\right),
			\label{eqB.3}\\
			|\partial_{d_i}\Pi_i|
			&\le C\varepsilon^{(N-2)/2}(1+|y-\widetilde a_i|^{2-N}),
			\label{eqB.4}\\
			|\partial_{\tau_{i,h}}\Pi_i|
			&\le C\left(\varepsilon^{(N-1)/2}
			+\varepsilon^{(N-2)/2}|y-\widetilde a_i|^{2-N}\right).
			\label{eqB.5}
		\end{align}
	\end{Lem}
	\begin{proof}
		The first two statements are Corollary~\ref{corB.2}.
		The last follows from
		$\partial_{\tau_{i,h}}=d_i\partial_{\widetilde\xi_{i,h}}$ at fixed $d_i$.
	\end{proof}
	
	\medskip
	\noindent\textbf{Auxiliary weighted quadratic estimate for the
		finite-dimensional projections.}  Here $\alpha$ and $R_\varepsilon$ are as in
	Section~\ref{sec7}, $\Pi_i$ is the projection defect of Corollary~\ref{corB.2},
	and $\phi_i$ denotes the inner correction.
	Let
	\[
	\bar q:=2^*-1=\frac{N+2}{N-2},
	\qquad f(t):=|t|^{\bar q-1}t,
	\]
	and, for $a>0$, define
	\[
	\mathcal N_a(h):=f(a+h)-f(a)-f'(a)h.
	\]
	For every scale or translation parameter $\eta$ used in the reduction,
	\[
	\bigl|\chi_i\,\partial_\eta\bigl(\widetilde P_\varepsilon U_{d_i,\widetilde\xi_i}\bigr)\bigr|
	\le C v_i^*.
	\]
	By Lemma~\ref{lem6.4},
	$\|\phi_i\|_{H_0^1}
	=O(\varepsilon^{(N-2\theta)/4}(\log\frac1\varepsilon)^{1/2})$.  Therefore
	\begin{equation}\label{eqB.quad}
		\int_{\widetilde\Omega_\varepsilon}
		|\mathcal N_{v_i^*}(\phi_i)|\,\bigl|\chi_i\,\partial_\eta\bigl(\widetilde P_\varepsilon U_{d_i,\widetilde\xi_i}\bigr)\bigr|
		=o\bigl(\varepsilon^{(N-2)/2}\bigr)
		\qquad\text{for every }N\ge5.
	\end{equation}
	Indeed, for $N\ge6$, the normalized quotient
	\[
	\frac{|1+x|^{\bar q-1}(1+x)-1-\bar qx}{x^2}
	\]
	is bounded on $\mathbb R$, with its continuous value at $x=0$. Hence
	\[
	|\mathcal N_a(h)|\le C a^{\bar q-2}|h|^2,
	\qquad a>0.
	\]
	For $N=5$, $\bar q=7/3$ and
	\[
	|\mathcal N_a(h)|
	\le C\bigl(a^{1/3}|h|^2+|h|^{7/3}\bigr).
	\]
	To estimate the parameter kernel, differentiate
	$-\Delta(\widetilde P_\varepsilon U_{d_i,\widetilde\xi_i})
	=U_{d_i,\widetilde\xi_i}^{2^*-1}$.
	Since the scale parameters remain in a fixed compact subset of
	$(0,+\infty)$, one has
	$|\partial_\eta U_{d_i,\widetilde\xi_i}|\le CU_{d_i,\widetilde\xi_i}$.
	Kato's inequality and the comparison principle then give
	\[
	|\partial_\eta(\widetilde P_\varepsilon U_{d_i,\widetilde\xi_i})|
	\le C\widetilde P_\varepsilon U_i,
	\]
	which is equivalent to the asserted kernel bound because $v_i^*$ differs
	from $\widetilde P_\varepsilon U_{d_i,\widetilde\xi_i}$ only by a fixed positive
	factor. For $N\ge6$ it follows that
	\[
	|\mathcal N_{v_i^*}(\phi_i)|\,\bigl|\chi_i\,\partial_\eta\bigl(\widetilde P_\varepsilon U_{d_i,\widetilde\xi_i}\bigr)\bigr|
	\le C(v_i^*)^{2^*-2}|\phi_i|^2.
	\]
	Therefore, by H\"older's and Sobolev inequalities,
	\[
	\int_{\widetilde\Omega_\varepsilon}
	(v_i^*)^{2^*-2}|\phi_i|^2
	\le C\|v_i^*\|_{L^{2^*}}^{2^*-2}
	\|\phi_i\|_{L^{2^*}}^2
	\le C\|\nabla\phi_i\|_2^2
	=O\bigl(\varepsilon^{(N-2\theta)/2}\bigr).
	\]
	For $N=5$, the same argument gives
	\[
	\int_{\widetilde\Omega_\varepsilon}
	|\mathcal N_{v_i^*}(\phi_i)|\,\bigl|\chi_i\,\partial_\eta\bigl(\widetilde P_\varepsilon U_{d_i,\widetilde\xi_i}\bigr)\bigr|
	\le C\left(\|\phi_i\|_{L^{2^*}}^2
	+\|\phi_i\|_{L^{2^*}}^{7/3}\right),
	\]
	and both terms are $o(\varepsilon^{3/2})$ when $0<\theta<1$.
	
	For the projection defect $\Pi_i$ of Corollary~\ref{corB.2}, the
	maximum principle applied to the differentiated
	harmonic boundary problem gives $|\partial_\eta\Pi_i|\le C\Pi_i$.
	Moreover,
	\begin{equation}\label{eqB.proj}
		\int\chi_iU_{d_i,\widetilde\xi_i}^{2^*-2}\Pi_i^2
		\le C\left[\varepsilon^{N-2}(R_\varepsilon)^{N-4}
		+\varepsilon^{N/2}\right].
	\end{equation}
	Insert $\Pi_i\le C\varepsilon^{(N-2)/2}(1+|y-\widetilde a_i|^{2-N})$ from
	Appendix~\ref{secB} into the integral: the regular part contributes
	$C\varepsilon^{N-2}(R_\varepsilon)^{N-4}$, the neighbourhood of the own hole
	contributes at most $C\varepsilon^{N-2}\int_{r_i\sqrt\varepsilon}^{C}r^{3-N}\,dr
	\le C\varepsilon^{N/2}$, and on the rest the bubble weight makes the integrand
	integrable, while the other hole lies outside $\operatorname{supp}\chi_i$.
	Both terms in \eqref{eqB.proj} are $o(\varepsilon^{(N-2)/2})$ because $\alpha<1/2$.
	
	\section{Continuous dependence of the solution on parameters}\label{secC}
	Fix the already chosen $R_0,R$ and a sufficiently small $\varepsilon>0$.
	Only $\mathbf p=(\mathbf d,\boldsymbol\tau)$ varies.  Constants in this
	appendix may depend on these fixed quantities.  All estimates are taken on
	a slightly larger compact parameter box containing
	$\mathcal P_{\gamma_*}$ in its interior.
	
	We first specify how functions on the moving regions are compared.
	Near a reference parameter $\widehat{\mathbf p}$ choose, for each component,
	\[
	\Theta_{i,\mathbf p}(y)=y+\rho_i(y)
	(\widetilde\xi_i(\mathbf p)-\widehat{\widetilde\xi}_i).
	\]
	The fixed smooth function $\rho_i$ is one near both outer spheres of
	$\widehat M_0,\widehat M_3$ (or $\widehat Q_0,\widehat Q_3$), zero near the small holes $\partial B_{r_i\sqrt\varepsilon}(\widetilde a_i)$
	and near $\partial\widetilde\Omega_\varepsilon$, and supported in $B_{4R}(\widetilde a_i)$.
	The strict separations fixed earlier allow this choice.  For nearby parameters
	these maps are diffeomorphisms of the full domain, map the indicated inner
	regions to their current counterparts, and tend to the identity in $C^2$.
	They also preserve the corresponding auxiliary outer domain, whose artificial
	hole is at the opposite center.  Near each matching sphere the map is just
	a translation.  Composition with these maps gives equivalent $H^1$, $L^q$,
	and trace norms.  Below, convergence on moving regions means convergence
	after this composition.
	
	For the first component write $D=\widehat M_3$ and $\Theta=\Theta_{1,\mathbf p}$.
	The coefficients and moments after this change of variables are
	\[
	\begin{aligned}
		J_{\mathbf p}&=\det \nabla\Theta,
		&A_{\mathbf p}&=J_{\mathbf p}(\nabla\Theta)^{-1}(\nabla\Theta)^{-T},\\
		W_{\mathbf p}&=J_{\mathbf p}(V_{1,\mathbf p}\circ\Theta),
		&B_{\mathbf p}^h&=J_{\mathbf p}(Z_{1,\mathbf p}^h\circ\Theta),\\
		m_{\mathbf p}^k(v)&=\int_D A_{\mathbf p}\nabla v\cdot
		\nabla(Z_{1,\mathbf p}^k\circ\Theta).
	\end{aligned}
	\]
	Thus the source is $F=J_{\mathbf p}(f\circ\Theta)$, and the equation is
	\[
	\mathcal L_{\mathbf p}v:=-\operatorname{div}(A_{\mathbf p}\nabla v)
	-W_{\mathbf p}v
	=F+\sum_{h=0}^Nc_hB_{\mathbf p}^h,
	\qquad v\in H_0^1(D),\quad m_{\mathbf p}^k(v)=0.
	\]
	Here $A_{\mathbf p}\to I$ in $C^1$ and $W_{\mathbf p}\to\widehat W$
	in $L^\infty$.  For fixed $\varepsilon$, the functions $Z^h_{1,\mathbf p}$ solve smooth
	Dirichlet problems on the fixed full domain.  Consequently their compositions
	converge in $H^1\cap L^q$, their images under $\mathcal L_{\mathbf p}$
	are bounded in $L^2$, and the displayed moments converge as bounded linear
	functionals on $H_0^1(D)$.
	
	\begin{Lem}\label{lemC.1}
		If $\mathbf p_n\to\widehat{\mathbf p}$ and $F_n\to\widehat F$ in
		$L^2(D)$, the corresponding zero-moment linear solutions satisfy
		\[
		\|v_n-\widehat v\|_{\widehat X_1}
		+|c_n-\widehat c|\longrightarrow0,
		\qquad
		\|v\|_{\widehat X_1}:=\|\nabla v\|_2+R^{-1}\|v\|_q,
		\quad q=\frac{2N}{N-4}.
		\]
		The same conclusion holds for the second component.
	\end{Lem}
	\begin{proof}
		Subtract a small combination of the functions $Z^h_{1,\mathbf p_n}\circ\Theta_n$
		from the reference solution,
		\[
		z_n=\widehat v-\sum_{h=0}^Na_{n,h}(Z_{1,\mathbf p_n}^h\circ\Theta_n),
		\qquad m_{\mathbf p_n}^k(z_n)=0.
		\]
		The coefficients solve the finite positive-definite matrix equation formed
		by the moments of these functions.  Since the reference moments are zero and
		the matrices converge, $a_n\to0$.  Hence $z_n\to\widehat v$ in
		$\widehat X_1$ and
		$\|\mathcal L_{\mathbf p_n}(\widehat v-z_n)\|_2\to0$.
		
		Now $w_n=v_n-z_n$ has the current zero moments.  Direct subtraction gives
		\begin{equation}\label{eqC.1}
			\begin{aligned}
				\mathcal L_{\mathbf p_n}w_n
				&=r_n+\sum_h(c_{n,h}-\widehat c_h)B_{\mathbf p_n}^h,\\
				r_n&=F_n-\widehat F+
				(\mathcal L_{\widehat{\mathbf p}}-\mathcal L_{\mathbf p_n})\widehat v
				+\mathcal L_{\mathbf p_n}(\widehat v-z_n)
				+\sum_h\widehat c_h(B_{\mathbf p_n}^h-\widehat B^h).
			\end{aligned}
		\end{equation}
		All source terms tend to zero in $L^2$.  For the principal-part term this
		uses just one interior estimate: $A_{\mathbf p_n}-I$ is supported on a
		fixed compact transition region $K\Subset D$, and the reference equation gives
		$\widehat v\in W^{2,2}$ near $K$.  Thus
		\[
		\|\operatorname{div}((A_{\mathbf p_n}-I)\nabla\widehat v)\|_2
		\le C\|A_{\mathbf p_n}-I\|_{C^1}\|\widehat v\|_{W^{2,2}(K)}\to0.
		\]
		Pushing \eqref{eqC.1} to the physical inner domain
		returns precisely the projected problem of Lemma~\ref{lem5.1};
		the physical source is $(r_n/J_{\mathbf p_n})\circ\Theta_n^{-1}$.
		That lemma therefore yields
		\[
		\|w_n\|_{\widehat X_1}+|c_n-\widehat c|\le C_R\|r_n\|_2\to0.
		\]
		Adding $z_n-\widehat v$ proves the assertion.
	\end{proof}
	
	\begin{Lem}\label{lemC.2}
		If the prescribed data converge in the local energy and boundary norms of
		$\Lambda_\varepsilon$ and in the $L^q$ norms defining $\mathrm{dist}$
		(Section~\ref{sec6}), the inner solutions converge in $H^1\cap L^q$ and the outer solutions
		converge in $H^1$, after pulling back by the coordinate maps $\Theta_{i,\mathbf p}$.  In particular, the matched
		solution, its multipliers, and all reduced residuals are continuous in
		$\mathbf p\in \mathcal P_{\gamma_*}$.
	\end{Lem}
	\begin{proof}
		Let $\mathbf p_n\to\widehat{\mathbf p}$.
		
		\textbf{Step 1: Inner solutions.}
		Let $a_n$ be the pulled-back harmonic lifting of boundary data converging in
		$H^{1/2}(\partial D)\cap L^\infty(\partial D)$.
		Its equation is $-\operatorname{div}(A_n\nabla a_n)=0$.
		Choose a lifting $e_n\to0$ in $H^1(D)$ of the boundary difference.
		Testing the difference by $a_n-\widehat a-e_n$ gives
		\[
		\|a_n-\widehat a\|_{H^1}
		\le C\bigl(\|e_n\|_{H^1}+\|A_n-I\|_\infty\|\nabla\widehat a\|_2\bigr)
		\longrightarrow0.
		\]
		The maximum principle bounds $a_n$ in $L^\infty$, so interpolation gives
		convergence in $L^q$ as well.  Also $m_n^k(a_n)=0$, by testing its weak
		equation with $Z^k_{1,\mathbf p_n}\circ\Theta_n$.
		
		For the homogeneous inner part the fixed-point map is defined by the
		projected inverse with source
		\[
		J_n\left[\widetilde l_{1,\mathbf p_n}\circ\Theta_n
		+\mathcal N_{1,\mathbf p_n}^{\Theta_n}(a_n+w)
		+(V_{1,\mathbf p_n}\circ\Theta_n)a_n\right].
		\]
		The superscript means that the background in the pointwise remainder is
		composed with $\Theta_n$.  The map can be evaluated for every $w$ in a
		small ball of $\widehat X_1$; its output, not its input, is required to have
		the current zero moments.  At a fixed $w$ the source displayed above converges in $L^2$,
		by the previous estimates and the nonlinear inequalities of Section~\ref{sec5}.
		Lemma~\ref{lemC.1} gives convergence of the maps at that input.
		The original contraction estimate, with a slightly smaller choice of
		$\varepsilon$, gives a common $\kappa<1$ on a ball containing both fixed points.
		Consequently,
		\[
		\|w_n-\widehat w\|_{\widehat X_1}
		\le\kappa\|w_n-\widehat w\|_{\widehat X_1}
		+\|\mathcal T_{1,\mathbf p_n}(\widehat w)
		-\mathcal T_{1,\widehat{\mathbf p}}(\widehat w)\|_{\widehat X_1}
		=\kappa\|w_n-\widehat w\|_{\widehat X_1}+o(1).
		\]
		Absorption proves inner continuity, including the multipliers.
		On smaller collars of $\partial B_{R_0}(\widehat{\widetilde\xi}_1)$,
		Lemma~\ref{lem5.10} first gives a uniform bound.  The scalar equation
		then gives uniform $C^{1,\alpha}$ bounds.  Compactness and the identified
		$L^q$ limit give $C^1$ convergence there, in particular uniform convergence
		of the traces.  The small-hole traces are zero.
		
		\textbf{Step 2: Outer solutions.}
		Pull each component back to the reference outer domain
		$\widetilde\Omega_\varepsilon\setminus(\widehat M_0\cup\widehat Q_0)$, and
		call the full components there $u_{i,n}$.  The prescribed local data converge in $H^1$
		and uniformly on their boundaries.  Bounded $H^1$ extensions give bounded
		energy competitors, and the minimizing property gives a uniform $H^1$ bound
		for $u_{i,n}$.  The exterior pointwise bounds give a uniform bound on that
		domain.  A subsequence therefore converges weakly in $H^1$, strongly in
		$L^2$, and almost everywhere.
		
		The coupling is zero on the supports of the coordinate changes, by
		Lemma~\ref{lem4.6}; where it can be nonzero both maps are the identity.
		Thus no different coordinates are mixed in the coupling term.  The pointwise
		bounds and almost-everywhere convergence imply $L^2$ convergence of the
		pulled-back nonlinear sources, denoted by $F_{i,n}$.  The weak limit has the
		correct boundary traces and solves the reference outer problem.  Its bounds
		pass to the limit, so outer uniqueness identifies it with $\widehat u_i$.
		
		For strong convergence, lift the difference of the boundary data to
		$e_{i,n}\to0$ in $H^1$ of the same domain and set
		$z_{i,n}=u_{i,n}-\widehat u_i-e_{i,n}\in H_0^1$.
		Testing the difference of the two weak equations by $z_{i,n}$ gives
		\[
		\|\nabla z_{i,n}\|_2
		\le C_{\varepsilon,R}\left(
		\|F_{i,n}-\widehat F_i\|_2
		+\|A_{i,n}-I\|_\infty\|\nabla\widehat u_i\|_2
		+\|e_{i,n}\|_{H^1}\right)\longrightarrow0.
		\]
		Together with the prescribed-data convergence this proves global $H^1$
		convergence of the outer functions.  On the matching collars their scalar
		equations give $C^1$ convergence as above, and hence convergence of the
		boundary traces in $H^{1/2}\cap L^\infty$.
		
		\textbf{Step 3: Matched fixed points.}
		Transport the prescribed data of the reference fixed point $\widehat\Phi$
		to the current regions.  Subtract a small combination of the functions $Z^h_{1,\mathbf p_n}\circ\Theta_n$
		to restore the zero moments, exactly as in Lemma~\ref{lemC.1}.
		These functions vanish on the spheres $\partial B_{R_0}(\widetilde a_i)$ and on the
		small holes $\partial B_{r_i\sqrt\varepsilon}(\widetilde a_i)$.
		If necessary multiply by
		\[
		t_n=\min\{1,\varepsilon^{(N-2\theta)/4}/B_n\},
		\]
		where $B_n$ is the largest of the four local norms constrained in
		$\Lambda_\varepsilon$; take $t_n=1$ when $B_n=0$.
		Convergence of the corrected data implies $\limsup B_n\le
		\varepsilon^{(N-2\theta)/4}$, so $t_n\to1$.
		We have therefore constructed $H_n\in\Lambda_\varepsilon(\mathbf p_n)$
		with the required local convergence.
		
		Set $E_n=S_{\mathbf p_n}(H_n)$ and let $\Phi_n$ be the current matched
		fixed point.  Steps 1--2 and preservation of prescribed data by $S$ give
		\[
		\mathrm{dist}_{\mathbf p_n}(S_{\mathbf p_n}(T_{\mathbf p_n}(E_n)),E_n)\to0.
		\]
		Both comparison points now belong to the same space.  Thus
		\[
		\mathrm{dist}_{\mathbf p_n}(\Phi_n,E_n)
		\le\tfrac12\mathrm{dist}_{\mathbf p_n}(\Phi_n,E_n)
		+\mathrm{dist}_{\mathbf p_n}(S_{\mathbf p_n}(T_{\mathbf p_n}(E_n)),E_n),
		\]
		and absorption proves convergence in the matching metric.
		Trace--Poincar\'e gives local $H^1$ convergence; applying Step 2 to the
		prescribed data of these fixed points gives global $H_0^1$ convergence.
		The same global $H_0^1$ convergence also holds without composition on the fixed full physical domain, by
		continuity of composition at the limiting $H^1$ function.
		
		\textbf{Step 4: Reduced residuals.}
		Step 1 gives continuity of all multipliers.  The matrices relating them to
		the reduced residuals are integrals of the smoothly varying functions $Z^h_{1,\mathbf p}$
		and of their derivatives on the fixed full domain.  Their entries are continuous,
		so every reduced residual is continuous.  This includes the boundary of the
		parameter box needed in the degree argument.
	\end{proof}
	
	\noindent\textbf{Acknowledgments}:
	This research was supported by the National Natural Science Foundation of China (Grant Nos. 12271539, 12371107).
	
	\noindent\textbf{Data availability}:
	No data was used for the research described in the article.

\end{document}